\documentclass{amsart}
\usepackage{amsmath,amssymb,amsthm,geometry,graphics,color,mdwlist,mathrsfs,bm}
\usepackage[all]{xy}


\usepackage{url}
\usepackage{xcolor}
\usepackage{fancybox}
\usepackage{ascmac}
\usepackage{amsfonts}
\usepackage{mathrsfs}
\usepackage[mathscr]{eucal}
\usepackage{graphics}
\usepackage{array,booktabs,float}
\usepackage[all]{xy}
\usepackage{tikz}
\usetikzlibrary{cd,bending,arrows.meta,positioning}
\usepackage{fancyhdr}
\newtheorem{thm}{Theorem}[section]
\newtheorem{prop}[thm]{Proposition}
\newtheorem{lem}[thm]{Lemma}
\newtheorem{cor}[thm]{Corollary}
\newtheorem{conj}[thm]{Conjecture}
\theoremstyle{definition}
\newtheorem{dfn}[thm]{Definition}
\newtheorem{eg}[thm]{Example}
\newtheorem{rem}[thm]{Remark}
\newtheorem{setup}[thm]{Setup}

\newtheorem{ques}[thm]{Question}
\newtheorem{nota}[thm]{Notation}
\makeatother

\usepackage{enumerate}
\usepackage[utf8]{inputenc}
\usepackage{enumitem}

\usepackage{hyperref}
\hypersetup{%
 setpagesize=false,%
 bookmarks=true,%
 bookmarksdepth=tocdepth,%
 bookmarksnumbered=true,%
 colorlinks=true,%
 linkcolor=teal,
 citecolor=orange,
 pdftitle={},%
 pdfsubject={},%
 pdfauthor={},%
 pdfkeywords={}}

\newcommand{\cH}{\mathcal{H}}

\newcommand{\bC}{\mathbb{C}}

\newcommand{\bG}{\mathbb{G}}
\newcommand{\bH}{\mathbb{H}}

\newcommand{\bR}{\mathbb{R}}

\newcommand{\bZ}{\mathbb{Z}}

\DeclareMathOperator{\ZZ}{\mathbb{Z}}
\newcommand{\bLam}{\bm{\Lambda}}
\newcommand{\bsL}{\bm{\Lambda}^{\st}}

\newcommand{\m}{\mathfrak{m}}

\newcommand{\p}{\mathfrak{p}}

\newcommand{\scrA}{\EuScript{A}}

\newcommand{\scrD}{\EuScript{D}}

\newcommand{\scrF}{\EuScript{F}}

\newcommand{\scrH}{\EuScript{H}}

\newcommand{\scrJ}{\EuScript{J}}
\newcommand{\scrK}{\EuScript{K}}

\newcommand{\scrT}{\EuScript{T}}
\newcommand{\scrU}{\EuScript{U}}
\newcommand{\scrV}{\EuScript{V}}

\newcommand{\scrX}{\EuScript{X}}
\newcommand{\scrY}{\EuScript{Y}}

\DeclareMathOperator{\oA}{\overline{\mathit{A}}}

\DeclareMathOperator{\oL}{\overline{\mathit{L}}}
\DeclareMathOperator{\oM}{\overline{\mathit{M}}}
\DeclareMathOperator{\oN}{\overline{\mathit{N}}}
\DeclareMathOperator{\oP}{\overline{\mathit{P}}}
\DeclareMathOperator{\oR}{\overline{\mathit{R}}}
\DeclareMathOperator{\oT}{\overline{\mathit{T}}}
\DeclareMathOperator{\oS}{\overline{\mathit{S}}}

\DeclareMathOperator{\oI}{\overline{\mathit{I}}}

\DeclareMathOperator{\oLam}{\overline{\Lambda}}

\DeclareMathOperator{\Ext}{Ext}
\DeclareMathOperator{\CM}{CM}

\DeclareMathOperator{\Ker}{Ker}
\DeclareMathOperator{\proj}{proj}
\DeclareMathOperator{\Hom}{Hom}

\DeclareMathOperator{\modcat}{mod}
\DeclareMathOperator{\depth}{depth}

\DeclareMathOperator{\rad}{rad}

\DeclareMathOperator{\pd}{pd}
\DeclareMathOperator{\id}{id}

\DeclareMathOperator{\Spec}{Spec}

\DeclareMathOperator{\add}{add}

\DeclareMathOperator{\End}{End}
\DeclareMathOperator{\fl}{fl}

\DeclareMathOperator{\height}{ht}
\DeclareMathOperator{\pp}{\mathfrak{p}}

\DeclareMathOperator{\mm}{\mathfrak{m}}
\DeclareMathOperator{\nn}{\mathfrak{n}}

\DeclareMathOperator{\stHom}{\underline{\Hom}}
\DeclareMathOperator{\thick}{thick}

\DeclareMathOperator{\modif}{mdf}
\DeclareMathOperator{\modifG}{mdfG}
\DeclareMathOperator{\MM}{MM}
\DeclareMathOperator{\NCCR}{NCCR}
\DeclareMathOperator{\MMG}{MMG}
\DeclareMathOperator{\NCCRG}{NCCRG}
\DeclareMathOperator{\pvd}{pvd}
\DeclareMathOperator{\per}{per}

\DeclareMathOperator{\stEnd}{\underline{\End}}
\DeclareMathOperator{\stCM}{\underline{\CM}}

\DeclareMathOperator{\Db}{D^b}
\DeclareMathOperator{\Kb}{K^b}
\DeclareMathOperator{\refl}{ref}

\DeclareMathOperator{\D}{D}
\DeclareMathOperator{\Modcat}{Mod}
\DeclareMathOperator{\Cok}{Cok}
\DeclareMathOperator{\EG}{EG}
\DeclareMathOperator{\reg}{reg}

\DeclareMathOperator{\twosilt}{2-silt}

\DeclareMathOperator{\twopresilt}{2-presilt}
\DeclareMathOperator{\twotilt}{2-tilt}
\DeclareMathOperator{\RHom}{\mathbf{R}Hom}

\DeclareMathOperator{\Lo}{\overset{\mathbf{L}}{\otimes}}
\DeclareMathOperator{\Stab}{Stab}

\DeclareMathOperator{\sg}{sg}

\DeclareMathOperator{\st}{st}

\DeclareMathOperator{\std}{MM}
\DeclareMathOperator{\stdStab}{\Stab^{\std}}
\DeclareMathOperator{\Cone}{Cone}
\DeclareMathOperator{\Free}{Free}

\DeclareMathOperator{\tors}{tors}
\DeclareMathOperator{\ftors}{f-tors}
\DeclareMathOperator{\msilt}{m-silt}
\DeclareMathOperator{\Fac}{Fac}
\DeclareMathOperator{\silt}{silt}
\DeclareMathOperator{\alg}{alg}

\DeclareMathOperator{\Filt}{Filt}
\DeclareMathOperator{\topmod}{top}

\def\gldim{\mathop{\rm gl.dim}\nolimits}

\def\pd{\mathop{\rm proj.dim}\nolimits}

\begin{document}

\title[NCCRs and stability conditions via reduction modulo regular sequences]{Non-commutative crepant resolutions and stability conditions via reduction modulo regular sequences}

\author{Issei Enomoto}
\email{isseienomoto.math@gmail.com}
\maketitle
\begin{abstract}
In this paper, we investigate the behavior of NCCRs, MMAs, and spaces of stability conditions under reduction modulo regular sequences. Let $R$ be a complete Gorenstein normal local domain, and let $t_1,\dots,t_m$ be an $R$-regular sequence. We set $\overline{R}:=R/(t_1,\ldots,t_m)$ and $\overline{(-)}:=(-)\otimes_R\overline{R}$, and assume that $\overline{R}$ is a normal domain of dimension at least two.

First, we prove that, for a Cohen--Macaulay modifying $R$-module $M$, if $\overline{M}$ gives an NCCR of $\overline{R}$ (respectively, if $\overline{M}$ is a maximal modifying $\overline{R}$-module), then $M$ gives an NCCR of $R$ (respectively, $M$ is a maximal modifying $R$-module). These results can be regarded as a noncommutative analogue of the fact that a Noetherian local ring is regular if its quotient by a regular element is regular.

Next, suppose that $\overline{R}$ is a three-dimensional isolated singularity and that there exists a Cohen--Macaulay maximal modifying $R$-module $M$ such that $\overline{M}$ is a maximal modifying $\overline{R}$-module. In this setting, we establish a reduction theory from higher-dimensional singularities $R$ to three-dimensional singularities $\overline{R}$. We first establish a natural bijection between maximal modifying $R$-modules and maximal modifying $\overline{R}$-modules that preserves mutation. Furthermore, using a construction similar to that of Iyama--Wemyss, we construct tilting modules over maximal modifying algebras of $R$. As an application, we obtain an Adachi--Iyama--Reiten-type correspondence between Cohen--Macaulay modifying $R$-modules and two-term silting complexes over $\underline{\operatorname{End}}_R(M)$. Finally, we prove that suitable subspaces of the spaces of Bridgeland stability conditions associated with $\operatorname{End}_R(M)$ and $\operatorname{End}_{\overline{R}}(\overline{M})$ are isomorphic.
\end{abstract}

\setcounter{tocdepth}{1}
\tableofcontents

\section{Introduction}
\subsection{Background}
A {\it non-commutative crepant resolution (NCCR)} is a non-commutative algebra introduced by Van den Bergh as a non-commutative analogue of a crepant resolution of a singularity. NCCRs are known to be related not only to singularity theory but also to various other areas, such as higher Auslander--Reiten theory \cite{Higher-dimensional-Auslander-Reiten-theory-on-maximal-orthogonal-subcategories,Cohen-Macaulay-representations-of-invariant-subrings}, additive categorification of cluster algebras \cite{Tilting-theory-and-cluster-combinatorics,Mutation-in-triangulated-categories-and-rigid-Cohen-Macaulay-modules} and Calabi--Yau algebras \cite{Fomin-Zelevinsky-mutation-and-tilting-modules-over-Calabi-Yau-algebras.,Silting-correspondences-and-Calabi-Yau-dg-algebras}. Moreover, a {\it maximal modification algebra (MMA)} is a generalization of an NCCR introduced by Iyama--Wemyss and corresponds, on the geometric side, to a minimal model
\cite{Maximal-modifications-and-Auslander-Reiten-duality-for-non-isolated-singularities}. In Wemyss's Homological Minimal Model Program (HomMMP), a bijection between basic Cohen--Macaulay maximal modifying generators over complete local cDV singularities and minimal models was established
\cite{Flops-and-clusters-in-the-homological-minimal-model-programme}. As this indicates, NCCRs and MMAs provide a framework for algebraically describing the geometric structures appearing in crepant resolutions and minimal models. Therefore, the existence and classification of NCCRs and MMAs for a given singularity are fundamental problems that have been extensively studied.

On the other hand, Bridgeland introduced {\it stability conditions} on triangulated categories as a generalization of stability for vector bundles and King's stability for quiver representations
\cite{Stability-conditions-on-triangulated-categories}. Under suitable finiteness assumptions, the space of stability conditions carries the structure of a complex manifold and provides an important invariant of a triangulated category. Spaces of stability conditions on Calabi--Yau triangulated categories exhibit interesting structures. For example, the spaces of stability conditions on Calabi--Yau triangulated categories arising from NCCRs and MMAs associated with three-dimensional flops and cDV singularities are known to be covering spaces of the complements of complexified "contracted Coxeter arrangements" \cite{Stability-conditions-for-3-fold-flops}.

A rich theory has been developed for NCCRs, MMAs, and spaces of stability conditions associated with three-dimensional cDV singularities. One of the characteristic features of cDV singularities is their close relationship with Du Val singularities through taking quotients by regular elements. This naturally leads to the question of whether the theory developed in dimension three can be extended to higher-dimensional singularities whose quotients by regular sequences are three-dimensional isolated singularities. Motivated by these observations, we study how NCCRs, MMAs, and their associated spaces of stability conditions change when passing to quotients by regular sequences. We first show that, for a Cohen–Macaulay modifying module, the NCCR and maximal modifying properties pass from its quotient by a regular element to the original module. We then show that, under suitable assumptions, maximal modifying modules over the original ring and its quotient are in bijection, compatibly with mutation, and construct the corresponding tilting modules. We also compare the spaces of stability conditions on the associated derived categories over the original ring and its quotient.

\subsection{Ascent for NCCRs and MMAs}
In this subsection, we use the following notation:
\begin{setup}\label{Intro setup ascent NCCR}
    \begin{enumerate}
        \item Let $(R,\mm)$ be a $d$-dimensional complete Gorenstein local normal domain. 
        \item Let $t\in \mm$ be a regular element such that  $\overline{R}:=R/tR$ is a normal domain of dimension $\ge2$.
    \end{enumerate}
\end{setup}
\begin{nota}
    We put 
        \[\overline{(-)}:=-\underset{R}{\otimes} \overline{R}:\modcat R\to \modcat \oR\quad\text{and}\] \[(-)^*:=\Hom_{\oR}(-,\oR):\modcat \oR\to \modcat \oR.\]

\end{nota}
We first establish ascent-type results for NCCRs and MM modules. These results can be regarded as a noncommutative analogue of the fact that a Noetherian local ring is regular if its quotient by a regular element is regular.
\begin{thm}[= Theorem \ref{NCCR +1}]\label{1.2}
    Let $M\in\refl R$ be a modifying module. The following statements hold.
    \begin{enumerate}
        \item If $\overline{M}^{**}$ gives an NCCR of $\overline{R}$, then $M$ gives an NCCR of $R$.
        \item If $\overline{M}^{**}$ is an MM $\overline{R}$-module, then $M$ is an MM $R$-module.
    \end{enumerate}
\end{thm}

If $M$ is Cohen--Macaulay, then $\oM\cong \oM^{**}$. Hence we obtain the following.

\begin{cor}[= Corollary \ref{CM NCCR +1}]\label{intro NCCR +1}
  Let $M\in\CM R$ be a modifying module. The following statements hold.
    \begin{enumerate}
        \item If $\overline{M}$ gives an NCCR of $\overline{R}$, then $M$ gives an NCCR of $R$.
        \item If $\overline{M}$ is an MM $\overline{R}$-module, then $M$ is an MM $R$-module.
    \end{enumerate}
\end{cor}

We expect that these theorems will also provide new methods for constructing NCCRs and MMAs over higher-dimensional singularities. We present concrete examples in subsection \ref{ex subsection}.

% As an application, we then study NCCRs of $cD_4$ singularities. We begin by proving the following result concerning the defining equations of $cD_4$ singularities.

% \begin{lem}[= \ref{cD4 form}]
%      Let $k$ be an algebraically closed field of $\ch k\ne 2$, and let $(R,\m,k)$ be a complete local $cD_4$ singularity. Then there exists $f\in k[[x,y,t]]$ such that
%     \begin{enumerate}
%     \item $f|_{t=0}=x^2y+y^3$ and 
%     \item $R\cong k[[x,y,z,t]]/(z^2+f)$.
%     \end{enumerate}
% \end{lem}
% In this case, $f$ factors into at most three irreducible elements. In what follows, we consider the case where $f$ factors into three irreducible elements.
% \begin{setup}
%     Throughout this subsection, in addition to Setup\ref{Intro setup ascent NCCR}, we make the following assumptions:
%     \begin{enumerate}
%     \item Let $k$ be an algebraically closed field of $\ch k\ne 2$. 
%     \item Let  $(R,\m,k)$ be a complete local $cD_4$ singularity isomorphic to $k[[x,y,z,t]]/(z^2+f)$ for $f\in k[[x,y,t]]$. 
%     \item We assume that $f$ has an irreducible factorization $f=f_1f_2f_3$ in $k[[x,y,t]]$. 
%     \item Let $C:=k[[x,y,t]]/(f_1,f_2,f_3)$.
%     \end{enumerate}
% \end{setup}
% \begin{thm}[= \ref{cD_4 NCCR}, \ref{cD4 NCCR 2}, \ref{cD4 NCCR summary}]
%     $R$ has an NCCR if and only if $C\cong k[[t]]$ or $C\cong k$. Moreover, in this case, the structure of $\NCCRG(R)$ and the quivers of all NCCRs of $R$ can be described.
% \end{thm}
\subsection{Reduction for NCCRs, tilting modules and 2-term silting complexes}
In this subsection, we use the following notation:
\begin{setup}\label{intro setup  reduction NCCR,tilt,2-silt}
    \begin{enumerate}
        \item Let $(R,\mm)$ be a $d$-dimensional complete Gorenstein local normal domain. 
        \item Let $t_1,\dots,t_m\in \mm$ be a regular sequence of $R$ for $m\in \ZZ_{\ge 0}$. We assume that \[\overline{R}:=R/(t_1,\dots,t_m)\] is a $3$-dimensional Gorenstein isolated singularity. 
        \item Let $M\in\MM(R)$ be an MM $R$-module such that $\overline{M}^{**}\in\MM(\oR)$.
    \end{enumerate}
\end{setup}

We prove that there is a natural bijection between maximal modifyings.

\begin{thm}[= Theorem \ref{MM(R)->MM(R/t) MM}]\label{intro MM(R)->MM(R/t) MM}
   There exist the following bijections, compatible with the Iyama--Wemyss mutation:
\[\begin{tikzcd}
	{\MM(R)} && {\MM(\oR)} & {\modif(R)} && {\modif(\oR)} \\
	{\MMG(R)} && {\MMG(\oR)} & {\modifG(R)} && {\modifG(\oR)}
	\arrow["{\overline{(-)}^{**}}"', from=1-1, to=1-3]
	\arrow["\sim", from=1-1, to=1-3]
	\arrow["{\overline{(-)}^{**}}"', from=1-4, to=1-6]
	\arrow["\sim", from=1-4, to=1-6]
	\arrow[hook', from=2-1, to=1-1]
	\arrow["{\overline{(-)}}"', from=2-1, to=2-3]
	\arrow["\sim", from=2-1, to=2-3]
	\arrow[hook', from=2-3, to=1-3]
	\arrow[hook', from=2-4, to=1-4]
	\arrow["{\overline{(-)}}"', from=2-4, to=2-6]
	\arrow[from=2-4, to=2-6]
	\arrow["\sim", from=2-4, to=2-6]
	\arrow[hook', from=2-6, to=1-6]
\end{tikzcd}\]

For $\MM$,$\MMG$,$\modif$,$\modifG$, see Definition\ref{2.14}.
\end{thm}

This is a particularly interesting result. Indeed, an analogue of the above theorem fails completely when a cDV singularity is reduced to a Du Val singularity by a regular element. In other words, while there are substantial differences in representation theory between dimensions two and three, even under favorable assumptions, this result suggests that, in the present setting, the representation theory remains unchanged in dimensions at least three. Indeed, we show below that several results previously known only in dimension three, such as the derived equivalence of MMAs and Adachi--Iyama--Reiten-type correspondences, extend to higher-dimensional singularities satisfying the assumptions of our setting.

Next, we state a remarkable result specific to the present setting, namely, the existence of two-term resolutions of modifying modules.

\begin{lem}[= Theorem \ref{MM resol general}]
    If $N\in\refl R$ is a modifying module, then $\overline{\Hom_R(M,N)}\cong \Hom_{\oR}(\oM^{**},\oN^{**})$. In particular, $\pd_{\bLam_M}\Hom_R(M,N)\le 1$.
\end{lem}

The above theorem fails completely for general higher-dimensional singularities. This is an interesting phenomenon specific to our setting.

\begin{nota}
    For $N\in\modif(R)$, we set
     \[\bLam_N:=\End_R(N),
        \quad \bLam_{\oN^{**}}:=\End_{\oR}(\oN^{**}),\]
    \[T_N:=\Hom_R(M,N),\quad T_{\oN^{**}}:=\Hom_{\oR}(\oM^{**},\oN^{**}).\]
\end{nota}

Next, as an application, we show that tilting $\bLam_M$-modules can be constructed in this setting by a method analogous to that in \cite{Maximal-modifications-and-Auslander-Reiten-duality-for-non-isolated-singularities}.

\begin{thm}[= Theorem \ref{tilting of NCCR}]\label{intro tilting of NCCR}
Let $N\in\refl (R)$. Then the following assertions hold.
\begin{enumerate}
    \item If $N\in\modif(R)$, then $T_N=\Hom_R(M,N)$ is a partial tilting $\bLam_M$-module.
    \item If $N\in\modif(R)$, then $T_{\oN^{**}}=\Hom_{\oR}(\oM^{**},\oN^{**})$ is a partial tilting $\bLam_{\oM^{**}}$-module.
    \item If $N\in\MM(R)$, then $T_N=\Hom_R(M,N)$ is a tilting $\bLam_M$-module.
    \item If $N\in\MM(R)$, then $T_{\oN^{**}}=\Hom_{\oR}(\oM^{**},\oN^{**})$ is a tilting $\bLam_{\oM^{**}}$-module.
\end{enumerate}
Moreover, if $N\in\MM(R)$, then there exists the following commutative diagram, in which all the vertical functors are equivalences of categories:
\[\begin{tikzcd}
	{\Db(\modcat\bLam_M)} && {\Db(\modcat\bLam_{\oM^{**}})} \\
	\\
	{\Db(\modcat\bLam_N)} && {\Db(\modcat\bLam_{\oN^{**}})}
	\arrow["{-\underset{R}{\Lo} \oR}", from=1-1, to=1-3]
	\arrow["{\RHom_{\bLam_M}(T_N,-)}"', from=1-1, to=3-1]
	\arrow["{\RHom_{\bLam_{\oM^{**}}}(T_{\oN^{**}},-)}", from=1-3, to=3-3]
	\arrow["{-\underset{R}{\Lo} \oR}"', from=3-1, to=3-3]
\end{tikzcd}\]
\end{thm}

The noncommutative Bondal--Orlov conjecture predicts that all MMAs are derived equivalent. This conjecture has been established in dimension three, but remains open in general dimensions. The above theorem shows that the noncommutative Bondal--Orlov conjecture holds in our setting. Moreover, the commutativity of the diagram is also a phenomenon specific to our setting. Indeed, it is known that the diagram in the theorem does not commute in general in the setting of cDV singularities and their Du Val sections.

\begin{nota}
    For $N\in\modif(R)$, we set
     \[\bsL_N:=\stEnd_R(N)
        \quad\text{and}\quad \bsL_{\oN^{**}}:=\stEnd_{\oR}(\oN^{**}).\]
\end{nota}

Next, as an application, we establish the following Adachi--Iyama--Reiten-type correspondence over $R$.

\begin{thm}[= Theorem \ref{MM(R)->MM(R/t) MM}, Theorem \ref{2-siltA->2-siltA/t}]\label{intro 2-silt corr}
    Assume that $M\in\CM R$. Then there exists the following commutative diagram compatible with mutation:
     \[
    \begin{tikzcd}
	{\MMG(R)} && {\MMG(\overline{R})} \\
	{\twosilt(\bsL_M)} && {\twosilt(\bsL_{\oM})}
	\arrow["\sim", from=1-1, to=1-3]
	\arrow["\sim", from=1-1, to=2-1]
	\arrow["\sim", from=1-3, to=2-3]
	\arrow["\sim", from=2-1, to=2-3]
\end{tikzcd}
\qquad
    \begin{tikzcd}
	{\modifG(R)} && {\modifG(\oR)} \\
	{\twopresilt(\bsL_M)} && {\twopresilt(\bsL_{\oM})}
	\arrow["\sim", from=1-1, to=1-3]
	\arrow["\sim", from=1-1, to=2-1]
	\arrow["\sim", from=1-3, to=2-3]
	\arrow["\sim", from=2-1, to=2-3]
\end{tikzcd}\]
Here, the correspondence is given by Theorem \ref{2-siltA->2-siltA/t}.
\end{thm}

For three-dimensional singularities, the correspondence between MM generators and two-term silting complexes over stable endomorphism algebras is well known. However, this correspondence fails completely for higher-dimensional singularities in general. Thus, this is another result specific to our setting. Moreover, the correspondence 
\[\overline{(-)}:\twosilt (\bsL_M)\xrightarrow{\sim}\twosilt(\bsL_{\oM})\]
is given by Kimura and Gnedin's reduction theorem, and the bijection in Theorem \ref{intro 2-silt corr} is compatible with this correspondence.

Finally, we discuss the spaces of stability conditions. For a ring $A$, we denote by $\fl A$ the category of finite length right $A$-modules. When $A$ is an MMA, the space of stability conditions on $\D^b(\fl A)$ has been extensively studied, particularly its subspace $\Stab^{\alg}$ consisting of stability conditions arising from algebraic t-structures. The following theorem follows from the more general Corollary \ref{intro Stab base change}.

\begin{thm}[= Corollary \ref{stEnd Stab reduction}]
There exists an isomorphism of algebraic stability spaces:
\[\overline{(-)}:\Stab^{\alg} \Db(\fl \bLam_M)\xrightarrow{\cong} \Stab^{\alg}\Db(\fl \bLam_{\oM^{**}}).\]
Moreover, if $M\in\CM R$, then there exists an isomorphism:
\[
\overline{(-)}:\Stab^{\alg} \Db(\fl \bLam_M^{\st})\xrightarrow{\cong} \Stab^{\alg}\Db(\fl \bLam_{\oM}^{\st}).\]
For $\Stab^{\alg}$, see Definition \ref{Stab^alg}.
\end{thm}

These results reveal interesting new phenomena. Moreover, they allow us to study higher-dimensional and three-dimensional singularities in relation to each other. Motivated by these results, we propose the notions of {\it universal singularities} and {\it universal MMAs}. We discuss these notions in more detail in subsection \ref{univ sing}.

\subsection{Tilting theory for stable endomorphism algebras}
In this subsection, we retain Setup \ref{intro setup  reduction NCCR,tilt,2-silt} and use the following notation:
\begin{setup}\label{intro setup tilting stEnd}
We assume that $\oR$ is an isolated cDV singularity over $\bC$ and $M\in\CM R$.
\end{setup}

In this subsection, we consider the tilting theory of stable endomorphism ring $\bsL_N$. It is known that the stable endomorphism algebras of MM generators share various properties with MMAs. One such property is that, under suitable assumptions, all these stable endomorphism algebras are derived equivalent, as shown by August (\cite{The-tilting-theory-of-contraction-algebras}). This result is natural from the viewpoint of the homological MMP and is fundamental to the theory. In this subsection, we generalize August's result to our setting.

\begin{nota}
For $N\in\MMG(R)$, we set
        \[\scrT_N:=\tau^{\ge-1}( \bsL_{N}\underset{\bLam_N}{\Lo}T_N\underset{\bLam_M}{\Lo} \bsL_{M})\quad\text{and}\quad  \scrT_{\oN}:=\tau^{\ge-1}( \bsL_{\oN}\underset{\bLam_{\oN}}{\Lo} T_{\oN}\underset{\bLam_{\oM}}{\Lo} \bsL_{\oM}).\]
\end{nota}

This complex is a two-sided analogue of the two-term silting complex over $\bsL_M$ obtained in Theorem \ref{2-siltA->2-siltA/t}.

\begin{thm}[= Theorem \ref{tilting of NCCR}, Theorem \ref{tilting for contraction algebra}, Proposition \ref{tilting composition}]\label{intro tilting}
For every $N\in\MMG(R)$, the following assertions hold:
\begin{enumerate}
    \item $\scrT_{\oN}$ is a two-sided tilting $(\bsL_{\oN},\bsL_{\oM})$-complex. 
    \item $\scrT_N$ is a two-sided tilting $(\bsL_N,\bsL_M)$-complex.
\end{enumerate}
Moreover, there exists the following commutative diagram, in which all the vertical functors are equivalences of categories:
 \[\begin{tikzcd}
	{\Db(\modcat \bsL_M)} && {\Db(\modcat \bLam_M)} & \\
	& {\Db(\modcat \bsL_{\oM})} && {\Db(\modcat \bLam_{\oM})} \\
	{\Db(\modcat \bsL_{N})} && {\Db(\modcat \bLam_{N})} \\
	& {\Db(\modcat \bsL_{\oN})} && {\Db(\modcat \bLam_{\oN})}
	\arrow[from=1-1, to=1-3]
	\arrow["{-\underset{R}{\Lo}\oR}"{description}, from=1-1, to=2-2]
	\arrow["{\RHom_{\bsL_{M}}(\scrT_N,-)}"{description, pos=0.3}, from=1-1, to=3-1]
	\arrow["{-\underset{R}{\Lo}\oR}"{description}, from=1-3, to=2-4]
	\arrow["{\RHom_{\bLam_M}({T_N},-)}"{description, pos=0.3}, from=1-3, to=3-3]
	\arrow[from=2-2, to=2-4]
	\arrow["{\RHom_{\bsL_{\oM}}(\scrT_{\oN},-)}"{description, pos=0.3}, from=2-2, to=4-2]
	\arrow["{\RHom_{\bLam_{\oM}}(T_{\oN},-)}"{description, pos=0.3}, from=2-4, to=4-4]
	\arrow[from=3-1, to=3-3]
	\arrow["{-\underset{R}{\Lo}\oR}"{description}, from=3-1, to=4-2]
	\arrow["{-\underset{R}{\Lo}\oR}"{description}, from=3-3, to=4-4]
	\arrow[from=4-2, to=4-4]
\end{tikzcd}\]
\end{thm}

\subsection{Reduction for stability conditions}
\begin{setup}\label{intro setup reduction stab}
    In this subsection, we use the following notation:
    \begin{enumerate}
        \item  Let $R$ be a complete local Noetherian non-positive commutative dg ring (see Definition \ref{dg def}), and $(H^0(R),\mm,k)$ be its degree-zero local ring and residue field. 
        \item Let $A$ be a dg $R$-algebra such that $H^{>0}A=0$ and $H^iA\in\modcat R$ for each $i\in\ZZ$.
    \end{enumerate}
\end{setup}

\begin{nota}
We set
\[\overline{(-)}:=-\underset{R}{\Lo}k: \D(A)\to \D(\oA),\]
 \[\pvd A:=\{X\in\D(A)\mid \bigoplus_{i\in\ZZ} H^i X\in\fl H^0 R\},\]
 \[\pvd \oA:=\{X\in\D(\oA)\mid \bigoplus_{i\in\ZZ} H^i X\in\fl H^0 R\}.\]
\end{nota}

We prove an isomorphism between the spaces of  algebraic stability conditions.

\begin{thm}[= Theorem \ref{Stab base change}]\label{1.18}
    There exists an isomorphism between the spaces of algebraic stability conditions:
    \[\Stab^{\alg}(\pvd A)\xrightarrow[]{\sim}\Stab^{\alg}(\pvd \oA). \]
    For $\Stab^{\alg}$, see Definition \ref{Stab^alg}.
\end{thm}

Let $(R,\mm,k)\to (S,\nn,k)$ be a local morphism (see Definition \ref{dg def}) of complete local Noetherian non-positive commutative dg rings. Applying Theorem \ref{1.18} to $(R,\mm,k)$ and $(S,\nn,k)$, we obtain the following.

\begin{cor}[= Corollary \ref{Stab general base change}]\label{intro Stab base change}
    There exists a natural isomorphism between the spaces of algebraic stability conditions:
    \[\Stab^{\alg}(\pvd A)\xrightarrow[]{\sim}\Stab^{\alg}(\pvd A\;\underset{R}{\Lo}\; S). \]
\end{cor}

\subsection*{Use of AI}
The author made extensive use of OpenAI's ChatGPT (GPT-6) in writing, translating, and revising the English text of this paper. It was also used to correct LaTeX code and minor typographical errors, and to assist in developing the proofs of Lemma \ref{iso +1} (2), Lemma \ref{CM +1} and Theorem \ref{Stab base change}. The author has carefully reviewed the AI-assisted content, verified the mathematical arguments, and takes full responsibility for the contents of this paper.

\subsection*{Conventions}
For a ring $A$ and right $A$-modules $M$ and $N$, we denote by
\begin{itemize}
    \item $\Modcat A$ the category of right $A$-modules
    \item $\modcat A$ the category of finitely generated right $A$-modules
    \item $\proj A$ the category of finitely generated projective right $A$-modules
    \item $\fl A$ the category of finite length right $A$-modules
    \item $\refl A$ the category of finitely generated reflexive right $A$-modules
    \item $\add M$ the smallest full subcategory of $\Modcat A$ containing $M$ closed under finite direct sums and direct summands
    \item $\Omega_A M$ the syzygy of $M$
    \item $\stHom_A(M,N)$ the stable Hom set.
\end{itemize}
For a Cohen--Macaulay local ring $R$, we denote by $\CM R$ the category of maximal Cohen--Macaulay $R$-modules.

We compose paths from right to left, following the convention for composition of maps. Thus, if $a:i\to j$ and $b:j\to k$, then $ba:i\to k$.
\section*{Acknowledgements}
I thank Professor Osamu Iyama for welcoming me into his research group and for his valuable guidance. I am grateful to  Ryu Tomonaga, Yuki Hirano and Wahei Hara for valuable discussions and advice and to Yuno Suzuki for many helpful comments and suggestions on the English and the presentation of this paper. 

\section{Preliminaries}
\subsection{Reflexive modules}

In this subsection, we recall some basic facts about reflexive modules. Let $R$ be a commutative Noetherian ring and $M$ a finitely generated $R$-module. We say that $M$ satisfies Serre's condition $(S_n)$ if the following inequality holds for every prime ideal $\pp\in\Spec R$:
\[\depth_{R_{\pp}}M_{\pp}\ge\min\{n,\height \pp\}.\]

The following lemma is elementary but very useful. We will use it repeatedly without further mention.

\begin{lem}\label{ht 1 check}
Let $R$ be a commutative Noetherian ring, and let $f:M\to N$ be a morphism in $\modcat R$. Assume that $M$ satisfies Serre's $(S_2)$-condition and that $N$ satisfies Serre's $(S_1)$-condition. If $f_{\pp}:M_{\pp}\to N_{\pp}$ is an isomorphism for every $\pp\in\Spec R$ with $\height \pp\le 1$, then $f$ is an isomorphism.
\end{lem}

\begin{lem}\label{ref S_2}
    Let $R$ be a commutative Noetherian normal domain. For $M\in\modcat R$, the following conditions are equivalent.
    \begin{enumerate}
        \item $M$ is a reflexive $R$-module.
        \item $M\in\Omega^2(\modcat R)$. 
        \item $M$ satisfies Serre's $(S_2)$-condition.
    \end{enumerate}
\end{lem}
% \begin{proof}
%     (1)$\Rightarrow$(2) Take a projective presentation $P_1\to P_0\to M^*\to 0$. Applying $(-)^*$, we obtain $0\to M^{**}\cong M\to P_0^*\to P_1^*$. 

%     (2)$\Rightarrow$(3) This follows from the depth lemma.

%     (3)$\Rightarrow$(1) This follows from Lemma \ref{ht 1 check}.
% \end{proof}

\begin{lem}\label{torl S_1}
    Let $R$ be a commutative Noetherian domain. For $M\in\modcat R$, the following conditions are equivalent.
    \begin{enumerate}
        \item $M$ is a torsionless $R$-module.
        \item $M\in\Omega(\modcat R)$. 
        \item $M$ satisfies Serre's $(S_1)$-condition.
    \end{enumerate}
\end{lem}
% \begin{proof}
%     (1)$\Rightarrow$(2)$\Rightarrow$(3) The argument is similar to that of Lemma \ref{ref S_2}.

%     (3)$\Rightarrow$(1) Let $K:=\Ker (M\to M^{**})$. We show $K=0$. Suppose not. We can choose $\pp\in\Ass K$. Since $R$ is a domain, $M\otimes_R \Frac R\to M^{**}\otimes_R \Frac R$ is an isomorphism; thus $\height\pp\ge 1$. However, since $M$ satisfies $(S_1)$, we have $\depth K_{\pp}\ge 1$. This contradicts the fact that $\pp\in\Ass K$. Therefore, $K=0$.
% \end{proof}

\begin{cor}\label{refl regular quot}
    Let $R$ be a commutative Noetherian normal domain and $t\in R$ such that $R/tR$ is a domain. If $M\in\modcat R$ is a reflexive $R$-module, then $M/tM$ is a torsionless $R/t$-module.
\end{cor}

Next, we first prepare a simple lemma for later use.

\begin{lem}\label{ref Hom lemma}
    Let $R$ be a commutative Noetherian normal domain and $M,N\in\refl R$, there is a natural isomorphism $\Hom_R(M,N)^{*}\cong \Hom_R(M^*,N^*)$.
\end{lem}
\begin{proof}
    There is a natural morphism $M^*\underset{R}{\otimes}N\to \Hom_R(M,N)$. Thus we have a natural morphism
    \[\Hom_R(\Hom_R(M,N),R)\to \Hom_R(M^*\underset{R}{\otimes} N,R)\cong\Hom_R(M^*,N^*).\]
    Therefore, the claim follows from Lemma \ref{ht 1 check}.
\end{proof}

% We first recall some basic facts about the divisor class group.

% \begin{prop}[{e.g. \cite{bourbaki1998commutative},\cite[2.5.6]{The-K-theory-of-cDV-singularities}}]\label{Groth group -> class group}
%     For a commutative Noetherian normal domain $R$, there exists a group isomorphism \[\det:K_0(\modcat R)/\langle [R],[X]\mid \dim X\le \dim R-2\rangle\to \Cl (R).\]
% \end{prop}
% \begin{lem}[{e.g. \cite[3.1]{Gorenstein-modifications-and-Q-Gorenstein-rings}}]\label{DITW fact}
%     Let $R$ be a commutative Noetherian normal domain. Suppose that $M,N\in\modcat R$ have positive rank and that $M$ is torsion-free. Then $[\det \Hom_R(M,N)]=(\rank M)[\det N]-(\rank N)[\det M]$ in $\Cl(R)$.
% \end{lem}

Finally, we recall that the reflexive equivalence.

\begin{prop}[{\cite[Proposition 2.4]{Fomin-Zelevinsky-mutation-and-tilting-modules-over-Calabi-Yau-algebras.}}]\label{refl equiv}
     Let $R$ be a commutative Noetherian normal domain and let $M\in\refl R$ be a non-zero reflexive module. Then the functor
     \[\Hom_R(M,-):\refl R\to \refl \End_R(M)\]
     is categorical equivalence.
\end{prop}

\subsection{Orders and Calabi-Yau algebras}
\begin{setup}In this subsection, we use the following notation:
    \begin{enumerate}
        \item Let $R$ be a complete Cohen--Macaulay local ring with canonical module $\omega_R$.
        \item Let $\Lambda$ be a module-finite $R$-algebra.
    \end{enumerate}
\end{setup}
First, we recall the notion of an order. An order can be viewed as a non-commutative generalization of a (commutative) Cohen-Macaulay ring, and it provides a framework in which one can develop Cohen---Macaulay representation theory, including Auslander--Reiten theory.
\begin{dfn}
    \begin{enumerate}
        \item $\Lambda$ is called an {\it $R$-order} if $\Lambda\in\CM R$.
        \item We define the {\it Cohen--Macaulay category $\CM\Lambda$ of $\Lambda$} by $\CM\Lambda:=\{M\in\modcat\Lambda\mid M\in\CM R\}$.
        \item We define the {\it canonical module $\omega_\Lambda$ of $\Lambda$} by $\omega_\Lambda:=\Hom_R(\Lambda,\omega_R)$.
    \end{enumerate}
\end{dfn}
Next, we recall the notions corresponding to Gorenstein and regular rings, which form important classes of Cohen--Macaulay rings. 
\begin{dfn}
Let $\Lambda$ be an $R$-order. We say
    \begin{enumerate}
        \item $\Lambda$ is {\it Gorenstein} if $\omega_{\Lambda}$ is projective as a $\Lambda$-module.
        \item $\Lambda$ is {\it symmetric} if $\Lambda\cong \omega_\Lambda$ as a $(\Lambda,\Lambda)$-bimodule.
        \item $\Lambda$ is {\it regular} if $\gldim\Lambda=d$. 
    \end{enumerate}
\end{dfn}
When $\Lambda=R$, being a Gorenstein order is equivalent to being a symmetric order, and these conditions are also equivalent to $R$ being Gorenstein as a commutative ring. However, for a general $\Lambda$, the notions of a Gorenstein order and a symmetric order are no longer equivalent.

Finally, we recall the notions of (singular) Calabi-Yau algebras.
\begin{dfn}
We say that:
\begin{enumerate}
    \item $\Lambda$ is {\it singular d-Calabi-Yau (=singular d-CY)} if there exists a functorial isomorphism $(\star)$ for all $X\in\Db(\fl\Lambda)$ and $Y\in\Kb(\proj\Lambda)$.

    \item $\Lambda$ is {\it d-Calabi-Yau (=d-CY)} if there exists a functorial isomorphism $(\star)$ for $X\in\Db(\fl\Lambda)$ and $Y\in\Db(\modcat\Lambda)$.
\[ \tag{$\star$}
    \Hom_{\D(\Modcat\Lambda)}(X,Y[d])\cong D\Hom_{\D(\Modcat \Lambda)}(Y,X) 
\]

where $D$ denotes the Matlis dual.
\end{enumerate}
\end{dfn}
We recall an important characterization of the (singular) Calabi-Yau property.
\begin{prop}[{\cite[Theorem 3.2]{Fomin-Zelevinsky-mutation-and-tilting-modules-over-Calabi-Yau-algebras.}}]\label{symm=sCY, regular=CY}
If $R$ is Gorenstein and $\Lambda$ is an $R$-order, then the following hold.
\begin{enumerate}
    \item $\Lambda$ is singular d-CY if and only if $\Lambda$ is symmetric.
    \item $\Lambda$ is d-CY if and only if $\Lambda$ is regular and symmetric.
\end{enumerate}
\end{prop}

\subsection{Non-commutative crepant resolutions}
\begin{setup}
    Throughout this subsection, let $R$ be a complete Gorenstein local normal domain of dimension $d:=\dim R\ge 2$, and let $M$ be a reflexive $R$-module.
\end{setup}

\begin{dfn}
    \begin{enumerate} 
        \item $M$ is called {\it modifying} if $\End_R(M)$ is an $R$-order. In this case, $\End_R(M)$ is called {\it modifying algebra} of $R$. 
        \item $M$ is called {\it maximal modifying (MM)} if 
        \[\add M=\{X\in\refl R\mid \End_R(M\oplus X)\in\CM R\}.\]
        In this case, $\End_R(M)$ is called {\it maximal modifying algebra (MMA)} of $R$. 
        % \item $M$ is called {\it CT} if
        % \[\add M=\{X\in\CM R\mid \Hom_R(M,X)\in\CM R\}=\{X\in\CM R\mid \Hom_R(X,M)\in\CM R\}.\]
        \item We say that $M$ gives a non-commutative crepant resolution (NCCR) if $\End_R(M)$ is a regular $R$-order. In this case, $\End_R(M)$ is called a non-commutative crepant resolution of $R$.
    \end{enumerate}
\end{dfn}
\begin{dfn}\label{2.14}
    We define
    \begin{itemize}
        \item $\modif (R):=\{M\in\refl R\mid \text{$M$ is a basic modifying $R$-module}\}/\cong$.
        \item $\modifG(R):=\{M\in\refl R\mid \text{$M$ is a basic modifying $R$-module such that $R\in\add M$}\}/\cong$.
        \item $\MM (R):=\{M\in\refl R\mid \text{$M$ is a basic MM $R$-module}\}/\cong$.
        \item $\MMG(R):=\{M\in\refl R\mid \text{$M$ is a basic MM $R$-module such that $R\in\add M$}\}/\cong$.
        \item $\NCCR(R):=\{M\in\refl R\mid \text{$M$ is  basic and gives an NCCR of $R$}\}/\cong$.
        \item $\NCCRG(R):=\{M\in\refl R\mid \text{$M$ is basic and gives an NCCR of $R$ such that $R\in\add M$}\}/\cong$.
        % \item $\CT(R):=\{M\in\CM R\mid \text{$M$ is a basic CT $R$-module}\}/\cong$.
    \end{itemize}
\end{dfn}

From the viewpoint of the homological MMP, modifying algebras, MMAs, and NCCRs can be regarded as noncommutative counterparts of crepant partial resolutions, minimal models, and crepant resolutions, respectively. Various results concerning these objects are known. In the following, we recall those needed in this paper.

\begin{prop}[{\cite[Section 5]{Maximal-modifications-and-Auslander-Reiten-duality-for-non-isolated-singularities}}]\label{MM NCCR CT summary}
   The following statements hold.
    \begin{enumerate}
        % \item If $M$ is a CT module, then $M$ gives an NCCR of $R$.
        \item If $M$ gives an NCCR, then $M$ is an MM $R$-module.
        \item Assume that  $R$ is a 3-dimensional Gorenstein ring and that $R$ has an NCCR. If $M$ is an MM module, then $M$ gives an NCCR.
        % \item Assume that $M$ is a generator (i.e. $R\in\add M$). If $M$ gives an NCCR, then $M$ is a CT module.
    \end{enumerate}
\end{prop}

% When $R$ has an isolated singularity, the notion of a {\it $(d-1)$-cluster tilting} module $M$ is defined by the following condition.
% \[\add M=\{X\in\CM R\mid \Ext_R^{i}(M,X)=0\ (1\le i\le d-2)\}=\{X\in\CM R\mid \Ext_R^{i}(X,M)=0\ (1\le i\le d-2)\}.\]
% Moreover, in this case, CT modules coincide with $(d-1)$-cluster tilting modules by \cite{Higher-dimensional-Auslander-Reiten-theory-on-maximal-orthogonal-subcategories}.

% \begin{lem}[\cite{Higher-dimensional-Auslander-Reiten-theory-on-maximal-orthogonal-subcategories}]\label{iso-sing Modif vs Ext}
%     Assume that $R$ is an isolated singularity. For $M,N\in\CM R$, the condition $\Hom_R(M,N)\in\CM R$ is equivalent to $\Ext_R^{i}(M,N)=0\ (1\le i\le d-2)$. In particular, being a CT module is equivalent to being a $(d-1)$-cluster tilting module when $R$ has an isolated singularity.
% \end{lem}

% The following proposition is a generalization of the above proposition. Recall that,  for a commutative Noetherian local ring $(A,\mm)$ and an $A$-module $M$, we define $\Sing(A):=\{\pp\in\Spec A\mid \gldim A_{\pp}=\infty\}$ and $\text{$\fl M:=\{x\in M\mid x\mm^i =0$ for some $i$ \}}$.
% \begin{prop}[\cite{Maximal-modifications-and-Auslander-Reiten-duality-for-non-isolated-singularities}]\label{non-iso-sing Modif vs Ext}
%     Assume that $\dim(\Sing R)\le 1$. For $M,N\in\CM R$, the condition $\Hom_R(M,N)\in\CM R$ is equivalent to $\Ext_R^{i}(M,N)=0\ (1\le i\le d-3)$ and $\fl \Ext_R^{d-2}(M,N)=0$.
% \end{prop}

Next, we recall {\it Iyama--Wemyss mutation} of modifying modules. This is an operation that produces new modifying modules from given ones. Geometrically, it corresponds to flops.

\begin{dfn}[\cite{Maximal-modifications-and-Auslander-Reiten-duality-for-non-isolated-singularities}]
Let $M\in\refl R$ be a modifying $R$-module and let $0\ne N\in\add M$. Let $N_0\to M$ be a right $\add N$-approximation of $M$ and let $N_1^{*}\to M^*$ be a right $\add N^{*}$-approximation of $M^{*}$, where $(-)^{*}:=\Hom_R(-,R)$. Take their kernels:

\[0\to K_0\to N_0\to M\quad\text{and}\quad 0\to K_1\to N_1^{*}\to M^*.\]

We define the {\it right mutation of $M$} at $N$ to be $\mu_N^{+}(M):=N\oplus K_0$ and  the {\it left mutation of $M$} at $N$ to be $\mu_N^{-}(M):=N\oplus K_1^{*}$.
\end{dfn}
\begin{prop}[{\cite[Section 6]{Maximal-modifications-and-Auslander-Reiten-duality-for-non-isolated-singularities}}]\label{mutation summary 1}
    The following statements hold.
    \begin{enumerate}
        \item $\mu_N^{+}(M)$ and $\mu_N^{-}(M)$ are modifying $R$-modules.
        \item $\mu_N^{+}\mu_N^{-}(M)\cong\mu_N^{-}\mu_N^{+}(M)\cong M$.
        \item If $M$ gives an NCCR, then so do $\mu_N^{+}(M)$ and $\mu_N^{-}(M)$.
    \end{enumerate}
\end{prop}
Let $M=\bigoplus\limits_i M_i$ be a decomposition into indecomposable modules. We define $\mu_i^+(M):=\mu_{\frac{M}{M_i}}^{+}(M)$ and  $\mu_i^-(M):=\mu_{\frac{M}{M_i}}^{-}(M)$.

\begin{prop}[{\cite[Section 6]{Maximal-modifications-and-Auslander-Reiten-duality-for-non-isolated-singularities}}]\label{mutation summary 2}
Assume $\dim R=3$. Let $M=\bigoplus\limits_i M_i$ be a basic maximal modifying module. Then
\begin{enumerate}
    \item $\End_R(M)/[\bigoplus\limits_{j\ne i}M_j]$ is artinian if and only if $\mu_i^+(M)=\mu_i^-(M)\ne M$.
    \item $\End_R(M)/[\bigoplus\limits_{j\ne i}M_j]$ is not artinian if and only if $\mu_i^+(M)=\mu_i^-(M)= M$.
    \item $\mu_i(M):=\mu_i^+(M)=\mu_i^-(M)$ is an MM module, and $\mu_i\mu_i(M)\cong M$.
\end{enumerate}    
\end{prop}
% In this case, we define $\EG(\MMG(R))$ as follows and call it the {\it exchange graph}. The vertex set is $\MMG(R)$, and there is an arrow from $M$ to $\mu_i(M)$ for each $i$.

% \begin{prop}[\cite{Reduction-of-triangulated-categories-and-maximal-modification-algebras-for-cA_n-singularities}]\label{mutation summary 3}
%     Assume $\dim R=3$. If $\EG(\MMG(R))$ has a finite connected component $C$, then $\EG(\MMG(R))=C$. 
% \end{prop}

% Finally, we recall Auslander--Reiten theory for NCCRs.
% \begin{prop}[{\cite{Higher-dimensional-Auslander-Reiten-theory-on-maximal-orthogonal-subcategories},\cite[Theorem 2.13]{Cohen-Macaulay-representations-of-invariant-subrings}}]
%     Assume that $M\in\refl R$ gives an NCCR of $R$. Let $\mathcal{C}:=\add M$. Then for any $X\in\ind \mathcal{C}$, there exists a $(d-1)$-almost split sequence in $\mathcal{C}$ of the following form.
%      \[X\to C_{d-2}\to \cdots \to C_0\to X.\]
% \end{prop}

\subsection{Space of stability conditions}
Finally, we recall stability conditions on triangulated categories.
\begin{setup}
    Throughout this subsection, $\scrT$ denotes a triangulated category  such that the
Grothendieck group $K_0(\scrT)$ is a finitely generated free $\bZ$-module.
\end{setup}
\begin{dfn}[\cite{Stability-conditions-on-triangulated-categories}]
Let $\scrA$ be an abelian category.
\begin{enumerate}
    \item A {\it stability condition} on $\scrA$ is a map (called a {\it central charge})
    \[Z:K_0(\scrA)\to \bC\]
    satisfying the following conditions:
    \begin{enumerate}
        \item For every $0\ne E\in \scrA$, $Z([E])\in\bH:=\{re^{i\pi\phi}\mid \text{$r\in\bR_{>0}$ and $0<\phi\le 1$}\}$.
        \item $Z$ satisfies the Harder--Narasimhan property (\cite{Stability-conditions-on-triangulated-categories}).
    \end{enumerate}
    \item A {\it stability condition} on $\scrT$ is a pair $(Z,\scrA)$, where $\scrA$ is the heart of a bounded t-structure of $\scrT$ and $Z$ is a stability condition on $\scrA$.
    \item A {\it space of stability conditions} of $\scrT$ is the set of stability conditions $(Z,\scrA)$ satisfying the support property (\cite{Stability-structures-motivic-Donaldson-Thomas-invariants-and-cluster-transformations}). We denote this space by $\Stab\scrT$.
\end{enumerate}
\end{dfn}
\begin{prop}[\cite{Stability-conditions-on-triangulated-categories}]
    The forgetful map
    \[\Stab\scrT\to \Hom_{\bZ}(K_0(\scrT),\bC)\]
    is a local isomorphism. In particular, $\Stab\scrT$ has a complex manifold structure.
\end{prop}

\begin{dfn}\label{Stab^alg}
Let $\scrA$ be an abelian category.
\begin{enumerate}
    \item $\scrA$ is {\it length} if every object of $\scrA$ has finite length. 
    \item Moreover, a stability condition $(Z,\scrA)\in\Stab \scrT$ is called {\it algebraic} if $\scrA$ is a length category with finitely many isomorphism classes of simple objects. We denote the set of algebraic stability conditions by $\Stab^{\alg}\scrT$. 
\end{enumerate}
\end{dfn}
\section{Ascent for NCCRs and MMs}
\begin{setup}
    In this section, we work under Setup \ref{Intro setup ascent NCCR}.
\end{setup}
\subsection{Ascent of regularity}

In this subsection, we study the ascent of certain properties for orders with respect to regular elements. Let $\Lambda$ be a module-finite $R$-algebra such that $t\in\reg\Lambda$.

The following lemmas are well known.
\begin{lem}\label{exactness}
    Let $0\to L\to M\to N\to 0$ be an exact sequence in $\modcat \Lambda$. If $t\in\reg N$, then $0\to \overline{L}\to \overline{M}\to \overline{N}\to 0$ is also exact.
\end{lem}
\begin{proof}
    This follows from the snake lemma.
\end{proof}
\begin{lem}[{e.g. \cite[Corollary 2.3.3]{Towards-a-theory-of-Bass-numbers-with-application-to-Gorenstein-algebras}}]\label{quot by u proj}
    Let $M$ be a finitely generated $\Lambda$-module such that $t$ is $M$-regular. Then $\pd_{\Lambda}M=\pd_{\overline{\Lambda}}\overline{M}$.
\end{lem}

\begin{lem}[{e.g. \cite[Corollary 2.3.4]{Towards-a-theory-of-Bass-numbers-with-application-to-Gorenstein-algebras}}]\label{fgd fid lift}
    If $\gldim\oLam$ is finite, then $\gldim\Lambda\le\gldim\overline{\Lambda}+1$.
\end{lem}
\begin{proof}
Let $\gldim\oLam=n$. Take $M\in\modcat\Lambda$ and consider a projective resolution $0\to \Omega^{n+1}M\to P_{n}\to \cdots \to P_1\to P_0\to M\to 0$. By  Lemma \ref{exactness}, the sequence  $0\to \overline{\Omega^{n+1} M}\to \overline{P_{n}}\to  \cdots \to \overline{P_1}\to \overline{\Omega M}\to 0$ is exact. Since $\gldim\overline{\Lambda}=n$, we see that $\overline{\Omega^{n+1} M}$ is a projective $\overline{\Lambda}$-module. Hence $\Omega^{n+1} M$ is a projective $\Lambda$-module. Therefore $\pd_{\Lambda}M\le n+1$.
\end{proof}

\begin{cor}\label{regular +1}
We assume that $\Lambda$ is an $R$-order. If $\overline{\Lambda}$ is a regular $\overline{R}$-order, then $\Lambda$ is also a regular $R$-order.
\end{cor}

% \begin{dfn}
%     weakly symmetric order
% \end{dfn}
% \begin{rem}
%     symmetric -> weakly symmetric
% \end{rem}
% \begin{prop}
% If $\overline{\Lambda}$ is a weakly symmetric $\overline{R}$-order, then $\Lambda$ is also a weakly symmetric $R$-order.
% \end{prop}
% \begin{proof}
    
% \end{proof}
\begin{ques}
    We assume that $\Lambda$ is an $R$-order. Does it follow that $\Lambda$ is a symmetric $R$-order if $\oLam$ is a symmetric $\oR$-order?
\end{ques}

\subsection{Ascent of NCCRs and MMs}

In this subsection, we give a proof of Corollary \ref{intro NCCR +1}. We begin by preparing several easy lemmas.
\begin{lem}\label{iso +1}
    Let $f:M\to N$ be a morphism in $\modcat R$. Assume that $t\in\reg M\cap\reg N$. 
    \begin{enumerate}
        \item If $\overline{f}$ is surjective, then $f$ is surjective.
        \item If $\overline{f}$ is injective, then $f$ is injective.
    \end{enumerate}
    In particular $\overline{f}$ is an isomorphism, then $f$ is an isomorphism.
\end{lem}
\begin{proof}
(1) This follows from NAK.

(2) Let $K:=\Ker f$. By $t\in\reg M$, we have $K\subseteq tM$. Moreover, by $t\in\reg N$, $K\subseteq tK$. Thus $K=0$ by NAK.
\end{proof}
\begin{cor}\label{iso** +1}
    Let $f:M\to N$ be a morphism in $\refl R$. If $\overline{f}^{**}:\overline{M}^{**}\to\overline{N}^{**}$ is an isomorphism, then $f$ is an isomorphism.
\end{cor}
\begin{proof}
    By Lemma \ref{ht 1 check}, it suffices to show that $f_{P}$ is an isomorphism for each $P\in\Spec R$ with $\height P\le 1$. Let $Q$ be a prime ideal in Lemma \ref{appendixlemma}. Then $R_{Q}$ is regular, since $R_{Q}/t R_{Q}=\oR_{Q/tR}$ is regular. Thus $M_{Q}$ and $N_{Q}$ are projective $R_{Q}$-module, and hence $f_{Q}$ is an isomorphism by Lemma \ref{iso +1}. Therefore, $f_{P}$ is also an isomorphism, since it is a localization of $f_{Q}$.
\end{proof}

The following lemma is a straightforward generalization of a known result.

\begin{lem}\label{Hom quot commu general}
    For $M,N\in \refl R$, there exists a natural isomorphism $\overline{\Hom_R(M,N)}^{**}\cong \Hom_{\overline{R}}(\overline{M}^{**},\overline{N}^{**})$.
\end{lem}
\begin{proof}
We prove only the case $m=1$. The general case follows by repeating the same argument. Let $t:=t_1$. 

Since $0\to N\xrightarrow{t}N\to \overline{N}\to 0$ is exact. Applying $\Hom_R(M,-)$, we have the exact sequence $0\to \Hom_R(M,N)\xrightarrow[]{t}\Hom_R(M,N)\to \Hom_{\overline{R}}(\overline{M},\overline{N})$. Thus we obtain an injection 
\[\overline{\Hom_R(M,N)}\to \Hom_{\overline{R}}(\overline{M},\overline{N}).\]
By composing with the natural morphism 
\[\Hom_{\overline{R}}(\overline{M},\overline{N})\to \Hom_{\overline{R}}(\overline{M}^{**},\overline{N}^{**})\quad f\mapsto f^{**},\] 
we obtain a morphism
\[
\overline{\Hom_R(M,N)}\to \Hom_{\overline{R}}(\overline{M},\overline{N})\to \Hom_{\overline{R}}(\overline{M}^{**},\overline{N}^{**}).
\] 
Moreover, $\Hom_{\overline{R}}(\overline{M}^{**},\overline{N}^{**})\in\refl \oR$, we have a natural morphism
\[
\varphi:\overline{\Hom_R(M,N)}^{**}\to \Hom_{\overline{R}}(\overline{M}^{**},\overline{N}^{**}).
\] 
By a standard argument, we obtain $\varphi_{\pp}$ is an isomorphism for each $\pp\in\Spec \oR$ with $\height\pp\le 1$. Therefore, the assertion follows from Lemma \ref{ht 1 check}.
\end{proof}

\begin{cor}[{e.g. \cite[Proposition 8.20]{IW2023tits}}]\label{Hom quot comm}
    For $M,N\in \refl R$, if $\Hom_R(M,N)\in \CM R$, then there exists a natural isomorphism $\overline{\Hom_R(M,N)}\cong \Hom_{\overline{R}}(\overline{M}^{**},\overline{N}^{**})$.
\end{cor}

% \begin{lem}\label{split epi +1}
%     Let $f:M\to N$ be a morphism in $\refl R$. Suppose that $\Hom_R(N,M)\in\CM R$. 
%     If $\overline{f}^{**}:\overline{M}^{**}\to\overline{N}^{**}$ is a split epi (resp. split mono), then $f$ is a split epi (resp. split mono).
% \end{lem}
% \begin{proof}
%     By assumption, there exists $\overline{g}^{**}:\overline{N}^{**}\to\overline{M}^{**}$ such that $\overline{f}^{**} \circ\overline{g}^{**}=\id_{\overline{N}^{**}}$. Since $\Hom_R(N,M)\in\CM R$, $\overline{g}^{**}$ admits a lift $g$ by Corollary \ref{Hom quot comm}. Since $\overline{f\circ g}^{**}=\overline{f}^{**}\circ\overline{g}^{**}=\id_{\oN^{**}}$, $f\circ g$ is an isomorphism by Corollary \ref{iso** +1}.
% \end{proof}

We now prove Theorem \ref{1.2}.
\begin{thm}\label{NCCR +1}
    Let $M\in\refl R$ be a modifying module. The following statements hold.
    \begin{enumerate}
        \item If $\overline{M}^{**}$ gives an NCCR of $\overline{R}$, then $M$ gives an NCCR of $R$.
        \item If $\overline{M}^{**}$ is an MM $\overline{R}$-module, then $M$ is an MM $R$-module.
    \end{enumerate}
\end{thm}
\begin{proof}
(1)\ This follows from Corollary \ref{Hom quot comm} and Corollary \ref{regular +1}.

(2)\ Let $N\in\refl R$ with $\End_R(M\oplus N)\in\CM R$. It suffices to show that $N\in\add M$. By Corollary \ref{Hom quot comm}, we have $\End_{\oR}(\oM^{**}\oplus\oN^{**})\cong \overline{\End_R(M\oplus N)}\in\CM\oR$. Since $\oM^{**}$ is maximal modifying, we have $\oN^{**}\in\add\oM^{**}$. Thus there exist morphisms $\overline{f}:\oN^{**}\to(\oM^{**})^{\oplus m}$ and $\overline{g}:(\oM^{**})^{\oplus m}\to\oN^{**}$ such that $\overline{g}\circ\overline{f}=\id_{\oN^{**}}$. By Corollary \ref{Hom quot comm}, there exist lifts $f$ and $g$ of $\overline{f}$ and $\overline{g}$, respectively. Since $\overline{g\circ f}=\overline{g}\circ\overline{f}$ is an isomorphism, $g\circ f$ is also an isomorphism by \ref{iso +1}. Therefore, $N\in\add M$.
\end{proof}
\begin{cor}\label{CM NCCR +1}
  Let $M\in\CM R$ be a modifying module. The following statements hold.
    \begin{enumerate}
        \item If $\overline{M}$ gives an NCCR of $\overline{R}$, then $M$ gives an NCCR of $R$.
        \item If $\overline{M}$ is an MM $\overline{R}$-module, then $M$ is an MM $R$-module.
    \end{enumerate}
\end{cor}
\begin{eg}
    If $\overline{R}$ is CM modifying-free (i.e. every CM modifying module is projective), then $R$ is also CM modifying-free.
\end{eg}
Other examples and applications will be discussed later.

\subsection{Iyama--Wemyss mutations}

In this subsection, we study the compatibility between the functor $\overline{(-)}^{**}$ and Iyama--Wemyss mutation. When mutating reflexive modules that are not necessarily Cohen--Macaulay, the argument becomes slightly more involved. We prepare the following lemmas for this purpose.

\begin{lem}\label{Cok t-regular}
Let $f:M\to N$ in $\modcat R$ such that $t\in\reg M\cap\reg N$. If $\overline{f}$ is injective, then $t\in\reg\Cok f$. 
\end{lem}
\begin{proof}
This follows from the Lemma \ref{iso +1} and snake Lemma.    
\end{proof}

Next, we make a simple observation on mutation.

\begin{prop}\label{general lemma}
    Let $f:M\to N$ be a morphism in $\refl R$. If $f_{P}$ is surjective for every $P\in\Spec R$ with $\height P\le 2$ and $t\in P$, then $\overline{\Ker f}^{**}\cong\Ker (\overline{f}^{**})$.
\end{prop}
\begin{proof}
There exists a natural morphism $\overline{\Ker f}\to\Ker (\overline{f}^{**})$. Since $\Ker (\overline{f}^{**})$ is reflexive, this morphism induces $\varphi:\overline{\Ker f}^{**}\to\Ker (\overline{f}^{**})$. It suffices to show that $\varphi_{\pp}$ is an isomorphism for each $\pp\in\Spec \oR$ with $\height \pp\le 1$.

Let $P\in\Spec R$ is a prime ideal corresponding to $\pp\in\Spec \oR$. By assumption, $f_P$ is surjective, thus $\Ker \overline{f_P}\cong \overline{\Ker f_P}$. Moreover, $R_P$ is regular local ring of dimension $\le 2$, so $M_P$ and $N_P$ are projective $R_P$-modules. Thus $\varphi_{\pp}$ is an isomorphism, as required.
\end{proof}

\begin{cor}\label{ker comm}
    Let $M,N\in\refl R$. If $f:M_N\to N$ is a right $\add M$-approximation, then $\overline{\Ker f}^{**}\cong \Ker (\overline{f}^{**})$.
\end{cor}
\begin{proof}
    By Proposition \ref{general lemma}, it suffices to show that $f_{P}$ is surjective for every $P\in\Spec R$ with $\height P\le 2$ and $t\in P$. Since $R_P$ is regular local ring of dimension $\le 2$, we have $\proj R_P=\refl R_P$. Therefore, the assertion follows from the fact that approximation are preserved under localization.
\end{proof}

\begin{cor}\label{mutation compatibility}
    Let $M\in\refl R$ be a modifying module and let $N\in\add M$. Then 
    \[\overline{\mu^{+}_{N}(M)}^{**}\cong \mu^{+}_{\overline{N}^{**}}(\overline{M}^{**})\quad\text{and}\quad \overline{\mu^{-}_{N}(M)}^{**}\cong \mu^{-}_{\overline{N}^{**}}(\overline{M}^{**}).\]
\end{cor}
\begin{proof}
The isomorphism
$\overline{\mu^{+}_{N}(M)}^{**}\cong \mu^{+}_{\overline{N}^{**}}(\overline{M}^{**})$
follows from Lemma \ref{approx =} and Corollary \ref{ker comm}.
Similarly, taking Lemma \ref{Hom quot commu general} into account, we obtain the isomorphism
$\overline{\mu^{-}_{N}(M)}^{**}\cong \mu^{-}_{\overline{N}^{**}}(\overline{M}^{**})$.
\end{proof}

\section{Reduction for NCCRs, tilting modules and two-term silting complexes}\label{NCCR-1}
\begin{setup}
In this section, we work under Setup \ref{intro setup  reduction NCCR,tilt,2-silt}.    
\end{setup}

We prove the following simple lemmas for later use.

\begin{lem}\label{approx =}
    Let $f:M\to N$ be a morphism in $\refl R$. Then the following assertions hold.
    \begin{enumerate}
        \item Suppose that $\Hom_R(M,M\oplus N)\in\CM R$, then 
        \[\text{$f$ is a right $\add M$-approximation $\iff$ $\overline{f}^{**}$ is a right $\add \overline{M}^{**}$-approximation.}\]
        \item  Suppose that $\Hom_R(M\oplus N,N)\in\CM R$, then 
        \[\text{$f$ is a left $\add N$-approximation $\iff$ $\overline{f}^{**}$ is a left $\add \overline{N}^{**}$-approximation.}\]
    \end{enumerate}
\end{lem}
\begin{proof}
    This follows from Corollary \ref{Hom quot comm} and NAK.
\end{proof}

\begin{lem}\label{quot indec isom}
    The following statements hold.
    \begin{enumerate}
        \item If $M\in\refl R$ is an indecomposable modifying module, then $\overline{M}^{**}\in\refl \overline{R}$ is also an indecomposable modifying module.
        \item For $M,N\in\refl R$, if $\oM^{**}\cong \oN^{**}$ and $\Hom_R(M,N)\in\CM R$, then $M\cong N$.
    \end{enumerate}
\end{lem}
\begin{proof}
    (1) This follows from Corollary \ref{Hom quot comm}, since the locality of the endomorphism ring is preserved.
    
    (2) Take an isomorphism $\overline{f}^{**}:\overline{M}^{**}\to\overline{N}^{**}$. By Corollary \ref{Hom quot comm}, there exists a lift $f:M\to N$ of $\overline{f}^{**}$. Then $f$ is an isomorphism by Corollary \ref{iso** +1}.
\end{proof}

\subsection{Reduction of modifying modules}
In this subsection, we construct the bijection for modifying modules appearing in Theorem \ref{intro MM(R)->MM(R/t) MM}. We begin by preparing several lemmas.

\begin{lem}
    For $L,N\in\modcat R$ with $t_1,\dots ,t_m$ is a regular sequence on $L$ there exists a following natural isomorphism for each $i\in\bZ_{\ge 0}$: 
    \[\Ext_R^i(L,\oN)\cong \Ext_{\oR}^i(\oL,\oN).\]
\end{lem}
\begin{proof}
    This follows from adjunction $\RHom_R(L,\oN)\cong \RHom_{\oR}(L\underset{R}{\Lo}\oR,\oN)$.
\end{proof}

Next, we recall the following.

\begin{prop}[{\cite[Theorem 4.12]{Maximal-modifications-and-Auslander-Reiten-duality-for-non-isolated-singularities}}]\label{MM resol}
    Let $M$ be an MM module over a 3-dimensional Gorenstein normal domain. If $N\in\refl R$ is a modifying module, then $\pd_{\End_R(M)}\Hom_R(M,N)\le 1$.
\end{prop}

Our first goal is to extend the above theorem to the present setting. The first key result is the following.

\begin{prop}\label{MM resol preliminary}
    For every $X\in\modif (\oR)$, there exists $N\in\modif (R)$ such that $\oN^{**}\cong X$ and $\overline{\Hom_R(M,N)}\cong \Hom_{\oR}(\oM^{**},\oN^{**})$. In particular, $\pd_{\bLam_M}\Hom_R(M,N)\le 1$.
\end{prop}
\begin{proof}
    \medskip{\it step1: construction of $N$.}
    By Proposition \ref{MM resol}, we have an exact sequence 
    \[0\to \oM_1^{**}\to \oM_0^{**}\xrightarrow{\varphi} X.\]
    with $\varphi$ is a right $\add \oM^{**}$-approximation. Applying $\Hom_{\oR}(\oM^{**},-)$, we have an exact sequence
\[\begin{tikzcd}
	0 & {\Hom_{\oR}(\oM^{**},\oM^{**}_1)} & {\Hom_{\oR}(\oM^{**},\oM^{**}_0)} & {\Hom_{\oR}(\oM^{**},X)} & 0 \\
	& {\overline{\Hom_R(M,M_1)}} & {\overline{\Hom_R(M,M_0)}}
	\arrow[from=1-1, to=1-2]
	\arrow[from=1-2, to=1-3]
	\arrow[from=1-3, to=1-4]
	\arrow[from=1-4, to=1-5]
	\arrow["\cong"', from=2-2, to=1-2]
	\arrow["{\overline{f}}", from=2-2, to=2-3]
	\arrow["\cong"', from=2-3, to=1-3]
\end{tikzcd}\]
    We can take a lift $f\in\Hom_{\End_R(M)}(\Hom_R(M,M_1),\Hom_R(M,M_0))$ of $\overline{f}$. Let $C:=\Cok f$. Then $\overline{C}\cong \Hom_{\oR}(\oM^{**},X)$ and $t\in\reg C$ by Lemma \ref{Cok t-regular}. Thus $C\in\refl \End_R(M)$. By Proposition \ref{refl equiv}, there is $N\in\refl R$, such that $C\cong \Hom_R(M,N)$. 

    \medskip{\it step2: $\oN^{**}\cong X$ and $\overline{\Hom_R(M,N)}\cong \Hom_{\oR}(\oM^{**},\oN^{**})$.}
    Since $\overline{\Hom_R(M,N)}\cong \Hom_{\oR}(\oM^{**},X)$ is reflexive, we obtain $ \overline{\Hom_R(M,N)}\cong  \overline{\Hom_R(M,N)}^{**}\cong \Hom_{\oR}(\oM^{**},\oN^{**})$ by Lemma \ref{Hom quot commu general}. Therefore, $\Hom_{\oR}(\oM^{**},X)\cong\Hom_{\oR}(\oM^{**},\oN^{**})$, we have $X\cong \oN^{**}$ by Proposition \ref{refl equiv}.

    \medskip{\it step3: $N$ is modifying.}
    We prove only the case $m=1$. The general case follows by repeating the same argument. Let $t:=t_1$. We show that $N\in\modif (R)$. Let $T:=\Hom_R(M,N)$ and let $\Lambda:=\End_R(M)$. Then $\oT=\overline{\Hom_R(M,N)}\cong\Hom_{\oR}(\oM^{**},X)$ and $\overline{\Lambda}=\overline{\End_R(M)}\cong \End_{\oR}(\oM^{**})$, thus $\End_{\oLam}(\oT)\cong \End_{\oR}(X)\in\CM \oR$. Hence it enough to show that $\overline{\Hom_\Lambda(T,T)}\cong \Hom_{\oLam}(\oT,\oT)$. We consider an exact sequence
    \[\Ext^1_{\Lambda}(T,T)\xrightarrow{t}\Ext^1_{\Lambda}(T,T)\to \Ext^1_{\Lambda}(T,\oT)=\Ext^1_{\oLam}(\oT,\oT)\]
    induced by $0\to T\xrightarrow{t}T\to \oT\to 0$. Since $\Ext^1_{\oLam}(\oT,\oT)=0$ by \cite[Proposition 4.14]{Maximal-modifications-and-Auslander-Reiten-duality-for-non-isolated-singularities}, we have $\Ext^1_{\Lambda}(T,T)=0$. Therefore $\overline{\Hom_\Lambda(T,T)}\cong \Hom_{\oLam}(\oT,\oT)$.
\end{proof}
\begin{cor}\label{modif surj}
    The map 
    \[\overline{(-)}^{**}:\modif (R)\to \modif (\oR)\]
    is surjective.
\end{cor}

Next, we show that the surjection in Corollary \ref{modif surj} restricts to modifying generators.

\begin{lem}\label{CM +1}
    Let $T\in\refl R$ such that $T_{P}\in\proj R_P$ for every $P\in\Spec R$ with $R_P$ is regular. If $\oT^{**}\in\CM \oR$, then $T\in\CM R$.
\end{lem}
\begin{proof}
    We prove only the case $m=1$. The general case follows by repeating the same argument. Let $t:=t_1$.

    Consider the exact sequence 
     \[0\to\oT\to \oT^{**}\to C\to 0\]
    and the induced exact sequence 
     \[\Ext^1_{\oR}(\oT^{**},\oR)\to\Ext^1_{\oR}(\oT,\oR)\to\Ext^2_{\oR}(C,\oR).\]
     By assumption, $C$ satisfies $\dim C\le \dim\oR-3$. Thus $\Ext^1_{\oR}(\oT^{**},\oR)=0=\Ext^2_{\oR}(C,\oR)$, we have $\Ext^1_{R}(T,\oR)=\Ext^1_{\oR}(\oT,\oR)=0$. Hence $\overline{\Hom_R(T,R)}\cong \Hom_{\oR}(\oT,\oR)$ by an exact sequence $0\to R\xrightarrow{t} R\to \oR$. By assumption, $\Hom_{\oR}(\oT,\oR)\in\CM\oR$; thus $\Hom_R(T,R)\in\CM R$. Since $T$ is reflexive, we obtain $T\in\CM R$.
\end{proof}

To verify the assumptions of the above lemma, we use the following fact.

\begin{prop}[{\cite[Theorem 4.4]{Maximal-orders}}]\label{modif RLR}
    Let $R$ ba a regular local ring and let $M\in\refl R$. Then 
    \[\text{$M\in\proj R \iff M\in\modif (R)$.}\]
\end{prop}

\begin{cor}\label{modif assumption}
    If $N\in\refl R$ is modifying, then $T=N$ satisfies the assumptions of the Lemma \ref{CM +1}.
\end{cor}
\begin{cor}\label{generator +1}
    Let $N\in\modif (R)$. If $\oR\in\add \oN^{**}$, then $R\in\add N$.
\end{cor}
\begin{proof}
    By Lemma \ref{quot indec isom}, the indecomposable decomposition of $N$ induces an indecomposable decomposition of $\oN^{**}$. It suffices to show that, if $\oR\cong \oN^{**}$, then $R\cong N$. By Corollary \ref{modif assumption} and Lemma \ref{CM +1}, $N\in\CM R$; thus $\oN\cong\oN^{**}$. The assertion follows immediately from this.
\end{proof}
\begin{rem}
    The assumption that $\oR$ is isolated singularity cannot be removed from Corollary \ref{generator +1}.
\end{rem}
\begin{proof}
    Let 
    \[R:=k[[x,y,z,w,u]]/(xy+zw)\quad\text{and}\quad \oR:=k[[x,y,z,u]]/(xy+z^2)\quad\text{and}\quad I:=(x^2,xz,z^2).\]
    Then $I\in\refl R$ and $I$ is modifying. Indeed, $R$ is a polynomial extension of the toric ring associated with the weights $(1,1,-1,-1)$, and its divisor class group is isomorphic to $\bZ$. In this case, $I$ correspond to the classes $2$ in the divisor class group. 

    Then, $\oI=(x^2,xz,xy)=x(x,y,z)\cong (x,y,z)$. Thus $\oI^*\cong \oR$. Therefore $\oI^{**}\cong \oR$, but $I\ncong R$.
\end{proof}
\begin{cor}\label{modif_0 surj}
     The map 
    \[\overline{(-)}:\modifG (R)\to \modifG (\oR)\]
    is surjective.
\end{cor}
\begin{proof}
This follows from Corollary \ref{modif surj} and Corollary \ref{generator +1}.
\end{proof}

Finally, we consider the injectivity of the map in Corollary \ref{modif surj}. 
\begin{prop}
    The map 
    \[\overline{(-)}^{**}:\modif (R)\to \modif (\oR)\]
    is bijective.
\end{prop}
\begin{proof}
Let $L,N\in\modif(R)$ with $\oL^{**}\cong \oN^{**}$. We show that $L\cong N$. By Lemma \ref{quot indec isom}, it suffices to show that $\Hom_R(L,N)\in\CM R$. This follows from Corollary \ref{modif assumption} and Lemma \ref{CM +1}.
\end{proof}

\begin{cor}\label{modif bij}
    There exist the following bijections:
\[\begin{tikzcd}
	{\modif(R)} && {\modif(\oR)} \\
	{\modifG(R)} && {\modifG(\oR)}
	\arrow["\sim", from=1-1, to=1-3]
	\arrow["{\overline{(-)}^{**}}"', from=1-1, to=1-3]
	\arrow[hook', from=2-1, to=1-1]
	\arrow["\sim", from=2-1, to=2-3]
	\arrow["{\overline{(-)}}"', from=2-1, to=2-3]
	\arrow[hook', from=2-3, to=1-3]
\end{tikzcd}\]
\end{cor}

\subsection{Reduction for NCCRs and MMs}

In this subsection, we prove Theorem \ref{intro MM(R)->MM(R/t) MM}. More precisely, by restricting the bijection in Corollary \ref{modif bij}, we obtain mutation-preserving bijections between maximal modifying $R$-modules and maximal modifying $\oR$-modules, and between maximal modifying generators over $R$ and $\oR$.

\begin{thm}\label{MM(R)->MM(R/t) MM}
    There exist the following bijections, compatible with the Iyama--Wemyss mutation:
\[\begin{tikzcd}
	{\MM(R)} && {\MM(\oR)} & {\modif(R)} && {\modif(\oR)} \\
	{\MMG(R)} && {\MMG(\oR)} & {\modifG(R)} && {\modifG(\oR)}
	\arrow["{\overline{(-)}^{**}}"', from=1-1, to=1-3]
	\arrow["\sim", from=1-1, to=1-3]
	\arrow["{\overline{(-)}^{**}}"', from=1-4, to=1-6]
	\arrow["\sim", from=1-4, to=1-6]
	\arrow[hook', from=2-1, to=1-1]
	\arrow["{\overline{(-)}}"', from=2-1, to=2-3]
	\arrow["\sim", from=2-1, to=2-3]
	\arrow[hook', from=2-3, to=1-3]
	\arrow[hook', from=2-4, to=1-4]
	\arrow["{\overline{(-)}}"', from=2-4, to=2-6]
	\arrow[from=2-4, to=2-6]
	\arrow["\sim", from=2-4, to=2-6]
	\arrow[hook', from=2-6, to=1-6]
\end{tikzcd}\]
\end{thm}
\begin{proof}
    \medskip{\it step1: well-definedness.}
    It suffices to show that $\overline{(-)}^{**}:\MM(R)\to \MM(\oR)$ is well-defined. Let $N\in\MM(R)$. Suppose that $\oN^{**}\oplus X\in\modif (\oR)$ for $X\in\refl\oR$. Then there exists $L\in\modif(R)$ such that $\oL^{**}\cong \oN^{**}\oplus X$. Moreover, by the injectivity of $\overline{(-)}^{**}:\modif(R)\to\modif(\oR)$ and Lemma \ref{quot indec isom}, there exists a decomposition $L\cong N\oplus L'$ such that $\overline{L'}^{**}\cong X$. Since $N\in\MM(R)$, we have $L'\in\add N$; thus $X\cong \overline{L'}^{**}\in\add\oN^{**}$.

    \medskip{\it step2: bijectivity.} The injectivity follows from Corollary \ref{modif bij}. The surjectivity  follows from Theorem \ref{NCCR +1} and Corollary \ref{modif bij}. 
    
    \medskip{\it step3: compatibility with mutations} This follows from Corollary \ref{mutation compatibility}.
\end{proof}

We next discuss the NCCR version of the above theorem. 

\begin{thm}\label{NCCR(R)->NCCR(R/t) NCCR}
    If $\overline{M}^{**}$ gives an NCCR of $\oR$, then there exist the following bijections, compatible with the Iyama--Wemyss mutation:
    \[\begin{tikzcd}
	{\NCCR(R)} && {\NCCR(\oR)} \\
	{\NCCRG(R)} && {\NCCRG(\oR)}
	\arrow["{\overline{(-)}^{**}}"', from=1-1, to=1-3]
	\arrow["\sim", from=1-1, to=1-3]
	\arrow[hook', from=2-1, to=1-1]
	\arrow["{\overline{(-)}}"', from=2-1, to=2-3]
	\arrow["\sim", from=2-1, to=2-3]
	\arrow[hook', from=2-3, to=1-3]
\end{tikzcd}\]
\end{thm}
\begin{proof}
    By Proposition \ref{MM NCCR CT summary}, $\NCCR(\oR)=\MM(\oR)$. Thus $\NCCR(R)=\MM(R)$ and $\NCCRG(R)=\MMG(R)$ by Corollary \ref{CM NCCR +1}. Therefore we obtain the claim by Theorem \ref{MM(R)->MM(R/t) MM}.
\end{proof}

\begin{rem}
If $\oR$ is not an isolated singularity, the theorem does not necessarily hold. For example, set
\[R:=k[[x,y,z,w,u]]/(xy+zw).\]
First, consider the regular element $t:=u$. By Theorem \ref{NCCR(R)->NCCR(R/t) NCCR}, we have $\NCCRG(R)\cong \NCCRG(R/(t))$. Hence, the exchange graph is given as follows:
\[\begin{tikzcd}
	{R\oplus(x,z)} & {R\oplus(x,w).}
	\arrow[no head, from=1-1, to=1-2]
\end{tikzcd}\]

Now consider another regular element $t':=z-w$. In this case, the exchange graph of $\NCCRG(R/(t'))$ is given by
\[\begin{tikzcd}
	{R/(t')\oplus(x,z)=R/(t')\oplus(x,w).}
	\arrow[no head, from=1-1, to=1-1, loop, in=55, out=125, distance=10mm]
\end{tikzcd}\]

These two exchange graphs are clearly not isomorphic.
\end{rem}

Finally, we state the following simple but important corollary.

\begin{cor}\label{NCCR-1 cor}
    Let $L,N\in\MMG(R)$. Then the following statements hold.
    \begin{enumerate}
        \item If $L=\bigoplus\limits_{i=1}^{\ell} L_i$ and $N=\bigoplus\limits_{i=1}^n N_i$ are indecomposable decompositions, then $\ell=n$.
        \item For each $i$, $\mu_i^-(N)=\mu_i^+(N)$ holds.
        \item If $\oM$ gives an NCCR of $\oR$, then $\MM(R)=\NCCR(R)$.
    \end{enumerate}
\end{cor}
\begin{proof}
    (1) This follows from Lemma \ref{quot by u proj}, Theorem \ref{MM(R)->MM(R/t) MM} and \cite{On-the-noncommutative-Bondal-Orlov-conjecture}.

    (2) This follows from Theorem \ref{MM(R)->MM(R/t) MM} and Proposition \ref{mutation summary 2}.

    (3) This follows from Theorem \ref{NCCR(R)->NCCR(R/t) NCCR} and Theorem \ref{MM(R)->MM(R/t) MM}. 
\end{proof}

\subsection{Reduction for tiltings}

In this subsection, we study the tilting theory of MMAs and NCCRs. By constructing natural tilting modules, we prove that all MMAs in this setting are derived equivalent. Thus the noncommutative Bondal--Orlov conjecture holds in our setting. First, we recall some basic facts about tilting modules.

\begin{dfn}
   Let $\Lambda$ be a ring and let $T\in\modcat\Lambda$. We consider the following conditions.
\begin{enumerate}
\item $\pd_{\Lambda}T\le 1$.
\item $\Ext_{\Lambda}^1(T,T)=0$.
\item There exists an exact sequence
$0\to \Lambda\to T_1\to T_0\to 0$
such that $T_i\in\add T$ for $i=0,1$.
\end{enumerate}
The module $T$ is called a {\it tilting module} (resp. {\it partial tilting module}) if it satisfies conditions (1), (2), and (3) (resp. conditions (1) and (2)).
\end{dfn}

The following is an important observation.

\begin{prop}\label{tilting +-1}
    Let $\Lambda$ be a Noetherian $R$-algebra and let $T\in\modcat \Lambda$. Assume that $t_1,\dots,t_m$ is a regular sequence on both $\Lambda$ and $T$. Then the following statements hold.
    \begin{enumerate}
        \item $T$ is a partial tilting $\Lambda$-module if and only if  $\oT$ is a partial tilting $\oLam$-module.
        \item $T$ is a tilting $\Lambda$-module if and only if  $\oT$ is a tilting $\oLam$-module.
    \end{enumerate}
\end{prop}
\begin{proof}
    We may assume that $m=1$ and let $t:=t_1$. By Lemma \ref{quot by u proj}, $\pd_{\Lambda}T\le 1$ if and only if $\pd_{\oLam}\oT\le 1$. Moreover, there exists an exact sequence
    \[\Ext^1_{\Lambda}(T,T)\xrightarrow{t}\Ext^1_{\Lambda}(T,T)\to \Ext^1_{\Lambda}(T,\oT)=\Ext^1_{\oLam}(\oT,\oT)\to \Ext^2_{\Lambda}(T,T),\]
    induced by $0\to T\xrightarrow{t}T\to \oT\to 0$.
    
    (1) First, suppose that $T$ is a partial tilting $\Lambda$-module. By assumption, $\Ext^{i}_{\Lambda}(T,T)=0$ for $i\ge 1$, thus we have $\Ext^1_{\oLam}(\oT,\oT)=0$. Conversely, suppose that $\oT$ is a partial tilting $\oLam$-module. Then $\Ext^1_{\oLam}(\oT,\oT)=0$; hence $\Ext^1_{\Lambda}(T,T)=0$ by NAK.
    
    (2) We suppose that $T$ is a partial tilting $\Lambda$-module. Then $\oT$ is a partial tilting $\oLam$-module by (1). In addition, applying $\overline{(-)}$ to the exact sequence $0\to \Lambda\to T_1\to T_0\to 0$, we have an exact sequence $0\to \oLam\to \overline{T_1}\to \overline{T_0}\to 0$. Hence $\oT$ is a tilting $\oLam$-module. 

    Conversely, suppose that $\oT$ is a tilting $\oLam$-module. By (1), $\Ext^1_{\Lambda}(T,T)=0$. Thus $\Hom_{\Lambda}(T,T)\to \Hom_{\oLam}(\oT,\oT)$ is surjective. Therefore, there exists a lift 
    \[0\to \Ker f\to T_1\xrightarrow{f} T_0\to 0\]
    of $0\to \oLam\to \overline{T_1}\xrightarrow{\overline{f}} \overline{T_0}\to 0$. By Lemma \ref{quot by u proj}, we see that $\Ker f\cong\Lambda$.
\end{proof}

Next, we extend the Proposition \ref{MM resol} to the present setting.

\begin{thm}\label{MM resol general}
    If $N\in\refl R$ is a modifying module, then $\overline{\Hom_R(M,N)}\cong \Hom_{\oR}(\oM^{**},\oN^{**})$. In particular, $\pd_{\bLam_M}\Hom_R(M,N)\le 1$.
\end{thm}
\begin{proof}
  By Proposition \ref{MM resol preliminary}, there exists $N'\in\modif (R)$ such that $\overline{N'}^{**}\cong \oN^{**}$ and $\overline{\Hom_R(M,N')}\cong \Hom_{\oR}(\oM^{**},\overline{N'}^{**})$. By Corollary \ref{modif bij}, $N'\cong N$, we have $\overline{\Hom_R(M,N)}\cong \Hom_{\oR}(\oM^{**},\oN^{**})$. Therefore, $\pd_{\bLam_M}\Hom_R(M,N)\le 1$ follows from Proposition \ref{MM resol} and Lemma \ref{quot by u proj}.
\end{proof}

Recall that for $N\in\modif(R)$, we set
     \[\bLam_N:=\End_R(N),
        \quad \bLam_{\oN^{**}}:=\End_{\oR}(\oN^{**}),\]
    \[T_N:=\Hom_R(M,N),\quad T_{\oN^{**}}:=\Hom_{\oR}(\oM^{**},\oN^{**}).\]
    
The following is a main result of this subsection.

\begin{thm}\label{tilting of NCCR}
Let $N\in\refl (R)$. Then the following assertions hold.
\begin{enumerate}
    \item If $N\in\modif(R)$, then $T_N=\Hom_R(M,N)$ is a partial tilting $\bLam_M$-module.
    \item If $N\in\modif(R)$, then $T_{\oN^{**}}=\Hom_{\oR}(\oM^{**},\oN^{**})$ is a partial tilting $\bLam_{\oM^{**}}$-module.
    \item If $N\in\MM(R)$, then $T_N=\Hom_R(M,N)$ is a tilting $\bLam_M$-module.
    \item If $N\in\MM(R)$, then $T_{\oN^{**}}=\Hom_{\oR}(\oM^{**},\oN^{**})$ is a tilting $\bLam_{\oM^{**}}$-module.
\end{enumerate}
Moreover, if $N\in\MM(R)$, then there exists the following commutative diagram, in which all the vertical functors are equivalences of categories:
\[\begin{tikzcd}
	{\Db(\modcat\bLam_M)} && {\Db(\modcat\bLam_{\oM^{**}})} \\
	\\
	{\Db(\modcat\bLam_N)} && {\Db(\modcat\bLam_{\oN^{**}})}
	\arrow["{-\underset{R}{\Lo} \oR}", from=1-1, to=1-3]
	\arrow["{\RHom_{\bLam_M}(T_N,-)}"', from=1-1, to=3-1]
	\arrow["{\RHom_{\bLam_{\oM^{**}}}(T_{\oN^{**}},-)}", from=1-3, to=3-3]
	\arrow["{-\underset{R}{\Lo} \oR}"', from=3-1, to=3-3]
\end{tikzcd}\]
\end{thm}
\begin{proof}
(2) and (4) This follows from Theorem \ref{MM(R)->MM(R/t) MM} and \cite[Proposition 4.14]{Maximal-modifications-and-Auslander-Reiten-duality-for-non-isolated-singularities}.

(1) and (3) By Theorem \ref{MM resol general}, $\overline{T_N}\cong T_{\oN^{**}}$ for every $N\in\modif (R)$. Therefore, the claim follows from (2), (4) and Proposition \ref{tilting +-1}.

Finally, we show the commutativity of the diagram. This follows from the following computation:
\begin{align}
    \RHom_{\bLam_M}&(\Hom_R(M,N),-)\underset{R}{\Lo} \oR\\
&\cong \RHom_{\bLam_M}(\Hom_R(M,N),-\underset{R}{\Lo}\oR)
&& \text{Since $\Hom_R(M,N)$ is perfect}\\
&\cong \RHom_{\bLam_{\oM}}(\Hom_R(M,N)\underset{R}{\Lo} \oR,-\underset{R}{\Lo} \oR)
&& \text{By adjunction} \\
&\cong \RHom_{\bLam_{\oM}}(\Hom_{\oR}(\oM,\oN),-\underset{R}{\Lo} \oR)
&& \text{By Theorem \ref{MM resol general}}
\end{align}
\end{proof}

\begin{rem}
    In the above proof, we used Proposition \ref{tilting +-1}, but the result can also be proved by an argument similar to that of \cite{Maximal-modifications-and-Auslander-Reiten-duality-for-non-isolated-singularities}.
\end{rem}

\subsection{Reduction for two-term silting complexes} 
Throughout this subsection, we assume that $M\in\CM R$. In this subsection, we study the stable endomorphism algebras of $M$ and $\oM^{**}$. In particular, we establish an Adachi--Iyama--Reiten-type correspondence in our setting and show that it is compatible with the corresponding three-dimensional Adachi--Iyama--Reiten-type correspondence.

Recall that for $N\in\modifG(R)$, we set
     \[\bsL_N:=\stEnd_R(N)
        \quad\text{and}\quad \bsL_{\oN^{**}}:=\stEnd_{\oR}(\oN^{**}).\]
        By Theorem \ref{MM resol general} $\pd_{\bLam_M}\Hom_R(M,N)\le 1$. We consider the projective resolution
\[0\to \Hom_R(M,M_1^N)\to \Hom_R(M,M_0^N)\to\Hom_R(M,N)\to 0.\] 
Moreover, this projective resolution induces the following projective resolution of $\Hom_{\oR}(\oM,\oN)$:
\[\begin{tikzcd}
	0 & {\overline{\Hom_R(M,M_1^N)}} & {\overline{\Hom_R(M,M_0^N)}} & {\overline{\Hom_R(M,N)}} & 0 \\
	0 & {\Hom_{\oR}(\oM,\oM_1^{\oN})} & {\Hom_{\oR}(\oM,\oM_0^{\oN})} & {\Hom_{\oR}(\oM,\oN)} & 0.
	\arrow[from=1-1, to=1-2]
	\arrow[from=1-2, to=1-3]
	\arrow["\cong", from=1-2, to=2-2]
	\arrow[from=1-3, to=1-4]
	\arrow["\cong", from=1-3, to=2-3]
	\arrow[from=1-4, to=1-5]
	\arrow["\cong", from=1-4, to=2-4]
	\arrow[from=2-1, to=2-2]
	\arrow[from=2-2, to=2-3]
	\arrow[from=2-3, to=2-4]
	\arrow[from=2-4, to=2-5]
\end{tikzcd}\]
Let 
\[P_N:=T_N\underset{\bLam_M}{\Lo}\bsL_M\cong(\stHom_R(M,M_1^N)\to \stHom_R(M,M_0^N)),\]\
\[P_{\oN}:=T_{\oN}\underset{\bLam_{\oM}}{\Lo}\bsL_{\oM}\cong(\stHom_{\oR}(\oM,\oM_1^{\oN})\to \stHom_{\oR}(\oM,\oM_0^{\oN})).\]
Note that, $P_{\oN}=P_N\underset{\bsL_M}{\otimes}\bsL_{\oM}=P_N\underset{\bsL_M}{\Lo}\bsL_{\oM}$.

\begin{lem}[e.g. {\cite[Proposition 4.2]{Tilting-and-silting-theory-of-Noetherian-algebras}}]\label{2-silt /I}
    Let $\Lambda$ be a ring and let $I\subseteq\Lambda $ be a two-sided ideal. Then there exists a map:
\[\begin{tikzcd}
	{\twopresilt(\Lambda)} && {\twopresilt(\Lambda/I)} \\
	{\twosilt(\Lambda)} && {\twosilt(\Lambda/I)}
	\arrow["{-\underset{\Lambda}{\otimes}\Lambda/I}", from=1-1, to=1-3]
	\arrow[hook', from=2-1, to=1-1]
	\arrow["{-\underset{\Lambda}{\otimes}\Lambda/I}", from=2-1, to=2-3]
	\arrow[hook', from=2-3, to=1-3]
\end{tikzcd}\]
\end{lem}

Although an assumption on the ideal is imposed in \cite{Tilting-and-silting-theory-of-Noetherian-algebras}, this assumption is not used in the proof of the above lemma.

\begin{thm}\label{2-siltA->2-siltA/t}
   There exists the  following commutative diagram compatible with mutation.

    \[
    \begin{tikzcd}
	{\MMG(R)} && {\MMG(\overline{R})} \\
	{\twosilt(\bsL_M)} && {\twosilt(\bsL_{\oM})}
	\arrow["\sim", from=1-1, to=1-3]
	\arrow["\sim", from=1-1, to=2-1]
	\arrow["\sim", from=1-3, to=2-3]
	\arrow["\sim", from=2-1, to=2-3]
\end{tikzcd}
\qquad
    \begin{tikzcd}
	{\modifG(R)} && {\modifG(\oR)} \\
	{\twopresilt(\bsL_M)} && {\twopresilt(\bsL_{\oM})}
	\arrow["\sim", from=1-1, to=1-3]
	\arrow["\sim", from=1-1, to=2-1]
	\arrow["\sim", from=1-3, to=2-3]
	\arrow["\sim", from=2-1, to=2-3]
\end{tikzcd}\]
The correspondence is given by the following:
\[\begin{tikzcd}
	N && \oN \\
	{P_N} && {\overline{P_N}\cong P_{\oN}}
	\arrow[maps to, from=1-1, to=1-3]
	\arrow[maps to, from=1-1, to=2-1]
	\arrow[maps to, from=1-3, to=2-3]
	\arrow[maps to, from=2-1, to=2-3]
\end{tikzcd}\]
\end{thm}
\begin{proof}
First, we show that $P_N$ and $P_{\oN}$ are silting (resp. presilting) for $N\in\MMG(R)$ (resp. $N\in\modifG(R)$). By Theorem \ref{tilting of NCCR}, $T_N$ and $T_{\oN}$ are tilting modules (resp. partial tilting modules). Thus
\[
(\Hom_R(M,M_1^N)\to \Hom_R(M,M_0^N)) \quad\text{and}\quad
(\Hom_{\oR}(\oM,\oM_1^{\oN})\to \Hom_{\oR}(\oM,\oM_0^{\oN}))
\]
are two-term silting complexes (resp. two-term presilting complexes), since they are quasi-isomorphic to $T_N$ and $T_{\oN}$, respectively. Therefore, the claim follows from Lemma \ref{2-silt /I}.

Next, we show  bijectivity. 
\begin{enumerate}
    \item By Theorem \ref{MM(R)->MM(R/t) MM}, $\MMG(R)\to \MMG(\oR)$ and $\modifG(R)\to \modifG(\oR)$ are bijective compatible with mutation.
    \item By \cite{tau-tilting-theory} and \cite[Theorem 2.18]{On-the-finiteness-of-the-derived-equivalence-classes-of-some-stable-endomorphism-rings}, $\MMG(\oR)\to \twosilt(\bsL_{\oM})$ and  $\modifG(\oR)\to \twopresilt(\bsL_{\oM})$ are bijective compatible with mutation.
    \item By \cite[Theorem 4.3]{Tilting-and-silting-theory-of-Noetherian-algebras} and \cite[Lemma 2.7]{Stability-over-cDV-singularities-and-other-complete-local-rings}, $\twosilt(\bsL_{M})\to \twosilt(\bsL_{\oM})$ and $\twopresilt(\bsL_{M})\to \twopresilt(\bsL_{\oM})$ are bijective compatible with mutation.
\end{enumerate}
Therefore, all the maps are  bijective compatible with mutation. 
\end{proof}

\subsection{Some remarks}

Finally, we examine the assumptions in Setup \ref{intro setup  reduction NCCR,tilt,2-silt}. In the setup, we assume that $M$ is maximal modifying. However, by Theorem \ref{NCCR +1}, it suffices to assume that $M$ is modifying. Thus the assumption in the setup amounts to the existence of a modifying module that lifts a maximal modifying module over $\oR$.

Nevertheless, requiring $M$ to be modifying is still a somewhat strong assumption. In this subsection, we prove that, when $M$ is Cohen--Macaulay, the modifying property of $M$ follows from that of $\oM$. Although this phenomenon is specific to the Cohen--Macaulay case, it is sufficient for practical applications.

\begin{lem}\label{HomCM+1}
    Let $L,N\in\CM R$. If $\Hom_{\overline{R}}(\overline{L},\overline{N})\in\CM\overline{R}$, then $\Hom_R(L,N)\in\CM R$.
\end{lem}
\begin{proof}
We prove only the case $m=1$. The general case follows by repeating the same argument. Let $t:=t_1$.

We consider the exact sequence
\[\Ext_R^1(L,N)\xrightarrow{t}\Ext_R^1(L,N)\to \Ext_{\overline{R}}^1(\overline{L},\overline{N}).\]
Since $\oR$ is a 3-dimensional isolated singularity and $\Hom_{\oR}(\oL,\oN)\in\CM\oR$, we have $\Ext_{\overline{R}}^1(\overline{L},\overline{N})=0$. Thus $\Ext_R^1(L,N)\xrightarrow{t}\Ext_R^1(L,N)$ is surjective; hence $\Ext_R^1(L,N)=0$ by NAK. Thus we have $\overline{\Hom_R(L,N)}\cong\Hom_{\oR}(\oL,\oN)$. Since $\Hom_{\overline{R}}(\overline{L},\overline{N})\in\CM \overline{R}$, we obtain $\Hom_R(L,N)\in\CM R$.
\end{proof}

Note that there is a counterexample to a natural generalization of this lemma.

\begin{rem}
    The following statement does not hold in general, as a counterexample shows.
    \begin{enumerate}
        \item Let $L,N\in\refl R$. If $\Hom_{\overline{R}}(\overline{L}^{**},\overline{N}^{**})\in\CM\overline{R}$, then $\Hom_R(L,N)\in\CM R$.
        \item Lemma \ref{HomCM+1} under the sole assumption that $\oR$ is 3-dimensional Gorenstein normal domain.  
    \end{enumerate}
\end{rem}
\begin{proof}
(1) Let 
\[R:=k[[x,y,z]], \quad\oR:=k[[x,y]]\quad\text{and}\quad L=N:=\Omega_R^2 k.\]

Since $\oR$ is 2-dimensional regular local ring, $\oL^{**}=\oN^{**}\in\proj \oR$; thus $\Hom_{\overline{R}}(\overline{L}^{**},\overline{N}^{**})\in\CM\overline{R}$. Suppose that $\Hom_R(L,N)\in\CM R$. Then $L$ is projective by Proposition \ref{modif RLR}, a contradiction. Thus $\Hom_R(L,N)\notin\CM R$.

(2) Let 
    \[R:=k[[x,y,z,w,u]]/(xy+zw)\quad\text{and}\quad \oR:=k[[x,y,z,u]]/(xy+z^2)\]
    and let 
    \[L:=(x,z)\quad\text{and}\quad N:=(x,w).\]
    
    Since $\oL=\oN=(x,z)$, we have $\Hom_{\oR}(\oL,\oN)\in\CM\oR$. But, $\Hom_R(L,N)\notin\CM R$. Indeed, $R$ is a polynomial extension of the toric ring associated with the weights $(1,1,-1,-1)$, and its divisor class group is isomorphic to $\bZ$. In this case, $L$ and $N$ correspond to the classes $1$ and $-1$ in the divisor class group, respectively. Moreover, it is known that, for $i\in\bZ$, the corresponding rank-one reflexive module is Cohen--Macaulay if and only if $i\in\{-1,0,1\}$. Therefore, $\Hom_R(L,N)\notin\CM R$.
\end{proof}
\begin{cor}\label{NCCR +1 aut}
    If there exists $M\in\CM R$ such that $\overline{M}\in\CM \oR$ gives an NCCR (resp. is an MM $\oR$-module), then $M$ gives an NCCR of $R$ (resp. is an MM $R$-module.)
\end{cor}
\begin{proof}
    This follows from Lemma \ref{HomCM+1} and Corollary \ref{CM NCCR +1}.
\end{proof}

% \begin{cor}\label{add =}
%     Let $L,N\in\CM R$ such that $\oL$ and $\oN$ are Cohen--Macaulay modifying modules. Then 
%     \[L\in\add N\Leftrightarrow \oL\in\add\oN.\]
% \end{cor}
% \begin{proof}
%     The implication $\Rightarrow$ is clear. We show the converse. By \ref{HomCM+1}, $L$ and $N$ is modifying. Therefore, the claim follows from \ref{HomCM+1} and \ref{quot indec isom}.
% \end{proof}

\section{Reduction for stability conditions via derived base change}
\subsection{Reduction theorem}

In this subsection, we prove that suitable subspaces of the spaces of stability conditions on derived categories are isomorphic under derived base change. Although this result is interesting, it follows almost directly from the Koenig--Yang correspondence and Fushimi's reduction theorem.

\begin{dfn}[\cite{Silting-theory-and-derived-base-change}]\label{dg def}
    Let $R$ and $S$ be commutative non-positive dg rings. 
    \begin{enumerate}
        \item $R$ is called {\it complete local Noetherian} if $H^0(R)$ is complete local Noetherian ring.
        \item A morphism $R\to S$ of complete local Noetherian non-positive commutative dg rings is called {\it local} if the induced ring morphism $H^0(R)\to H^0(S)$ is local.
    \end{enumerate}
\end{dfn}

\begin{setup}
    In this subsection, we work under Setup \ref{intro setup reduction stab}.
\end{setup}

First, we recall some known results on the silting base change theorem.

\begin{prop}[{\cite[Theorem B]{Silting-theory-and-derived-base-change}}]\label{silt derived base change}
    The functor $\overline{(-)}=-\underset{R}{\Lo} k$ induces a poset isomoprphism
    \[\overline{(-)}:(\silt(\per A),\le)\xrightarrow[]{\sim}(\silt(\per \oA),\le).\]
    Here $\per A:=\thick_{\D(A)}A_A$ and $\per \oA:=\thick_{\D(\oA)}\oA_{\oA}$  is a perfect derived category.  
\end{prop}

We next recall some known results on the  Koenig--Yang correspondence in this setting.

\begin{prop}[{\cite{Silting-objects-simple-mindedcollections-t-structures-and-co-t-structures-for-finite-dimensional-algebras,Silting-theory-and-derived-base-change}}]\label{KY corr}
    There is a bijection between the following sets:
    \begin{enumerate}
        \item isomorphism class of basic silting objects in $\per A$,
        \item isomorphism class of basic simple minded collections in $\pvd A$,
        \item bounded co-$t$-structures on $\per A$,
        \item algebraic $t$-structures on $\pvd A$.
    \end{enumerate}
\end{prop}
\begin{proof}
(1)$\leftrightarrow$(2) follows from \cite[Theorem A]{Silting-theory-and-derived-base-change}. 
(1)$\leftrightarrow$(3) follows from \cite[Corollary 5.9]{Auslander-Buchweitz-context-and-co-t-structures}.
(2)$\leftrightarrow$(4) follows from \cite{Silting-objects-simple-mindedcollections-t-structures-and-co-t-structures-for-finite-dimensional-algebras}.   
\end{proof}

For a basic silting object $S\in\per A$, we denote the corresponding algebraic $t$-structure on $\pvd A$ by $(\scrD^{\le0}_S,\scrD^{\ge 0}_S)$. Moreover, we set $\scrH_S:=\scrD^{\le0}_S\cap\scrD^{\ge 0}_S\simeq \fl \End_{A}(S)$. Note that, by Proposition \ref{KY corr}, we have $\Stab^{\alg}(\pvd A)=\{(\scrH_S,Z)\in\Stab(\pvd A)\mid S\in\silt (\per A)\}.$

Let $\iota:\pvd \oA\to \pvd A$ be the forgetful functor. This functor induces an isomorphism of Grothendieck groups $\iota_{\ast}:K_0(\pvd \oA)\xrightarrow[]{\sim}K_0(\pvd A)$. For $Z\in\Hom_{\bZ}(K_0(\pvd A),\bC)$, we set $\overline{Z}:=Z\circ\iota_{\ast}$.

\begin{thm}\label{Stab base change}
     There exists an isomorphism between the spaces of algebraic stability conditions:
    \[\overline{(-)}:\Stab^{\alg}(\pvd A)\xrightarrow[]{\sim}\Stab^{\alg}(\pvd \oA)\quad (\scrH_S,Z)\mapsto (\scrH_{\oS},\overline{Z}). \]
\end{thm}
\begin{proof}
    By Proposition \ref{silt derived base change}, $\overline{(-)}:\Stab^{\alg}(\pvd A)\xrightarrow[]{}\Stab^{\alg}(\pvd \oA)$ is bijective. Therefore, it suffices to show that $\overline{(-)}$ is an isometry.

    Let $S,T\in\silt(\per A)$  and let $\sigma=(\scrH_S,Z)=(\mathcal{P}_S,Z),\tau=(\scrH_T,W)=(\mathcal{P}_T,W)\in\Stab^{\alg} (\pvd A)$ where $\mathcal{P}_S,\mathcal{P}_T$ are corresponding slicing of $\pvd A$. Moreover, slicing of $\pvd A$ corresponding to $\overline{\sigma},\overline{\tau}$ denote by $\mathcal{P}_{\oS},\mathcal{P}_{\oT}$.

\medskip{\it step1: $\mathcal{P}_S(I)=\Filt(\iota(\mathcal{P}_{\oS}(I)))$ for every interval $I\subseteq\bR$, where $\Filt(\iota(\mathcal{P}_{\oS}(I)))$ is the extension closure in $\mathcal{P}_S(I)$ of $\iota(\mathcal{P}_{\oS}(I))$.}
It suffices to prove the assertion in the case $I=\{\phi\}$ where $\phi\in\bR$. A inclusion $\supseteq$ is clear. We show $\subseteq$. Since $\scrH_S$ is length, every object in $\mathcal{P}_S(\phi)$ is filtered by $\sigma$-stable objects in $\mathcal{P}_S(\phi)$. Thus it suffices to show that every $\sigma$-stable object in $\mathcal{P}_S(\phi)$ belong to $\Filt(\iota(\mathcal{P}_{\oS}(\phi)))$. Every $\sigma$-stable object is brick, and the actions of the element of $\mm$ lie in the radical, hence are zero. Therefore, the claim follows.

\medskip{\it step2: $d(\mathcal{P}_S,\mathcal{P}_T)=d(\mathcal{P}_{\oS},\mathcal{P}_{\oT})$.} We recall that 
\[d(\mathcal{P}_S,\mathcal{P}_T):=\sup_{0\ne E\in \pvd A}\{|\phi_{\mathcal{P}_S}^+(E)- \phi_{\mathcal{P}_{T}}^+(E)|,|\phi_{\mathcal{P}_S}^-(E)- \phi_{\mathcal{P}_{T}}^-(E)|\},\]
\[d(\mathcal{P}_{\oS},\mathcal{P}_{\oT}):=\sup_{0\ne F\in \pvd \oA}\{|\phi_{\mathcal{P}_{\oS}}^+(F)- \phi_{\mathcal{P}_{\oT}}^+(F)|,|\phi_{\mathcal{P}_{\oS}}^-(F)- \phi_{\mathcal{P}_{\oT}}^-(F)|\}.\]

By step1, the functor $\iota$ sends the Harder--Narasimhan filtration of an object with respect to $\overline{\sigma}$ to its Harder--Narasimhan filtration with respect to $\sigma$. Hence, for any non-zero object $F\in\pvd \oA$, we have $\phi_{\sigma}^{\pm}(\iota(F))=\phi_{\overline{\sigma}}^{\pm}(F)$. The same statement holds for $\tau$ and $\overline{\tau}$. Therefore, by the definition of the distance between slicings, we obtain $d(\mathcal{P}_{\oS},\mathcal{P}_{\oT})\le d(\mathcal{P}_S,\mathcal{P}_T)$. We show the converse inequality. Recall that, by \cite[Lemma 6.1]{Stability-conditions-on-triangulated-categories}, we have
\[d(\mathcal{P}_S,\mathcal{P}_T)=\inf \{\varepsilon\in\bR_{\ge 0}\mid \text{$\mathcal{P}_T(\phi)\subset \mathcal{P}_S([\phi-\varepsilon,\phi+\varepsilon])$ for all $\phi\in\bR$}\},\]
\[d(\mathcal{P}_{\oS},\mathcal{P}_{\oT})=\inf \{\varepsilon\in\bR_{\ge 0}\mid \text{$\mathcal{P_{\oT}}(\phi)\subset \mathcal{P}_{\oS}([\phi-\varepsilon,\phi+\varepsilon])$ for all $\phi\in\bR$}\}.\]
Thus, the inequality $d(\mathcal{P}_{\oS},\mathcal{P}_{\oT})\ge d(\mathcal{P}_S,\mathcal{P}_T)$ follows from step1.

\medskip{\it step3: $||U||_{\sigma}=||\overline{U}||_{\overline{\sigma}}$ for every $U\in\Hom_{\bZ}(K_0(\pvd A),\bC)$.}
We recall that 
\[\text{$||U||_{\sigma}:=\sup\{\frac{|U(E)|}{|Z(E)|}\mid \text{$E$ is $\sigma$-semistable}\}$ and $||\overline{U}||_{\overline{\sigma}}:=\sup\{\frac{|\overline{U}(F)|}{|\overline{Z}(F)|}\mid \text{$F$ is $\overline{\sigma}$-semistable}\}.$}\]
By step1, if $F$ is $\overline{\sigma}$-semistable, then $\iota(F)$ is also $\sigma$-semistable. Hence, the inequality $||U||_{\sigma}\ge||\overline{U}||_{\overline{\sigma}}$ follows. We now prove the converse. Foe $E\in\mathcal{P}(\phi)=\langle\iota(\overline{\mathcal{P}}(\phi))\rangle$, there exists $F_1,\dots,F_s\in\overline{\mathcal{P}}(\phi)$ such that $[E]=\sum\limits_{i=1}^{s}[\iota(F_i)]$. Thus $|U(E)|\le\sum\limits_{i=1}^{s}|\overline{U}(F_i)|$. Moreover, since $\overline{Z}(F_1),\dots,\overline{Z}(F_s)$ lie on a single line, $|Z(E)|=\sum\limits_{i=1}^s|\overline{Z}(F_i)|$. Therefore $\frac{|U(E)|}{|Z(E)|}\le\frac{\sum\limits_{i=1}^s |\overline{U}(F_i)|}{\sum\limits_{i=1}^s |\overline{Z}(F_i)|}\le ||\overline{U}||_{\overline{\sigma}}$, we obtain the inequality $||U||_{\sigma}\le||\overline{U}||_{\overline{\sigma}}$.
\end{proof}
\begin{cor}\label{Stab general base change}
    Let $(R,\mm,k)\to (S,\nn,k)$ be a local morphism of complete Noetherian local rings. Then there exists a natural isomorphism between the spaces of algebraic stability conditions:
    \[\Stab^{\alg}(\pvd A)\xrightarrow[]{\sim}\Stab^{\alg}(\pvd A\;\underset{R}{\Lo}\; S). \]
\end{cor}
\subsection{Corollaries}
\begin{setup}
    In this subsection, we work under Setup \ref{intro setup  reduction NCCR,tilt,2-silt}.
\end{setup}

In this subsection, we apply Corollary \ref{Stab general base change} to the spaces of stability conditions on the derived categories associated with MMAs.

The following lemma is particularly interesting.

\begin{lem}\label{regular sequence}
    If $N\in\CM R$ is a Cohen--Macaulay modifying $R$-module, then $t_1,\dots,t_m$ is a regular sequence on $\stHom_R(M,N)$.
\end{lem}
\begin{proof}
Suppose that $m=1$ and let $t:=t_1$. We use the equivalence $\stCM R\simeq \sg R:=\Db(\modcat R)/\Kb(\proj R)$ (e.g. \cite{Quotients-of-triangulated-categories-and-equivalences-of-Buchweitz-Orlov-and-Amiot-Guo-Keller}). 

First, note that two functors $\RHom_{\oR}(\oR.-)$ and $-\underset{R}{\Lo}\oR$ preserved $\D^b (\modcat)$ and $\Kb(\proj)$:
\[
\begin{tikzcd}
	{\Db(\modcat R)} && {\Db(\modcat \oR)} \\
	{\Kb(\proj R)} && {\Kb(\proj \oR).}
	\arrow["{-\underset{R}{\Lo}\oR}"', from=1-1, to=1-3]
	\arrow[hook, from=2-1, to=1-1]
	\arrow[hook, from=2-3, to=1-3]
	\arrow["{-\underset{R}{\Lo}\oR}"', from=2-1, to=2-3]
\end{tikzcd}
\quad
\begin{tikzcd}
	{\Db(\modcat R)} && {\Db(\modcat \oR)} \\
	{\Kb(\proj R)} && {\Kb(\proj \oR)}
	\arrow["{\RHom_{\oR}(\oR.-)}", from=1-3, to=1-1]
	\arrow[hook, from=2-1, to=1-1]
	\arrow["{\RHom_{\oR}(\oR.-)}", from=2-3, to=2-1]
	\arrow[hook, from=2-3, to=1-3]
\end{tikzcd}
\]
Therefore, these functors induce functors between the singularity categories, and the induced functors form an adjoint pair:
\[\begin{tikzcd}
	{\sg R} && {\sg\oR}
	\arrow[""{name=0, anchor=center, inner sep=0}, "{-\underset{R}{\Lo}\oR}", shift left=2, from=1-1, to=1-3]
	\arrow[""{name=1, anchor=center, inner sep=0}, "{\RHom_{\oR}(\oR.-)}", shift left=2, from=1-3, to=1-1]
	\arrow["\dashv"{anchor=center, rotate=-90}, draw=none, from=0, to=1]
\end{tikzcd}\]

Next, we consider an exact sequence 
\[\Hom_{\sg R}(M,N[-1])\to \Hom_{\sg R}(M,\oN[-1])=\Hom_{\sg \oR}(\oM,\oN[-1])\to \Hom_{\sg R}(M,N)\xrightarrow{t}\Hom_{\sg R} (M,N)\]
induced by $0\to N\to N\to \oN\to 0$. Since $\Omega_{R} N$ is modifying by \cite[Theorem 4.3]{Maximal-modifications-and-Auslander-Reiten-duality-for-non-isolated-singularities}, we can apply Theorem \ref{MM resol general} to conclude that $\Hom_{R}(M,\Omega_R N)\to \Hom_{\oR}(\oM,\Omega_{\oR}\oN)$ is surjective. Thus $\Hom_{\sg R}(M,N[-1])\to \Hom_{\sg \oR}(\oM,\oN[-1])$ is surjective; hence, we have $t \in\reg \stHom_R(M,N)$. The desired conclusion follows by repeatedly applying the above argument.
\end{proof}

As a corollary, we obtain the following. 

% \begin{cor}
%     $\bsL_M=\stEnd_R(M)$ is an order over $k[[x_1,\cdots ,x_m]]$, where $x_i$ acts as multiplication $t_i$.
% \end{cor}

\begin{cor}
    $\bsL_M\underset{R}{\Lo}\oR=\bsL_M\underset{R}{\otimes}\oR=\bsL_{\oM}$ holds.
\end{cor}

\begin{cor}\label{stEnd Stab reduction}
    There exists a following isomorphisms:
    \begin{enumerate}
        \item $\overline{(-)}:(\silt(\per \bsL_M),\le)\xrightarrow[]{\sim}(\silt(\per \bsL_{\oM}),\le)$ is an isomoprphism of posets.
        \item $\overline{(-)}:\Stab^{\alg} \Db(\fl \bLam_M)\xrightarrow{\cong} \Stab^{\alg}\Db(\fl \bLam_{\oM^{**}})$ is an isomoprphism of topological spaces.
        \item If $M\in\CM R$, then $\overline{(-)}:\Stab^{\alg}(\Db(\fl\bsL_{M}))\xrightarrow[]{\sim}\Stab^{\alg}(\Db(\fl \bsL_{\oM}))$ is an isomoprphism of topological spaces.
    \end{enumerate}
\end{cor}
\begin{proof}
    This follows from Proposition \ref{silt derived base change} and Corollary \ref{Stab general base change}.
\end{proof}

\section{t-structures for silting discrete Noetherian algebras}

\begin{setup}
    In this subsection, we work under Setup \ref{intro setup reduction stab}.
\end{setup}

In this section, we extend \cite[Theorem A]{Contractibility-of-the-stability-manifold-for-silting-discrete-algebras} to the setting of Noetherian algebras. We say that $A$ is {\it silting discrete} if $\per A$ is silting discrete triangulated category (see \cite{Classifying-tilting-complexes-over-preprojective-algebras-of-Dynkin-type}).

For a basic silting object $S\in\per A$, we denote the corresponding algebraic $t$-structure on $\pvd A$ by $(\scrD^{\le0}_S,\scrD^{\ge 0}_S)$. That is, 
\[\scrD_S^{\le 0}:=\{X\in\pvd A\mid \Hom_A(S,X[>0]=0)\}\quad\text{and}\quad
\scrD_S^{\ge 0}:=\{X\in\pvd A\mid \Hom_A(S,X[<0]=0)\}.\]

Next, we prepare same facts needed for the proof.

\begin{lem}
For $S,T\in\silt(\per A)$, $S\le T$ if and only if
$\scrD_S^{\le 0}\subseteq\scrD_T^{\le 0}$.
\end{lem}
\begin{proof}
Let $B:=\RHom_A(S,S)$.
Since $S$ is silting, $B$ is
non-positive. By derived Morita theory, there is
an equivalence
\[
\RHom_A(S,-):\D(A)\xrightarrow{\sim}\D(B).
\]
Moreover, this equivalence restricts to an equivalence
\[
\RHom_A(S,-):\pvd A\xrightarrow{\sim}\pvd B.
\]
Thus, we may assume that $S=A$ without loss of
generality.

First, suppose that $A\le T$. For $X\in\scrD_A^{\le 0}$, $X\underset{A}{\Lo}\RHom_A(T,A)$ concentrated in degree $\le 0$ by assumption. Thus $\RHom_A(T,X)\cong X\underset{A}{\Lo}\RHom_A(T,A)$ concentrated in degree $\le 0$, we obtain $X\in\scrD_T^{\le 0}$

Conversely, suppose that
$\scrD_A^{\le0}\subseteq\scrD_T^{\le0}$. Then $\topmod H^0 (A)\in\scrD_A^{\le0}$, and hence $\topmod H^0(A)\underset{A}{\Lo}\RHom_A(T,A)\cong \RHom_A(T,\topmod H^0(A))$ concentrated in degree $\le0$. Therefore, $\RHom_A(T,\topmod H^0(A))$ concentrated in degree $\le0$, and hence $A\le T$.
\end{proof}

\begin{cor}\label{IJY result}
    For $S\in\silt(\per A)$, the bijection between (1) and (4) in Proposition \ref{KY corr} restrict to the following bijection:
\[
\begin{array}{c@{\qquad}c@{\qquad}c}
\left\{
\text{$T\in\silt(\per A)\mid S[1]\le T\le S$}\right\}
&
\overset{1\text{-}1}{\longleftrightarrow}
&
\left\{
\begin{array}{c}
\text{algebraic }t\text{-structures $(\scrD^{\le 0},\scrD^{\ge 0})$ in $\pvd A$}\\
\text{such that }\scrD_S^{\le-1}\subseteq\scrD^{\le 0}\subseteq\scrD_S^{\le0}
\end{array}
\right\}
\end{array}
\]
\end{cor}

\begin{lem}[{e.g. \cite[Proposition 2.3]{Stability-conditions-torsion-theories-and-tilting}}]\label{HRS corr}
    Let $\scrT$ be a triangulated category and let $(\scrX,\scrY)$ be a $t$-structure with heart $\scrH$. Then there exists the following bijection:

\[
\begin{array}{c@{\qquad}c@{\qquad}c}
\left\{
t\text{-structure }(\scrU,\scrV)
\text{ with }\scrX[1]\subseteq\scrU\subseteq\scrX
\right\}
&
\overset{1\text{-}1}{\longleftrightarrow}
&
\left\{
\text{torsion pair }(\scrT,\scrF)\text{ in }\scrH
\right\}
\\[6pt]
(\scrU,\scrV)
&
\longrightarrow
&
(\scrH\cap\scrU,\scrH\cap\scrV[-1])
\\[6pt]
(\scrX[1]\ast\scrT,\scrF[1]\ast\scrY)
&
\longleftarrow
&
(\scrT,\scrF).
\end{array}
\]
\end{lem}

We say that a $t$-structure $(\scrU,\scrV)$ is {\it intermediate with respect to $(\scrX,\scrY)$} if $\scrX[1]\subseteq\scrU\subseteq\scrX$.

\begin{prop}[{\cite[Theorem 4.3 and Theorem 5.4]{Tilting-and-silting-theory-of-Noetherian-algebras}}]\label{Kimura result}
   For Noetherian algebra $\Lambda$ over complete Noetherian local ring, there exists the following diagram:
\[\begin{tikzcd}
	{\tors(\fl \Lambda)} && {\tors(\overline{\Lambda})} \\
	{\ftors(\Lambda)} && {\ftors(\overline{\Lambda})} \\
	{\msilt(\Lambda)} && {\msilt(\overline{\Lambda})} \\
	{\twosilt(\Lambda)} && {\twosilt(\overline{\Lambda})}
	\arrow["\sim", from=1-1, to=1-3]
	\arrow["{(-)\cap \modcat\overline{\Lambda}}"', from=1-1, to=1-3]
	\arrow["{\cap\fl \Lambda}", hook, from=2-1, to=1-1]
	\arrow["\sim", from=2-1, to=2-3]
	\arrow["{(-)\cap \modcat\overline{\Lambda}}"', from=2-1, to=2-3]
	\arrow[hook, from=2-3, to=1-3]
	\arrow["\sim"', from=3-1, to=2-1]
	\arrow["{\Fac (-)}", from=3-1, to=2-1]
	\arrow["\sim", from=3-1, to=3-3]
	\arrow["{M\mapsto M/\mm M}"', from=3-1, to=3-3]
	\arrow["\sim", from=3-3, to=2-3]
	\arrow["{\Fac (-)}"', from=3-3, to=2-3]
	\arrow["\sim"', from=4-1, to=3-1]
	\arrow["{H^0(-)}", from=4-1, to=3-1]
	\arrow["\sim", from=4-1, to=4-3]
	\arrow["{P\mapsto P\underset{\Lambda}{\Lo}\overline{\Lambda}}"', from=4-1, to=4-3]
	\arrow["\sim", from=4-3, to=3-3]
	\arrow["{H^0(-)}"', from=4-3, to=3-3]
\end{tikzcd}\]
Here, $\mm$ is the maximal ideal of the base ring and $\oLam:=\Lambda/\mm\Lambda$.
\end{prop}

We now prove the main theorem. The proof is largely parallel to that of \cite{Contractibility-of-the-stability-manifold-for-silting-discrete-algebras}, except that at certain points we use the above theorem to reduce to the finite dimensional case.

\begin{lem}\label{comm}
    Let $S\in\per A$ be a basic silting object, $\Lambda:=\End_{A}(S)$ and let $\scrH_S:=\scrD^{\le0}_S\cap\scrD^{\ge 0}_S$. Then there exists the following commutative diagram:
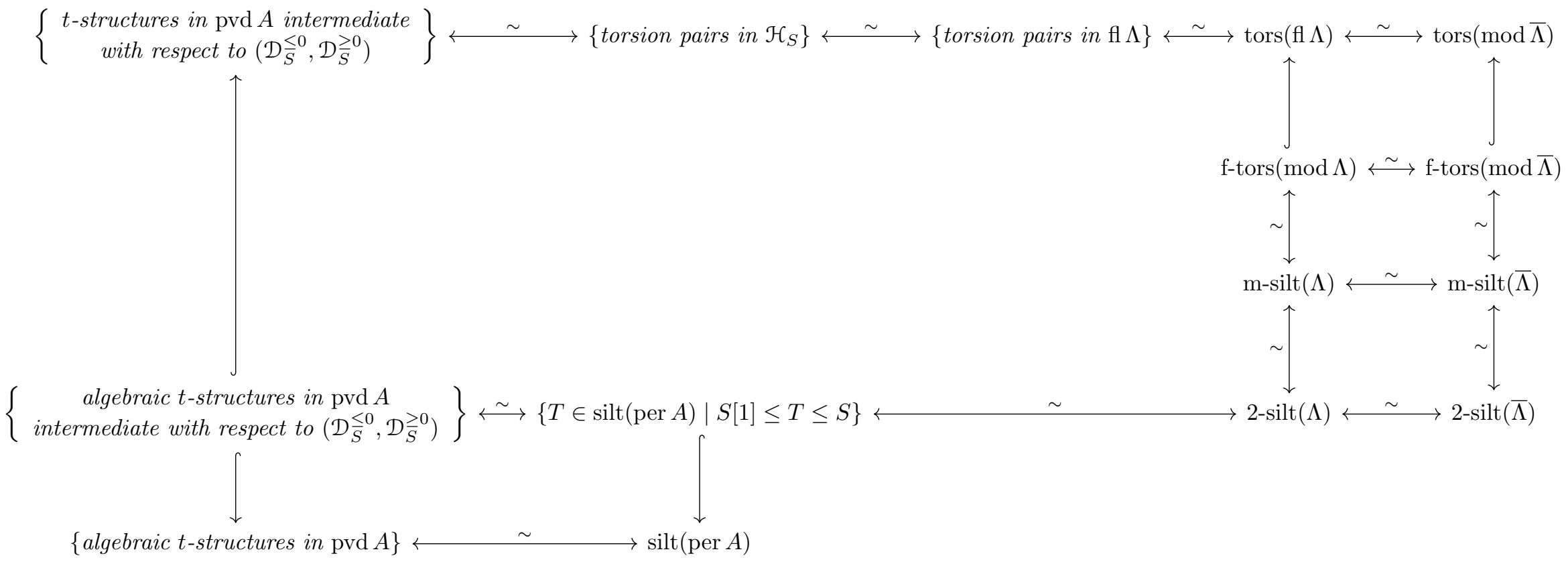
\begin{center}
\resizebox{\linewidth}{!}{%
\begin{tikzcd}[
  ampersand replacement=\&,
  column sep=2em,
  row sep=3em
]
\left\{
\begin{array}{c}
t\text{-structures in $\pvd A$ intermediate}\\
\text{with respect to }(\scrD_S^{\leq0},\scrD_S^{\geq0})
\end{array}
\right\}
\arrow[r, leftrightarrow, "\sim"]
\&
\left\{
\text{torsion pairs in $\scrH_S$}
\right\}
\arrow[r, leftrightarrow, "\sim"]
\&
\left\{
\text{torsion pairs in $\fl\Lambda$}
\right\}
\arrow[r, leftrightarrow, "\sim"]
\&
\tors(\fl\Lambda)
\arrow[r, leftrightarrow, "\sim"]
\&
\tors(\modcat\overline\Lambda)
\\
\&
\&
\&
\ftors(\modcat\Lambda)
\arrow[u, hook]
\arrow[r, leftrightarrow, "\sim"]
\&
\ftors(\modcat\overline{\Lambda})
\arrow[u, hook]
\\
\&
\&
\&
\msilt(\Lambda)
\arrow[u, leftrightarrow, "\sim"]
\arrow[r, leftrightarrow, "\sim"]
\&
\msilt(\overline{\Lambda})
\arrow[u, leftrightarrow, "\sim"]
\\
\left\{
\begin{array}{c}
\text{algebraic }t\text{-structures in $\pvd A$}\\
\text{intermediate with respect to }(\scrD_S^{\leq0},\scrD_S^{\geq0})
\end{array}
\right\}
\arrow[uuu, hook]
\arrow[d, hook]
\arrow[r, leftrightarrow, "\sim"]
\&
\left\{
T\in\silt(\per A)
\mid
S[1]\leq T\leq S
\right\}
\arrow[d, hook]
\arrow[rr, leftrightarrow, "\sim"]
\&
\&
\twosilt(\Lambda)
\arrow[u, leftrightarrow, "\sim"]
\arrow[r, leftrightarrow, "\sim"]
\&
\twosilt(\overline{\Lambda})
\arrow[u, leftrightarrow, "\sim"]
\\
\left\{
\text{algebraic }t\text{-structures in $\pvd A$}
\right\}
\arrow[r, leftrightarrow, "\sim"]
\&
\silt(\per A)
\&
\&
\&
\end{tikzcd}%
}
\end{center}
\end{lem}
\begin{proof}
By the natural equivalence $\scrH_S\simeq \fl\Lambda$, we obtain the bijection 
\[\{\text{torsion pairs in $\scrH_S$}\}\leftrightarrow\{\text{torsion pairs in $\fl\Lambda$}\}.\] The other bijections follow from the above facts.

We now check the commutativity. The lower-left square commutes by Corollary \ref{IJY result}, and the three squares on the right commute by Proposition \ref{Kimura result}. Finally, we check the commutativity of the large square in the middle. We take $T\in\silt(\per A)$ with $S_1\xrightarrow{d}S_0\to T\to\quad(S_1,S_0\in\add S)$. $\scrD_T^{\le 0}\cap\scrH_S=\Fac H^0T\cap\fl\Lambda$ under equivalence $\Phi:\scrH_S\xrightarrow{\sim} \fl\Lambda$ is what we need to prove. For $N\in\scrH_S$, 
there exists the following exact sequence:
\[\Hom_{A}(S_0,N)\xrightarrow{\Hom_{A}(d,N)} \Hom_{A}(S_1,N)\to\Hom_{A}(T,N[1])\to \Hom_{A}(S_0,N[1])=0. \]
Thus
\begin{align}
    N\in\scrD_T^{\le0}
    &
    \Leftrightarrow
    \Hom_{A}(T,N[>0])=0\\
    &
    \Leftrightarrow
    \Hom_{A}(T,N[1])=0
    &&
    \text{By $S[1]\le T\le S$ and $N\in\scrH_S$}\\
    &
    \Leftrightarrow
    \text{$\Hom_{A}(d,N)$ is surjective}\\
    &
    \Leftrightarrow
    \text{$\Hom_{\Lambda}(\Phi(d),\Phi(N))$ is surjective}\\
    &
    \Leftrightarrow 
    \Phi(N)\in\Fac H^0(\Phi(T))
    &&
    \text{Since $\Phi(T)$ is two-term silting of $\Lambda$}.
\end{align}
\end{proof}
\begin{lem}\label{lem1}
    Suppose that $A$ is silting discrete. Let $S\in\per A$ be a basic silting object and let $(\scrD^{\le0},\scrD^{\ge 0})$ be a $t$-structure on $\pvd A$. If $(\scrD^{\le0},\scrD^{\ge 0})$ is intermediate with respect to $(\scrD^{\le0}_S,\scrD^{\ge 0}_S)$, then  $(\scrD^{\le0},\scrD^{\ge 0})$ is algebraic.
\end{lem}
\begin{proof}
Since $A$ is silting discrete, $\{T\in\silt(\per A)\mid S[1]\le T\le S\}$ is finite set. Thus $\overline{\Lambda}$ is $\tau$-tilting finite by Lemma \ref{comm}.  Thus $\tors(\overline{\Lambda})=\ftors(\overline{\Lambda})$; therefore, the claim follows from Lemma \ref{comm}.
\end{proof}
\begin{lem}\label{Lem2}
     The following statements hold:
     \begin{enumerate}
         \item Let $S\in\per A$ be a basic silting object and let $(\scrD^{\le0},\scrD^{\ge 0})$ be a bounded $t$-structure on $\pvd A$. Then there exist integers $m\ge n$ such that $\scrD^{\le 0}_S[m]\subseteq \scrD^{\le 0}\subseteq \scrD^{\le 0}_S[n]$.
         \item Suppose that $A$ is silting discrete. Then there exists a basic silting object $S\in\per A$ such that $\scrD^{\le 0}_S[1]\subseteq \scrD^{\le 0}\subseteq \scrD^{\le 0}_S$.
     \end{enumerate}
\end{lem}
\begin{proof}
    The proof is similar to \cite[Lemmas 8 and 9]{Contractibility-of-the-stability-manifold-for-silting-discrete-algebras}.
\end{proof}
\begin{thm}\label{all bdd t-str is algebraic}
     If $A$ is silting discrete, then any bounded $t$-structure in $\pvd A$ is algebraic.  
\end{thm}
\begin{proof}
    This follows from Lemma \ref{Lem2} and Lemma \ref{lem1}.
\end{proof}

As in \cite{Contractibility-of-the-stability-manifold-for-silting-discrete-algebras}, we expect the space of stability conditions over silting discrete dg algebra to be contractible, although we have not yet proved this.
\section{Compound DuVal singularities}
\begin{setup}In this section, we work under Setup \ref{intro setup tilting stEnd}.
\end{setup}

In this section, we prove Theorem \ref{intro tilting}. The derived equivalence of the stable endomorphism algebras of maximal modifying generators is not difficult to establish, since August's proof applies in our setting using the reduction theorem from section \ref{NCCR-1}. The nontrivial part is the commutativity of the diagram, which requires some additional calculations.

\begin{nota}
Let $N\in\MMG(R)$ and let $N=\bigoplus\limits_{i=0}^nN_i$ be an indecomposable decomposition with $N_0=R$. For $1\le i\le n$, we set
    \[P_i^N:=\stHom_R(N,N_i)\quad P_i^{\oN}:=\stHom_{\oR}(\oN,\oN_i).\]
    These are indecomposable projective modules. Moreover, we set 
    \[S_i^N:=\text{simple $\bsL_N$-module corresponding to $P_i^{N}$}\]
    \[S_i^{\oN}:=\text{simple $\bsL_{\oN}$-module corresponding to $P_i^{\oN}$}.\]
    Note that $S_i^N\cong S_i^{\oN}$ by $\overline{\bsL_N}\cong \bsL_{\oN}$ and $t\in\rad\bsL_{N}$. Moreover, let
    \[\iota:\Db(\fl \bsL_{\oM})\to \Db(\fl \bsL_M)\]
    be the forgetful functor. Finally, we set
    \[\scrK_N:=K_0(\Db(\fl \bsL_{N}))=\bigoplus\limits_{i=1}^n\bZ [S_i^N]
    \]
    \[\scrK_{\oN}:=K_0(\Db(\fl \bsL_{\oN}))=\bigoplus\limits_{i=1}^n\bZ [S_i^{\oN}].
    \]
    Note that these are naturally isomorphic under the functor $\iota$.
\end{nota}
\subsection{Tilting theory for stable endomorphism algebras}

First, following August's argument, we prove Theorem \ref{intro tilting} in a special case. Namely, we consider the case where $M$ and $N$ are related by a simple mutation.

We recall some basic facts about two-sided tilting complexes.
\begin{dfn}
    Let $\Lambda$ and $\Gamma$ be rings, and let $X$ be a complex of $(\Gamma,\Lambda)$-bimodules. Then $X$ is called a {\it two-sided tilting $(\Gamma,\Lambda)$-complex} if 
    \begin{enumerate}
        \item The map $\Gamma\to \RHom_{\Lambda}(X,X)$ is a quasi-isomorphism and
        \item X is quasi-isomorphic to some complex $T\in\Kb(\proj\Lambda)$ and  
        \item $\thick(T)=\Kb(\proj\Lambda)$.
    \end{enumerate}
\end{dfn}

Next, we prove the lemmas needed for the proof.

\begin{lem}\label{hs}
    $R$ is a hypersurface.
\end{lem}
\begin{proof}
This follows from easy commutative ring theory.
\end{proof}
\begin{lem}
Any two-term silting complex over $\bsL_M$ is tilting.
\end{lem}
\begin{proof}
    Let $P\in\twosilt (\bsL_M)$. We consider an exact sequence 
    \[\Hom_{\bsL_M}(P,P[i])\xrightarrow{t}\Hom_{\bsL_M}(P,P[i])\to \Hom_{\bsL_{\oM}}(\oP,\oP[i])\]
    induced by the triangle $P\to P\to \overline{P}\to P[1]$ in $\D(\bsL_{M})$. Since $\bsL_{\oM}$ is symmetric, $\twosilt (\bsL_{\oM})=\twotilt(\bsL_{\oM})$. Thus, we have $\Hom_{\bsL_{\oM}}(\oP,\oP[i])=0$ for each $i\ne 0$. Therefore, $\Hom_{\bsL_{M}}(P,P[i])=0$ for each $i\ne 0$ by NAK.
\end{proof}
\begin{prop}\label{August Lemma}
    Suppose that $N=\mu_i M$ for some $i$, and let \[0\to M_i\xrightarrow{b_i} V_i\xrightarrow{d_i} K_i\to 0\] be an exchange sequence. Let \[P_N:=(\stHom_R(M,M_i)\xrightarrow{b_i\cdot -}\stHom_R(M,V_i))\oplus(0\to \bigoplus\limits_{j\ne i}\stHom_R(M,M_j)).\] Then 
    \begin{enumerate}
    \item $P_N$ is a tilting complex over $\bsL_M$.
    \item $\mathscr{T}_N\cong P_N$ in $\Db(\bsL_M)$.
    \end{enumerate}
\end{prop}
\begin{proof}
    By Theorem \ref{2-siltA->2-siltA/t}, the same argument as for \cite[Theorem 2.16 and Proposition 3.5]{The-tilting-theory-of-contraction-algebras} applies. Note that the three-dimensional assumption in \cite{The-tilting-theory-of-contraction-algebras} is not essential. The crucial point is that the ring is a hypersurface, and the same argument applies in the present setting.
\end{proof}

Recall that for $N\in\MMG(R)$, we set
        \[\scrT_N:=\tau^{\ge-1}( \bsL_{N}\underset{\bLam_N}{\Lo}T_N\underset{\bLam_M}{\Lo} \bsL_{M})\quad\text{and}\quad  \scrT_{\oN}:=\tau^{\ge-1}( \bsL_{\oN}\underset{\bLam_{\oN}}{\Lo} T_{\oN}\underset{\bLam_{\oM}}{\Lo} \bsL_{\oM}).\]
        
We prove the main theorem in a special case.

\begin{thm}\label{tilting for contraction algebra}
    If $N=\mu_i M$ for some $1\le i\le n$, then $\scrT_N$ is a two-sided tilting $(\bsL_N,\bsL_M)$-complex, and there exists the following commutative diagram:
    \[\begin{tikzcd}
	{\Db(\modcat \bsL_M)} && {\Db(\modcat \bLam_M)} & \\
	& {\Db(\modcat \bsL_{\oM})} && {\Db(\modcat \bLam_{\oM})} \\
	{\Db(\modcat \bsL_{N})} && {\Db(\modcat \bLam_{N})} \\
	& {\Db(\modcat \bsL_{\oN})} && {\Db(\modcat \bLam_{\oN})}
	\arrow[from=1-1, to=1-3]
	\arrow["{-\underset{R}{\Lo}\oR}"{description}, from=1-1, to=2-2]
	\arrow["{\RHom_{\bsL_{M}}(\scrT_N,-)}"{description, pos=0.3}, from=1-1, to=3-1]
	\arrow["{-\underset{R}{\Lo}\oR}"{description}, from=1-3, to=2-4]
	\arrow["{\RHom_{\bLam_M}({T_N},-)}"{description, pos=0.3}, from=1-3, to=3-3]
	\arrow[from=2-2, to=2-4]
	\arrow["{\RHom_{\bsL_{\oM}}(\scrT_{\oN},-)}"{description, pos=0.3}, from=2-2, to=4-2]
	\arrow["{\RHom_{\bLam_{\oM}}(T_{\oN},-)}"{description, pos=0.3}, from=2-4, to=4-4]
	\arrow[from=3-1, to=3-3]
	\arrow["{-\underset{R}{\Lo}\oR}"{description}, from=3-1, to=4-2]
	\arrow["{-\underset{R}{\Lo}\oR}"{description}, from=3-3, to=4-4]
	\arrow[from=4-2, to=4-4]
\end{tikzcd}\]
\end{thm}
\begin{proof}
    Let $X:= \bsL_{N}\underset{\bLam_N}{\Lo}T_N\underset{\bLam_M}{\Lo} \bsL_{M}$. We prove only the case $m=1$. The general case follows by repeating the same argument. Let $t:=t_1$.
    
\medskip{\it step1: $\oR\underset{R}{\Lo}X\cong 
\bsL_{\oN}\underset{\bLam_{\oN}}{\Lo}T_{\oN}\underset{\bLam_{\oM}}{\Lo} \bsL_{\oM} $ in $\D(\bsL_{\oM})$.}
    This follows from the following computation:
    \begin{align}
        \oR\underset{R}{\Lo} X
        &=(\oR\underset{R}{\Lo} \bsL_N)\underset{\bLam_N}{\Lo} T_{N}\underset{\bLam_M}{\Lo}  \bsL_M\\
        &=(\oR\underset{R}{\otimes} \bsL_N)\underset{\bLam_N}{\Lo} T_{N}\underset{\bLam_M}{\Lo}  \bsL_M
        && \text{By Lemma \ref{regular sequence}}\\
        &=\bsL_{\oN}\underset{\bLam_N}{\Lo} T_{N}\underset{\bLam_M}{\Lo}  \bsL_M\\
        &=\bsL_{\oN}\underset{\bLam_{\oN}}{\Lo}(\bLam_{\oN}\underset{\bLam_N}{\Lo} T_{N})\underset{\bLam_M}{\Lo} \bsL_M\\
        &=\bsL_{\oN}\underset{\bLam_{\oN}}{\Lo}(\bLam_{\oN}\underset{\bLam_N}{\otimes} T_{N})\underset{\bLam_M}{\Lo} \bsL_M\\
        &=\bsL_{\oN}\underset{\bLam_{\oN}}{\Lo}T_{\oN}\underset{\bLam_M}{\Lo} \bsL_M\\
        &=\bsL_{\oN}\underset{\bLam_{\oN}}{\Lo}T_{\oN}\underset{\bLam_{\oM}}{\Lo}(\bLam_{\oM}\underset{\bLam_M}{\Lo} \bsL_M)\\
        &=\bsL_{\oN}\underset{\bLam_{\oN}}{\Lo}T_{\oN}\underset{\bLam_{\oM}}{\Lo}(\bLam_{\oM}\underset{\bLam_M}{\otimes} \bsL_M)
        && \text{By Lemma \ref{regular sequence}}\\
        &=\bsL_{\oN}\underset{\bLam_{\oN}}{\Lo}T_{\oN}\underset{\bLam_{\oM}}{\Lo} \bsL_{\oM}
    \end{align}

\medskip{\it step2: $\oR\underset{R}{\Lo} \scrT_N\cong \scrT_{\oN}$ in $\D(\bsL_{\oM})$.}  By step1, we have the following commutative diagram:
\[\begin{tikzcd}
	{\oR\underset{R}{\Lo} \tau ^{\le -2}X} & {\oR\underset{R}{\Lo} X} & {\oR\underset{R}{\Lo} \tau ^{\ge -1}X} & {} \\
	{\tau ^{\le -2}Y} & Y & {\tau ^{\ge -1}Y} & {}
	\arrow[from=1-1, to=1-2]
	\arrow[dashed, from=1-1, to=2-1]
	\arrow[from=1-2, to=1-3]
	\arrow["\cong", from=1-2, to=2-2]
	\arrow[from=1-3, to=1-4]
	\arrow["\alpha", dashed, from=1-3, to=2-3]
	\arrow[from=2-1, to=2-2]
	\arrow[from=2-2, to=2-3]
	\arrow[from=2-3, to=2-4]
\end{tikzcd}\]
Here $Y:=\bsL_{\oN}\underset{\bLam_{\oN}}{\Lo}T_{\oN}\underset{\bLam_{\oM}}{\Lo} \bsL_{\oM}$. Note that $\scrT_{\oN}=\tau^{\ge -1}Y$ by definition. We show that $\alpha$ is a quasi-isomorphism.

First, we show that $H^i\alpha$ is an isomorphism for $i\ge-1$. Consider an exact sequence:
\[\begin{tikzcd}
	{H^i(\oR\underset{R}{\Lo} \tau ^{\le -2}X)} & {H^i(\oR\underset{R}{\Lo} X)} & {H^i(\oR\underset{R}{\Lo} \tau ^{\ge -1}X)} & {H^{i+1}(\oR\underset{R}{\Lo} \tau ^{\le -2}X)}
	\arrow[from=1-1, to=1-2]
	\arrow[from=1-2, to=1-3]
	\arrow[from=1-3, to=1-4]
\end{tikzcd}\]
If $i\ge -1$, then $H^i(\oR\underset{R}{\Lo} \tau ^{\le -2}X)=H^{i+1}(\oR\underset{R}{\Lo} \tau ^{\le -2}X)=0$. Hence $H^i(\oR\underset{R}{\Lo} \tau ^{\ge -1}X)\cong H^i(\oR\underset{R}{\Lo} X)\cong H^i(Y)\cong H^i(\tau^{\ge-1}Y)$; therefore, $H^i\alpha$ is an isomorphism for $i\ge-1$. 

Next, we show that $H^i\alpha$ is an isomorphism for $i<-2$. There exists an exact sequence:
\[H^i(\tau^{\ge -1}X)\to H^i(\oR\underset{R}{\Lo}\tau^{\ge -1}X)\to H^{i+1}(\tau^{\ge -1}X)\xrightarrow{t} H^{i+1}(\tau^{\ge -1}X)
\]
induced by a triangle $\tau^{\ge -1}X\to \tau^{\ge -1}X\to \oR\underset{R}{\Lo}\tau^{\ge -1}X\to $. Since $H^i(\tau^{\ge -1}X)=H^{i+1}(\tau^{\ge -1}X)=0$, we have $H^{i}( \oR\underset{R}{\Lo}\tau^{\ge -1}X)=0$; in particular, $H^i\alpha$ is an isomorphism for $i<-2$. 

Finally, we show that $H^{-2}\alpha$ is an isomorphism. It suffices to show that $H^{-2}(\oR\underset{R}{\Lo}\tau^{\ge -1}X)=0$. First,
\begin{align}
    H^{-1}(\tau^{\ge -1}X)
    &=H^{-1}(X)\\
    &\cong H^{-1}(P)
    && \text{By Proposition \ref{August Lemma}}\\
    &\cong\Ker(\stHom_R(M,M_i)\to \stHom_R(M,V_i))\\
    &\cong\Cok(\stHom_R(M,\Omega V_i)\to \stHom_R(M,\Omega K_i)).
\end{align}
Moreover, by Lemma \ref{hs}, $\Omega^2\cong \id$. Hence $\stHom_R(M,\Omega V_i)\cong \stHom_R(\Omega M,V_i)=\Ext_R^1(M,V_i)=0$. Thus  $H^{-1}(\tau^{\ge -1}X)\cong \stHom_R(M,\Omega K_i)$; hence, we obtain $t\in\reg H^{-1}(\tau^{\ge -1}X)$ by Lemma \ref{regular sequence}. Therefore $H^{-2}(\oR\underset{R}{\Lo}\tau^{\ge -1}X)=0$.
The above argument shows that $\alpha$ is a quasi-isomorphism.

\medskip{\it step3: $\scrT_N$ is a two-sided tilting $(\bsL_N,\bsL_M)$-complex.}
By Proposition \ref{August Lemma}, it suffices to show that $\bsL_N\to \Hom_{\bsL_M}(\scrT_N,\scrT_N)$ is an isomorphism. Since
\begin{align}
\Hom_{\bsL_M}(\scrT_N,\scrT_N)\underset{R}{\Lo}\oR
&\cong\RHom_{\bsL_M}(\scrT_N,\scrT_N)\underset{R}{\Lo}\oR 
&& \text{Since $\scrT_N$ is tilting}\\
& \cong \RHom_{\bsL_{\oM}}(\scrT_N\underset{R}{\Lo}\oR,\scrT_N\underset{R}{\Lo}\oR)
&& \text{Since $\scrT_N$ is perfect}\\
& \cong \RHom_{\bsL_{\oM}}(\scrT_{\oN},\scrT_{\oN})
&& \text{By step 2}\\
& \cong\Hom_{\bsL_{\oM}}(\scrT_{\oN},\scrT_{\oN})
&& \text{Since $\scrT_{\oN}$ is tilting}
\end{align}
Moreover,
\begin{align}
    \bsL_{\oN}
    &\cong\bsL_{N}\underset{R}{\Lo}\oR
    && \text{By Lemma \ref{regular sequence}}\\
    &\to \Hom_{\bsL_{M}}(\scrT_N,\scrT_N)\underset{R}{\Lo}\oR\\
    & \cong\Hom_{\bsL_{\oM}}(\scrT_{\oN},\scrT_{\oN})
    && \text{By the above computation}
\end{align}
is an isomorphism by \cite[Lemma 3.9]{The-tilting-theory-of-contraction-algebras}. Therefore, by derived NAK, $\bsL_N\to \Hom_{\bsL_M}(\scrT_N,\scrT_N)$ is an isomorphism. 

\medskip{\it step4: commutativity of the diagram.}
The commutativity of the top and bottom faces is clear. The front and back faces are commutative by an argument similar to that in \cite[Theorem 3.2]{The-tilting-theory-of-contraction-algebras}. The commutativity of the right face follows from Theorem \ref{tilting of NCCR}, while that of the left face follows from Step~2 and an argument similar to the proof of Theorem \ref{tilting of NCCR}.

\end{proof}

\subsection{Preliminaries on hyperplane arrangements}
In this subsection, we give the complete proof of Theorem \ref{intro tilting}. As in August's proof, we study the behavior of compositions of the derived equivalences arising from simple mutations. For this purpose, we recall the hyperplane arrangement associated with the singularity.

We fix an indecomposable decomposition $M=\bigoplus\limits_{i=0}^nM_i$ with $M_0=R$. For $L\in\MMG(R)$ and a sequence $p=(i_{\ell},\dots,i_1)\quad (1\le i_j\le n)$, we put $\mu_pL:=\mu_{i_{\ell}}\cdots\mu_{i_1}L$. By Theorem \ref{tilting of NCCR}, there exists a composition of standard equivalences:
\[\Psi_{p}:=\Psi_{i_{\ell}}\circ\cdots\circ \Psi_{i_1}:\Db(\modcat \bLam_L)\xrightarrow{\sim} \Db(\modcat \bLam_{\mu_{p} L}).\]
This equivalence restricts to the following equivalence:
\[\Psi_{p}:\Db(\fl \bsL_{L})\xrightarrow{\sim} \Db(\fl \bsL_{\mu_{p} L}).\]
Moreover, this equivalence induces an isomorphism of Grothendieck groups: 
\[\Psi^{K}_p:\scrK_L\to \scrK_{\mu_{p} L}.\]
We call the sequence $p=(i_{\ell},\dots,i_1)$ a {\it positive path} from $L$ to $\mu_p L$, and call $\ell$ the {\it length} of $p$. If $p$ is a positive path from $L$ to $\mu_p L$, we write $p:L\to \mu_p L$. For each $N\in\MMG(R)$, set
\[
\Theta_N:=\Hom_{\bZ}(\scrK_N,\bR),
\]
and let $e_1^N,\dots,e_n^N$ be the dual basis of
$[S_1^N],\dots,[S_n^N]$. For a positive path $p:M\to N$, define
\[
C_+^N:=\sum_{i=1}^n\bR_{>0}e_i^N
\quad
C_p:=(\Psi_p^K)^*(C_+^N)\subseteq\Theta_M.
\]
Here $(\Psi_p^K)^*(\theta)=\theta\circ\Psi_p^K$.
Next, we show that $C_p$ is independent of the choice of $p$. For this purpose, we use Lemma \ref{C_L}, Lemma \ref{^2=1} and Corollary \ref{CL cor}. Since their proofs are analogous to those in \cite{Stability-conditions-for-contraction-algebras,Stability-conditions-for-3-fold-flops}, we omit them.
\begin{lem}\label{C_L}
    Let \[0\to K_i\to \bigoplus\limits_{j\ne i}L_j^{\oplus c_{ij}}\to L_i \to 0\]
    be an exchange sequence. Then $\Psi^K_i:\scrK_L\xrightarrow{\sim}\scrK_{\mu_i L}$ sends:
    \begin{equation}
  [S_j^L] \longmapsto
  \begin{cases}
    -[S_i^{\mu_i L}] & \text{if } j=i,\\
    c_{ij}[S_i^{\mu_i L}]+[S_j^{\mu_i L}] & \text{if } j\neq i.
  \end{cases}
\end{equation}
\end{lem}
\begin{rem}
Taking the dual of the map in Lemma \ref{C_L}, we obtain
\[
(\Psi_i^K)^*:\Theta_{\mu_iL}\xrightarrow{\sim}\Theta_L
\]
which sends
\[
e_j^{\mu_iL}\longmapsto
\begin{cases}
-e_i^L+\displaystyle\sum_{k\ne i}
c_{ik}e_k^L & \text{if } j=i,\\[4pt]
e_j^L & \text{if } j\ne i.
\end{cases}
\]
\end{rem}

By Corollary \ref{NCCR-1 cor}, we have $\mu_i\mu_iL\cong L$. Thus 
\[\begin{tikzcd}
	{\Db(\modcat \bsL_{L})} & {\Db(\bsL_{\mu_i L}).}
	\arrow["{\Psi_i}", shift left, from=1-1, to=1-2]
	\arrow["{\Psi_i}", shift left, from=1-2, to=1-1]
\end{tikzcd}\]
Therefore, we obtain the following four isomorphisms of Grothendieck groups:
\[\begin{tikzcd}
	{\scrK_L} & {\scrK_{\mu_i L}}
	\arrow["{\Psi^K_i}", shift left, from=1-1, to=1-2]
	\arrow["{\Psi^K_i}", shift left, from=1-2, to=1-1]
    \end{tikzcd}
    \quad 
    \begin{tikzcd}
    {\scrK_L} & {\scrK_{\mu_i L}}
	\arrow["{(\Psi_i^{-1})^K}", shift left, from=1-1, to=1-2]
	\arrow["{(\Psi_i^{-1})^K.}", shift left, from=1-2, to=1-1]
    \end{tikzcd}\]

\begin{lem}\label{^2=1}
    The above four isomorphisms are given by the same matrix, and this matrix squares to the identity. 
\end{lem}

\begin{cor}\label{CL cor}
    In particular, for any two paths $p,q:M\to N$, we have $C_p=C_q$. Thus, we denote 
    \[C_N:=C_p=C_q.\]
\end{cor}

Note that $\MMG(R)$ is connected by Theorem \ref{MM(R)->MM(R/t) MM} and \cite{IW2023tits}. We define the {\it finite mutation cone} of $R$ by 
\[\Cone(R):=\bigcup\limits_{N\in\MMG(R)} C_N\subseteq \bR^n.\]

The preceding argument applies verbatim to $\oR$, and hence $C_N=C_{\oN}$ follows from Lemma \ref{C_L}. Thus, $\Cone(R)=\Cone(\oR)$ holds. Therefore, by Theorem \ref{MM(R)->MM(R/t) MM} and \cite{IW2023tits}, we obtain the following results.

\begin{prop}
    The following statements hold.
    \begin{enumerate}
        \item The open decomposition $\bigcup\limits_{N\in\MMG(R)} C_N$ gives the chambers of the $\scrJ$-finite hyperplane arrangement. (cf. \cite{IW2023tits}.)
        \item For $N,L\in\MMG(R)$, $N\ncong L$ if and only if $C_N\cap C_L=\varnothing$.
        \item For $N,L\in\MMG(R)$, $N\cong \mu_iL$ for some $i\ne 0$ if and only if $\overline{C_N}$ and $\overline{C_L}$ share a codimension-one wall where $\overline{C_N}$ and $\overline{C_L}$ is the closure.
    \end{enumerate}
    Therefore, we obtain a finite simplicial hyperplane arrangement $\scrH_R:=\bR^n\setminus \bigcup\limits_{N\in\MMG(R)} C_N$.
\end{prop}

We next recall some basics of a hyperplane arrangement. Let  $\scrH:=\scrH_R$.

\begin{dfn} \label{def:Gamma graph}
The graph $\Gamma_{\cH}$ of oriented arrows is defined as follows.  The vertices of $\Gamma_{\cH}$ are the chambers $\{C_N\}_{N\in\MMG(R)}$ of $\bR^n\setminus\scrH$. For chambers $C_N$ and $C_L$, there is an arrow $C_N\to C_L$ if and only if $C_N$ and $C_L$ are adjacent.
\end{dfn}

A {\it positive path of length $\ell$} in $\Gamma_{\cH}$ is a formal symbol
\[
p=a_{\ell}\circ \cdots\circ a_2\circ a_1,
\]
such that there exists a sequence of vertices $v_0,\cdots,v_{\ell}$ of $\Gamma_{\cH}$ and there exist arrows $a_i: v_{i-1}\to v_i$ in $\Gamma_{\cH}$. Set $s(p):=  v_0$, $t(p):=  v_{\ell}$, and write $p: s(p)\to t(p)$. The notation $\circ$ should remind us of composition, but we will often drop the $\circ$'s in future.  If 
$q=b_m\circ\hdots\circ b_2 \circ b_1$ is another positive path with $t(p)=s(q)$, we consider the formal symbol
\[
q\circ p:=  b_m\circ\cdots\circ b_2 \circ b_1
\circ
a_{\ell}\circ \cdots\circ a_2\circ a_1,
\]
and call it the {\it composition} of $p$ and $q$. 

\begin{dfn}
A positive path is called {\it reduced} if it does not cross any hyperplane twice. 
\end{dfn}

Let  $\mathrm{Free}(\Gamma_\cH)$ be the free category on the graph $\Gamma_\cH$, where morphisms are directed paths. We define an equivalence relation $\sim$ on $\Free(\Gamma_\scrH)$ to be the smallest equivalence relation satisfying the following conditions:
\begin{enumerate}
    \item If $p\sim q$, then $s(p)=s(q)$ and $t(p)=t(q)$.
    \item If $p,q$ are positive minimal paths with $s(p)=s(q)$ and $t(p)=t(q)$, then $p\sim q$.
    \item If $p\sim q$, then $upr\sim uqr$ for all positive paths $u$ and $r$ with $t(r)=s(p)=s(q)$ and $s(u)=t(p)=t(q)$.
\end{enumerate}
We define $\bG^{+}$ to be the quotient category of $\Free(\Gamma_{\scrH})$ by the equivalence relation $\sim$:
\[
\bG^{+}_\scrH:= \Free(\Gamma_\scrH)/\sim. 
\]
\begin{dfn} \label{def: Deligne groupoid completion}
Let $\bG_\scrH$ be the groupoid obtained from $\bG^+_\scrH$ by formally adjoining inverses to all morphisms. We call $\bG_\scrH$ the {\it Deligne groupoid}.
\end{dfn} 
The following is well-known \cite{Les-immeubles-des-groupes-de-tresses-generalises,Universal-cover-of-Salvettis-complex-and-topology-of-simplicial-arrangements-of-hyperplanes,The-covers-of-a-complexified-real-arrangement-of-hyperplanes-and-their-fundamental-groups,Topology-of-the-complement-of-real-hyperplanes-in-Cn}.
\begin{prop}\label{ver gp}
For any vertex $v$ in the Deligne groupoid, $\End_{\bG}(v)\cong\pi_1(\mathbb{C}^n \backslash \scrH_\mathbb{C})$.
\end{prop}

Since $\bG_{\scrH}$ is determined solely by the hyperplane arrangement $\scrH$, the following is immediate.

\begin{prop}
There exists a natural equivalence:
\[\overline{(-)}:
\bG_{\scrH_{R}}\to \bG_{\scrH_{\oR}}
.\]
\end{prop}
We conclude this section by proving a result concerning compositions of the functors $\Psi_i$. This result is an analogue of \cite[Theorem 4.12]{The-tilting-theory-of-contraction-algebras}.

\begin{prop}\label{tilting composition}
    Let $p=(C_M\xrightarrow{a_{i_1}}C_{\mu_{i_1}M}\xrightarrow{a_{i_2}}\cdots \xrightarrow{a_{i_{\ell}}}C_{\mu_{i_\ell}\cdots\mu_{i_1}M}=C_N)$ be a positive path in $\Gamma_{\scrH}$. If $p$ is minimal, then 
    \[\Psi_{p}:=\Psi_{i_\ell}\circ\cdots \Psi_{i_1}\cong\RHom_{\bsL_M}(\scrT_N,-):\Db(\modcat\bsL_M)\xrightarrow{\sim}\Db(\modcat\bsL_N).\]
\end{prop}
\begin{proof}
For each $j$, let 
\begin{itemize}
    \item $N_{(j)}:=\mu_{i_j}\cdots\mu_{i_1}M$
    \item $\bLam_{(j)}:=\End_R(N_{(j)})$ and $\overline{\bLam}_{(j)}:=\End_{\oR}(\oN_{(j)})$
    \item $\bsL_{(j)}:=\stEnd_R(N_{(j)})$
    \item $T_{(j)}:=\Hom_R(N_{(j-1)},N_{(j)})$ and $\oT_{(j)}:=\Hom_{\oR}(\oN_{(j-1)},\oN_{(j)})$
    \item $\scrT_{(j)}:=\tau^{\ge -1}(\bsL_{(j)}\underset{\bLam_{(j)}}{\Lo}T_{(j)}\underset{\bLam_{(j-1)}}{\Lo}\bsL_{(j-1)})$
\end{itemize}

\medskip{\it step1: The natural map $T_{(\ell)}\underset{\bLam_{(\ell-1)}}{\Lo}\cdots\underset{\bLam_{(1)}}{\Lo}T_{(1)}\to T_N$ is a quasi-isomorphism of $(\bLam_{N},\bLam_M)$-bimodule complexes.} Applying $-\underset{R}{\Lo}\oR$ to this morphism, we obtain the natural morphism
\[\oT_{(\ell)}\underset{\overline{\bLam}_{(\ell-1)}}{\Lo}\cdots\underset{\overline{\bLam}_{(1)}}{\Lo}\oT_{(1)}\to T_{\oN}\]
The above natural morphism is an isomorphism by \cite[Theorem 4.6]{Faithful-actions-from-hyperplane-arrangements}. Hence, by derived NAK, the original morphism is also an isomorphism.

\medskip{\it step2: There exists a quasi-isomorphism $\scrT_{(\ell)}\underset{\bsL_{(\ell-1)}}{\Lo}\cdots\underset{\bsL_{(1)}}{\Lo}\scrT_{(1)}\cong \scrT_N$ of $(\bsL_{N},\bsL_M)$-bimodule complexes.}
The proof is similar to that in \cite[Theorem 4.12]{The-tilting-theory-of-contraction-algebras}. Note that the auxiliary lemmas needed for this argument can likewise be proved in the present setting by the same method.
\end{proof}
\begin{rem}
    We note that the step 2 in Proposition \ref{tilting composition} cannot be proved in the same way as step 1, since equality $\scrT_N\underset{R}{\Lo}\oR\cong \scrT_{\oN}$ has not yet been established for arbitrary $N$.
\end{rem}
For $\alpha=[p]\in\bG$, set 
\[\Psi_{\alpha}:=\Psi_p.\]
By Proposition \ref{tilting composition}, this is well-defined.

\subsection{Stability spaces}
In this subsection, we consider the $\Stab(\Db(\fl \bsL_{M}))$. First, the following definition is analogous to that in \cite{Stability-conditions-on-noncommutative-crepant-resolutions-of-3-dimensional-isolated-singularities}.
\begin{dfn}\label{MM stab dfn}
    \begin{enumerate}
        \item Let $\scrA$ be a heart of a bounded $t$-structure on $\Db(\fl \bsL_{M})$. We call $\scrA$ {\it maximal modifying (MM)} if there exists $N\in\MMG(R)$ and $\alpha\in\Hom_{\bG}(C_+,C_N)$, such that $\scrA=(\Psi_\alpha)^{-1}(\fl \bsL_{N})$.
        \item We say that a stability condition $(\scrA,Z)\in\Stab\Db(\fl \bsL_{M})$ is {\it maximal modifying (MM)} if $\scrA$ is MM. We denote the set of such stability conditions by $\Stab^{\std}\Db(\fl \bsL_{M})$. 
    \end{enumerate}
    
\end{dfn}

Next, we show that $\stdStab \Db(\fl \bsL_{M})=\Stab^{\alg}\Db(\fl \bsL_{M})=\Stab\Db(\fl \bsL_{M})$ and $\stdStab \Db(\fl \bsL_{\oM})=\Stab\Db(\fl \bsL_{\oM})$. 

\begin{prop}\label{alg=all}
    $\bsL_M$ is silting discrete. In particular, $\Stab^{\alg}\Db(\fl \bsL_{M})=\Stab\Db(\fl \bsL_{M})$. 
\end{prop}
\begin{proof}
    By \cite[Theorem 3.12]{On-the-finiteness-of-the-derived-equivalence-classes-of-some-stable-endomorphism-rings}, $\bsL_{\oM}$ is silting discrete. Moreover, by Corollary \ref{stEnd Stab reduction}, $\overline{(-)}:(\silt(\per \bsL_M),\le)\xrightarrow[]{\sim}(\silt(\per \bsL_{\oM}),\le)$ is an isomoprphism of posets. Therefore, $\bsL_M$ is also silting discrete. Finally, $\Stab^{\alg}\Db(\fl \bsL_{M})=\Stab\Db(\fl \bsL_{M})$ follows from Theorem \ref{all bdd t-str is algebraic}.
\end{proof}

\begin{rem}
Using Theorem \ref{MM(R)->MM(R/t) MM}, we can also prove that $\bsL_M$ is silting discrete by exactly the same argument as \cite[Theorem 3.12]{On-the-finiteness-of-the-derived-equivalence-classes-of-some-stable-endomorphism-rings}.
\end{rem}

More strongly, the following holds. This result is important because it gives an explicit classification of bounded $t$-structures.

\begin{prop}\label{classification of t-str}
    $\stdStab \Db(\fl \bsL_{M})=\Stab\Db(\fl \bsL_{M})$ holds.
\end{prop}
\begin{proof}
    Using Proposition \ref{KY corr} and Theorem \ref{all bdd t-str is algebraic}, the same argument as in the proof of \cite[Proposition 3.3]{Stability-conditions-for-contraction-algebras} is applied.
\end{proof}

Moreover, silting-discreteness yields the following.

\begin{cor}\label{cor stability}
There exists an isomoprphism 
\[\overline{(-)}:\Stab \Db(\fl \bsL_{M})\to \Stab\Db(\fl \bsL_{\oM}).\]
\end{cor}
\begin{proof}
    This follows from Proposition \ref{alg=all} and Corollary \ref{stEnd Stab reduction}.
\end{proof}

\section{Examples and problems}
\subsection{Modifying modules}
In this subsection, we construct modifying modules for hypersurface singularities of a special form.
\begin{setup}
In this subsection, we use the following notation:
    \begin{enumerate}
        \item Let $f_i\ (1\le i\le m)$ and $g_j\ (1\le j\le n)$ be irreducible elements of $S:=k[[x_1,\dots,x_d]]$. Assume that, for any $i$ and $j$, $f_i$ and $g_j$ have no common factor. 
        \item Let $f:=f_1\cdots f_m$, $g:=g_1\cdots g_n$ and let $R:=S/(f+g)$.
    \end{enumerate}
\end{setup}

\begin{lem}\label{hyp surf mofif}
    Assume that $R$ is a normal domain. For any decompositions $abc=f$ and $pq=g$, the following statements hold.
    \begin{enumerate}
        \item $(a,p)\in\CM R$.
        \item $bcp^{-1}:(a,p)\xrightarrow{\cong} (bc,q)$.
        \item $\Hom_R((a,p),(ab,p))\cong (b,p)$.
        \item $\Hom_R((ab,p),(a,p))\cong (ac,p)p^{-1}$.
        % \item $[(a,p)]+[(b,p)]=[(ab,p)]$ in $\Cl(R)$.
    \end{enumerate}
\end{lem}
\begin{proof}
The proof is similar to that of \cite[Lemmas 5.3 and 5.4]{Reduction-of-triangulated-categories-and-maximal-modification-algebras-for-cA_n-singularities}. 
\end{proof}

\subsection{Some hypersurface singularities}\label{ex subsection}
\begin{setup}
    In this subsection, we use the following notations:
    \begin{itemize}
       \item Let $k$ be a field and $S:=k[[t_1,\dots ,t_{d-2},x]]$ for $d\in\bZ_{\ge 2}$.
       \item For $1\le i\le n$, let $f_i\in k[[t_1,\dots ,t_{d-2}]]$ be irreducible elements, and let $F_i:=f_i+x\in S$.
       \item Let $R:=S[[u,v]]/(F_1\cdots F_n+uv)$.
       \item For a permutation $\sigma\in S_n$, let $M_{\sigma}:=\bigoplus\limits_{i=1}^n (F_{\sigma(1)}\cdots F_{\sigma(i)},u)$
   \end{itemize}
\end{setup}
\begin{thm}\label{higher cA NCCR}
   The following statements hold.
   \begin{enumerate}
       \item $R$ is a normal domain.
       \item For each $\sigma$, $M_{\sigma}\in\CM R$ gives an NCCR of $R$.
   \end{enumerate}
\end{thm}
\begin{proof}
(1) Since $t_1,\dots,t_{d-2}$ is a regular sequence, this follows from Corollary \ref{normal +1}.

(2) By Lemma \ref{hyp surf mofif}, $M_{\sigma}$ is a Cohen--Macaulay modifying $R$-module. We consider a regular sequence $t_1,\dots ,t_{d-2}$. Since
$R/(t_1,\dots,t_{d-2})=k[[x,u,v]]/(x^n+uv)$
has an NCCR given by
$M_\sigma/(t_1,\dots,t_{d-2})
=\bigoplus_{i=1}^n(x^i,u)$,
applying Theorem \ref{NCCR +1} inductively to the regular sequence
$t_1,\dots,t_{d-2}$, we obtain the claim.
\end{proof}

We consider a toric singularity that provides an example of Theorem \ref{higher cA NCCR}.

\begin{eg}\label{toric eg}
    Let $R:=\bC[[x_1,\dots,x_n,u,v]]/(x_1\cdots x_n+uv)$. This is the toric singularity associated with the quasi-symmetric representation determined by the following weight vectors:
\[
\begin{pmatrix}
    1 \\ 1 \\ 1 \\ \vdots \\ 1 
\end{pmatrix}
\quad
\begin{pmatrix}
    1 \\ 0 \\ 0 \\ \vdots \\ 0
\end{pmatrix}
\quad
\begin{pmatrix}
    0 \\ 1 \\ 0 \\ \vdots \\ 0
\end{pmatrix}
\quad \cdots\quad
\begin{pmatrix}
    0 \\ 0 \\ 0 \\ \vdots \\ 1
\end{pmatrix}
\quad
\begin{pmatrix}
    -1 \\ -1 \\ -1 \\ \vdots \\ -1
\end{pmatrix}
\quad
\begin{pmatrix}
    -1 \\ 0 \\ 0 \\ \vdots \\ 0
\end{pmatrix}
\quad
\begin{pmatrix}
    0 \\ -1 \\ 0 \\ \vdots \\ 0
\end{pmatrix}
\quad \cdots\quad
\begin{pmatrix}
    0 \\ 0 \\ 0 \\ \vdots \\ -1
\end{pmatrix}
\]

By Theorem \ref{higher cA NCCR}, $\bigoplus\limits_{i=1}^n(x_{\sigma(1)}\cdots x_{\sigma(i)},u)$ gives an NCCR of $R$ for every permutation $\sigma\in S_n$. Moreover, Theorem \ref{NCCR(R)->NCCR(R/t) NCCR} and \cite[Theorem 5.15]{Reduction-of-triangulated-categories-and-maximal-modification-algebras-for-cA_n-singularities} show that these exhaust $\NCCRG(R)$. When $n=3$, the exchange graph is given by Figure \ref{Ex Graph}.

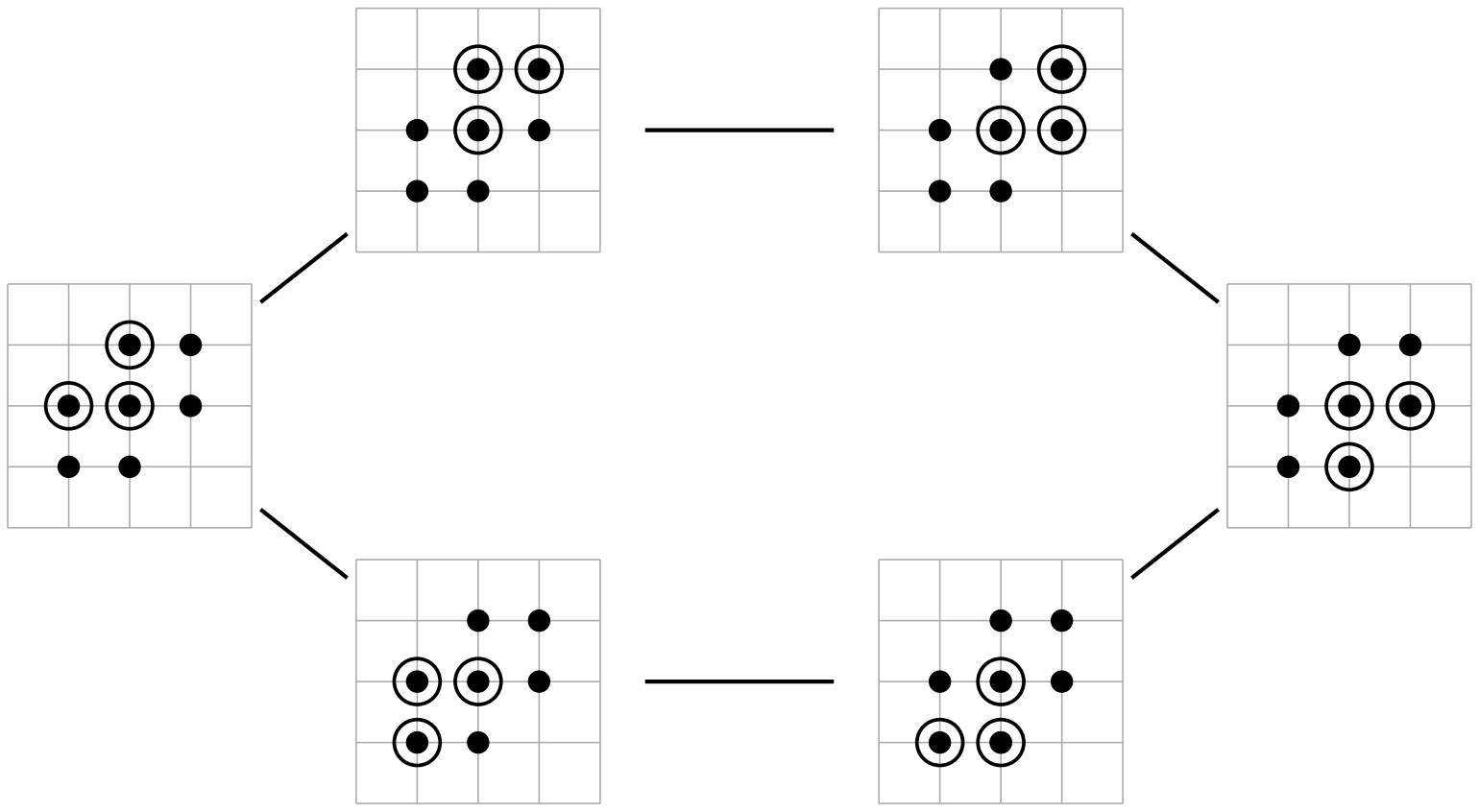
\begin{figure}[htbp]\label{Ex Graph}
\centering

\begin{tikzpicture}[
    gridline/.style={gray!65, line width=0.3pt},
    edge/.style={black, line width=0.8pt}
]

%------------------------------------------------
% mini figure
% #1 : name
% #2 : center
% #3 : black points
% #4 : circled points
%------------------------------------------------
\newcommand{\panel}[4]{%
  \begin{scope}[shift={#2}, x=0.42cm, y=0.42cm]

    % grid
    \foreach \i in {-2,-1,0,1,2}{
      \draw[gridline] (\i,-2) -- (\i,2);
      \draw[gridline] (-2,\i) -- (2,\i);
    }

    % ordinary black points
    \foreach \p in {#3}{
      \fill \p circle (2.2pt);
    }

    % circled black points
    \foreach \p in {#4}{
      \fill \p circle (2.2pt);
      \draw[line width=0.7pt] \p circle (4.5pt);
    }

    % center used for connecting the panels
    \coordinate (#1) at (0,0);

  \end{scope}
}

%================================================
% six mini figures
%================================================

% upper left
\panel{A}{(-1.8,1.9)}
{
  (-1,-1),
  (0,-1),
  (-1,0),
  (1,0)
}
{
  (0,0),
  (0,1),
  (1,1)
}

% upper right
\panel{B}{(1.8,1.9)}
{
  (-1,-1),
  (0,-1),
  (-1,0),
  (0,1)
}
{
  (0,0),
  (1,0),
  (1,1)
}

% right
\panel{C}{(4.2,0)}
{
  (-1,-1),
  (-1,0),
  (0,1),
  (1,1)
}
{
  (0,-1),
  (0,0),
  (1,0)
}

% lower right
\panel{D}{(1.8,-1.9)}
{
  (-1,0),
  (0,-1),
  (0,1),
  (1,1),
  (1,0)
}
{
  (-1,-1),
  (0,0),
  (0,-1)
}

% lower left
\panel{E}{(-1.8,-1.9)}
{
  (0,-1),
  (1,0),
  (0,1),
  (1,1)
}
{
  (-1,-1),
  (-1,0),
  (0,0)
}

% left
\panel{F}{(-4.2,0)}
{
  (-1,-1),
  (0,-1),
  (1,0),
  (1,1)
}
{
  (-1,0),
  (0,0),
  (0,1)
}

%================================================
% exchange graph
%================================================

\draw[edge, shorten <=1.15cm, shorten >=1.15cm]
  (A) -- (B);

\draw[edge, shorten <=1.15cm, shorten >=1.15cm]
  (B) -- (C);

\draw[edge, shorten <=1.15cm, shorten >=1.15cm]
  (C) -- (D);

\draw[edge, shorten <=1.15cm, shorten >=1.15cm]
  (D) -- (E);

\draw[edge, shorten <=1.15cm, shorten >=1.15cm]
  (E) -- (F);

\draw[edge, shorten <=1.15cm, shorten >=1.15cm]
  (F) -- (A);

\end{tikzpicture}

\caption{Here $\bullet$ denotes a CM weight and $\circ$ indicates an NCCR}
\end{figure}
In this case, all the NCCRs given above are isomorphic, and their quivers with relations are given as follows:
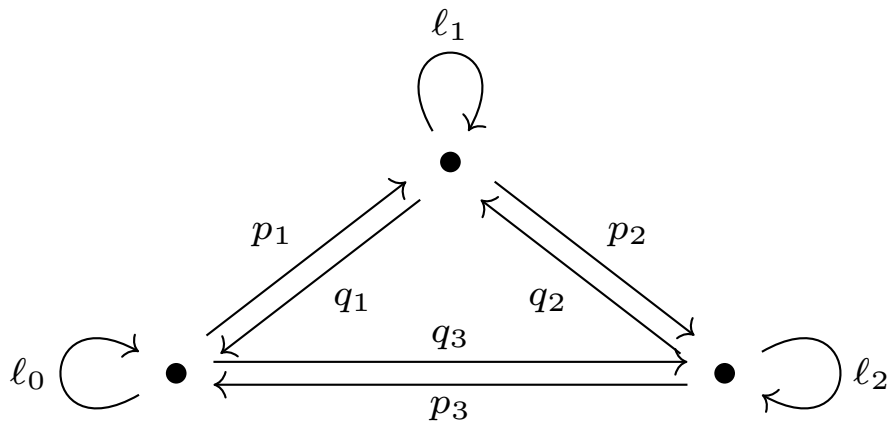
\[
\begin{array}{cc}
\begin{array}{c}
\begin{tikzcd}[ampersand replacement=\&, row sep=2.8em, column sep=3.8em]
  \& \bullet
    \arrow[shift left=.5ex, "q_1", dl]
    \arrow[shift left=.5ex, "p_2", dr]
    \arrow["\ell_1", from=1-2, to=1-2, loop,
           in=60, out=120, distance=8mm]
\\
  \bullet
    \arrow[shift left=.5ex, "p_1", ur]
    \arrow[shift left=.5ex, "q_3", rr]
    \arrow["\ell_0", from=2-1, to=2-1, loop,
           in=150, out=210, distance=8mm]
  \&\&
  \bullet
    \arrow[shift left=.5ex, "q_2", ul]
    \arrow[shift left=.5ex, "p_3", ll]
    \arrow["\ell_2", from=2-3, to=2-3, loop,
           in=-30, out=30, distance=8mm]
\end{tikzcd}
\end{array}
&
{\scriptstyle
\begin{array}{l}
p_1\ell_0=q_2p_2p_1,
\qquad
\ell_1p_1=p_1p_3q_3,
\\[0.3em]
\ell_0q_1=q_1q_2p_2,
\qquad
q_1\ell_1=p_3q_3q_1,
\\[0.3em]
p_2\ell_1=q_3p_3p_2,
\qquad
\ell_2p_2=p_2p_1q_1,
\\[0.3em]
\ell_1q_2=q_2q_3p_3,
\qquad
q_2\ell_2=p_1q_1q_2,
\\[0.3em]
p_3\ell_2=q_1p_1p_3,
\qquad
\ell_0p_3=p_3p_2q_2,
\\[0.3em]
\ell_2q_3=q_3q_1p_1,
\qquad
q_3\ell_0=p_2q_2q_3.
\end{array}}
\end{array}
\]
\end{eg}

% \subsection{Quotient singularities}
% \begin{setup}
%     In this subsection, we use the following notation:
%     \begin{enumerate}
%         \item Let $S$ be a complete Gorenstein nonrmal local domain containg a field $k$. 
%         \item Let $G$ be a finite group acting on $S$ such that $|G|\in k^{\times}$.
%         \item Let $R:=S^G$ be an invariant subring.
%         \item Let $t\in R$ be a regular element.
%     \end{enumerate}
% \end{setup}

% Since $t\in R$, the action of $G$ on $S$ induces an action on $S/tS$.

% \begin{prop}\label{/t ^G comm}
%     $R/tR\cong (S/tS)^G$ holds.
% \end{prop}
% \begin{proof}
    
% \end{proof}

% Recall that, the the action of $G$ on $S$ is {\it small} if... 

% We set
% \begin{itemize}
%     \item Let $\refl^G S$ be the category of reflexive $G$-equivariant $S$-modules.
%     \item Let $\NCCR^G(S):=\{M\in\refl^G S\mid M\in\NCCR(S)\}$.
%     \item Let $\NCCR_0^G(S):=\{M\in\refl^G S\mid M\in\NCCRG(S)\}$.
% \end{itemize}
% Similarly for $S/tS$.

% \begin{cor}
%     Assume that the action of $G$ on $S$ and $S/tS$ are small.  If $M\in\NCCR_0^G(S)$ and $M/tM\in \NCCRG(S/tS)$, then $M\in\NCCR_0^G(R)$ and that $M/tM\in\NCCR(R/tR)$.
% \end{cor}

% \begin{eg}
%      Let $S$ be $\cdots$, and define a $G:=C_2$-action on $S$ as follows:
%      \[.\]
%      Moreover, let 
%      \[f:= \quad\text{and}\quad S:=P/fP=\frac{\bC[[x,y,z,u,v]]}{}.\]
%      Then $f\in P^G$, we obtain
%      \[S^G=(P/tP)^G=\frac{\bC[[]]}{}.\]
%      \[\]
% \end{eg}

\subsection{Universal singularities}\label{univ sing}

Finally, we introduce the notion of universal singularities and universal MMAs. We begin with the following observation. 

\begin{eg}[$cA$ singularities]
    Fix $n>1$. Consider $n$ distinct irreducible element $f_1,\dots,f_n\in (x,y)\subseteq k[[x,y]]$. For the isolated $cA$ singularities 
\[R_{f_1,\dots,f_n}:=\frac{\bC[[x,y,u,v]]}{(f_1\cdots f_n+uv)},\]
the "cluster theories" determined by $R_{f_1,\dots,f_n}$ are all isomorphic. More precisely, there are isomorphisms between exchange quivers of MMs, spaces of stability conditions, and so on. One explanation for this phenomenon is that the "cluster theories" of an isolated cDV singularity is controlled by the hyperplane arrangement associated with the singularity, and in the present situation the hyperplane arrangements associated with the singularities $R_{f_1,\dots,f_n}$ are all identical.

However, Theorem \ref{MM(R)->MM(R/t) MM} and Corollary \ref{cor stability} provide another explanation for this phenomenon. Let 
\[R:=\frac{\bC[[x,y,z_1,\dots,z_n,u,v]]}{(z_1\cdots z_n+uv)}.\]
By choosing a regular sequence $(t_i:=z_i-f_i)_{1\le i\le n}$ on $R$ appropriately, each of the singularities $R_{f_1,\dots,f_{n}}$ can be realized as a quotient $R/(t_1,\ldots,t_{n})$. Moreover, this setting satisfies the assumptions of Setup \ref{intro setup  reduction NCCR,tilt,2-silt}. Thus, Theorem \ref{MM(R)->MM(R/t) MM} and Corollary \ref{cor stability} yield isomorphisms between exchange quivers of MMs, spaces of stability conditions, and so on. In this sense, $R$ may be regarded as a "universal singularity" controlling the family $\{R_{f_1,\dots,f_n}\}_{f_1,\dots,f_n}$ of isolated $cA$ singularities (Figure \ref{ovewview univaarsal sing}).
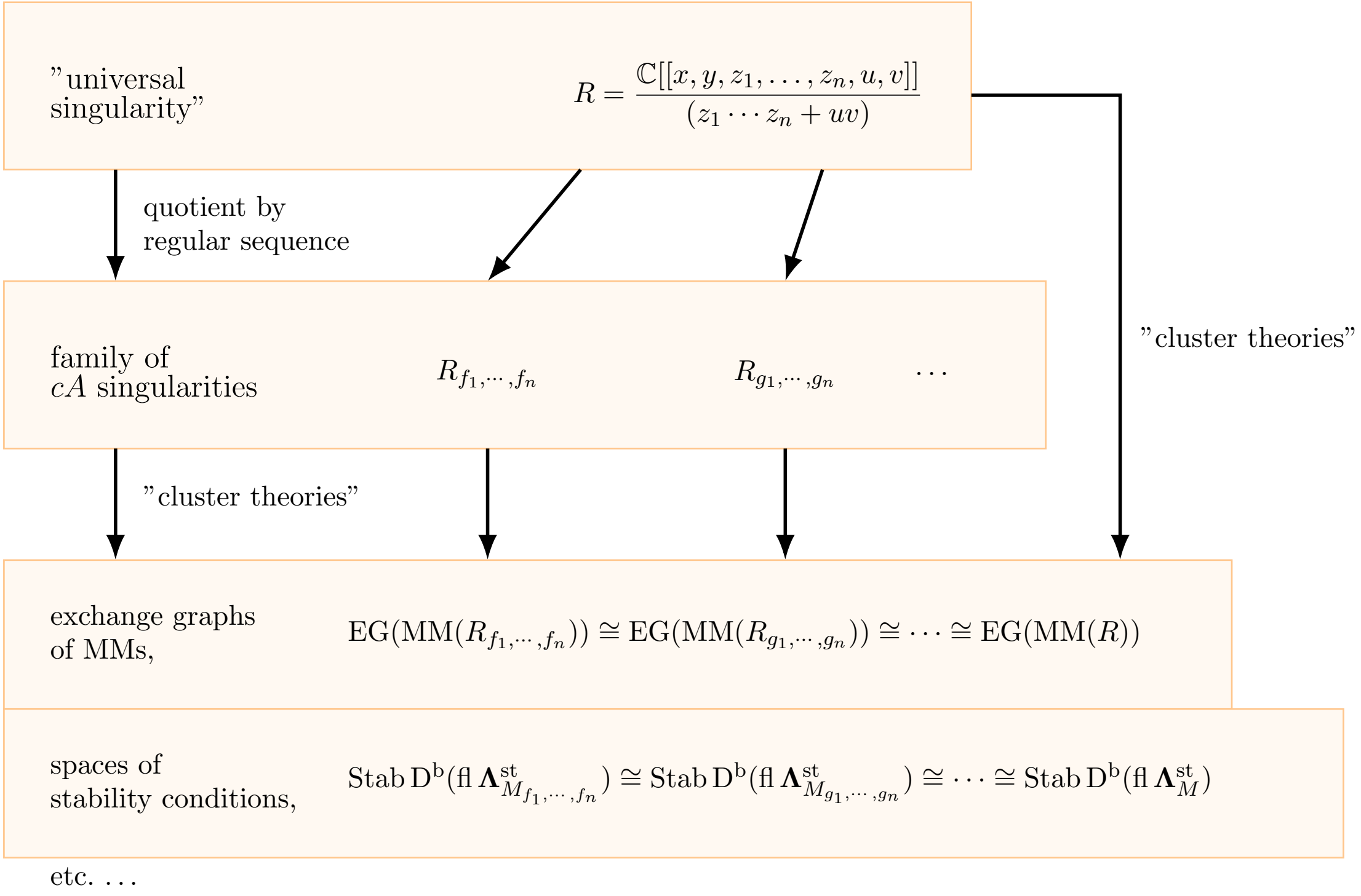
\begin{figure}[http]
\begin{tikzpicture}[scale=0.8,transform shape,x=1cm,y=1cm,>=Latex,
  every node/.style={font=\normalsize},
  arrow/.style={->,line width=1.05pt},
  box/.style={draw=orange!45,fill=orange!5,line width=.5pt},
  math/.style={font=\large}]

% Highlighted regions.
\path[box] (1,11.5) rectangle (14,13.75);
\path[box] (1,7.75)  rectangle (15,10);
\path[box] (1,4.25)  rectangle (17.5,6.25);
\path[box] (1,2.25)  rectangle (19,4.25);

% Top: universal singularity.
\node[align=left,anchor=west,font=\Large] at (1.5,12.5)
  {"universal\\[-1mm]singularity"};
\node[math] (R) at (11,12.5)
  {$\displaystyle R=\frac{\bC[[x,y,z_1,\ldots,z_n,u,v]]}
  {(z_1\cdots z_n+uv)}$};

% Middle: the family of cA singularities.
\node[align=left,anchor=west,font=\Large] (family) at (1.5,8.75)
  {family of\\[-1mm]$cA$ singularities};
\node[math] (Ra) at (7.5,8.75) {$R_{f_1,\cdots,f_n}$};
\node[math] (Rb) at (11.5,8.75) {$R_{g_1,\cdots,g_n}$};
\node[math] at (13.5,8.75) {$\cdots$};

% Quotient arrows and annotation.
\draw[arrow] (2.5,11.5) -- (2.5,10);
\node[align=left,anchor=west,font=\large] at (2.75,10.75)
  {quotient by\\regular sequence};
\draw[arrow] (8.75,11.5) -- (7.5,10);
\draw[arrow] (12,11.5) -- (11.5,10);

% NCCR arrows.
\draw[arrow] (2.5,7.75) -- (2.5,6.25);
\node[anchor=west,font=\large] at (2.75,7.12) {"cluster theories"};
\draw[arrow] (7.5,7.75) -- (7.5,6.25);
\draw[arrow] (11.5,7.75) -- (11.5,6.25);

% Long top-right implication.
\draw[arrow] (14,12.5) -- (16,12.5) -- (16,6.25);
\node[anchor=west,font=\large] at (16.15,9.25) {"cluster theories"};

% Exchange-quiver row.
\node[align=left,anchor=west,font=\large] at (1.5,5.25)
  {exchange graphs\\of MMs,};
\node[math,anchor=west] at (5.5,5.25)
  {$\EG(\MM(R_{f_1,\cdots,f_n}))
   \cong \EG(\MM(R_{g_1,\cdots,g_n}))
   \cong \cdots \cong
   \EG(\MM(R))$};

% Stability-condition row.
\node[align=left,anchor=west,font=\large] at (1.5,3.25)
  {spaces of\\stability conditions,};
\node[math,anchor=west] at (5.5,3.25)
  {$\Stab\Db (\fl \bsL_{M_{f_1,\cdots,f_n}})
   \cong \Stab\Db (\fl \bsL_{M_{g_1,\cdots,g_n}})
   \cong \cdots \cong
   \Stab\Db (\fl \bsL_{M})$};
\node[anchor=west,font=\large] at (1.5,2) {etc.\ $\ldots$};

% The bounding box is intentionally limited to the diagram itself so that the
% resulting PDF can be inserted directly with \includegraphics[width=\linewidth].
\end{tikzpicture}
\caption{Universal singularity}
\label{ovewview univaarsal sing}
\end{figure}

\end{eg}

A similar structure also appears in the following example of a $cD_4$ singularity. Moreover, in this setting, various results are known concerning the ring that plays the role of a universal singularity. (e.g. \cite{Threefold-flops-via-matrix-factorization,The-length-classification-of-threefold-flops-via-noncommutative-algebras,Quivers-from-matrix-factorizations})

\begin{eg}[length 2 flops]

By \cite[Theorem 8.11]{Local-normal-forms-of-noncommutative-functions}, the three-dimensional cDV singularities admitting a simple $D_4$ flop have been completely classified up to analytic isomorphism. More precisely, for $n\in\mathbb{Z}_{>0}$ and $m\in\{1,\ldots,2n\}\cup\{\infty\}$, let $g_m(s):=s^m$ if $m<\infty$, and $g_\infty(s):=0$. Then every such singularity is isomorphic to

$$
R_{n,m}:=\frac{\mathbb{C}[[r,s,u,v]]}{(u^2+sv^2-r(s^{2n+1}+(r+g_m(s))^2))}.
$$

Moreover, distinct pairs $(n,m)$ give non-isomorphic singularities. In particular, the case $m=\infty$ recovers the Laufer family.

On the other hand, Curto and Morrison constructed the universal flop of length two \cite{Threefold-flops-via-matrix-factorization}, whose base singularity is given by

$$
R:=\frac{\mathbb{C}[[r,s,t,u,v,w,z]]}{(u^2-rw^2+2zvw-sv^2+(rs-z^2)t^2)}.
$$

For each pair $(n,m)$ as above, consider the elements $t_1:=z$, $t_2:=w-r-g_m(-s)$ and $t_3:=t-(-s)^n$ of $R$. Then $t_1,t_2,t_3$ form an $R$-regular sequence, and we have

$$
R/(t_1,t_2,t_3)\cong R_{n,m}.
$$

Indeed, the isomorphism follows by substituting $z=0$, $w=r+g_m(-s)$ and $t=(-s)^n$, followed by replacing $s$ with $-s$. Moreover, the quotient is a three-dimensional hypersurface. Since $R$ is a six-dimensional Cohen--Macaulay local ring, the three elements form an $R$-regular sequence.

Consequently, every three-dimensional cDV singularity admitting a simple $D_4$ flop can be obtained from the universal singularity $R$ by quotient by regular sequence.

Moreover, let $M=R\oplus N\in \MMG(R)$ be the
MM generator associated with the universal flop,
where $N$ is indecomposable of rank two.
For each $(n,m)$, put $M_{n,m}:=
(M/(t_1,t_2,t_3)M)$. By Corollary \ref{Hom quot comm} and Lemma \ref{quot indec isom}, $M_{(n,m)}$ is modifying and has two non-isomorphic indecomposable direct summands. Since every basic MM generator over $R_{n,m}$
has precisely two indecomposable summands, we see that $M_{(n,m)}\in\MMG(R_{(n,m)})$.

Thus, this setting satisfies the assumptions of Setup \ref{intro setup  reduction NCCR,tilt,2-silt}. Consequently, the "cluster theories" associated with the element of the family $\{R_{(n,m)}\}_{n,m}$ are all isomorphic. 
\end{eg}

These examples naturally lead to the following question.

\begin{ques}\label{question}
    For a family of cDV singularities, what kind of
    "universal singularity" can one consider, and
    in what sense is it natural?
\end{ques}

A {\it contracted extended Dynkin datum} is a pair $(X,\mathscr{J})$, where $X$ is an ADE extended Dynkin diagram and $\mathscr{J}$ is
a subset of the vertex set of $X$.

\begin{dfn}
Let $(X,\mathscr{J})$ be a contracted extended Dynkin datum.
    \begin{enumerate}
        \item A pair $(R,M)$ with $M\in\MM(R)$ is
        {\it of type $(X,\mathscr{J})$} if there exists
        a regular element $t\in R$ such that $\oR$ is a simple singularity of type $X$ and 
        $\End_{\oR}(\oM^{**})\cong
        e_{\mathscr{J}}\Pi_Xe_{\mathscr{J}}$.
        \item $R$ is {\it of type $(X,\mathscr{J})$}
        if there exists $M\in\MM(R)$ such that
        $(R,M)$ is of type $(X,\mathscr{J})$.
    \end{enumerate}
\end{dfn}

As one possible formulation of Question \ref{question},
we propose the following conjecture.

\begin{conj}
Fix a contracted extended Dynkin datum $(X,\mathscr{J})$.
Then there exist a complete Gorenstein normal local
domain $R_{(X,\mathscr{J})}$ and
$M_{(X,\mathscr{J})}\in\MMG(R_{(X,\mathscr{J})})$
such that the following holds:

For any isolated cDV singularity $R$ of type
$(X,\mathscr{J})$, there exist $M\in\MMG(R)$
and a regular sequence $t_1,\dots,t_m$
in $R_{(X,\mathscr{J})}$ such that $(R,M)$ is of type $(X,\mathscr{J})$,  $\overline{R_{(X,\mathscr{J})}}\cong R$ and $\overline{M_{(X,\mathscr{J})}}\cong M$. Here, the $\overline{(-)}$ denotes modulo
$(t_1,\dots,t_m)$.
\end{conj}

Universal flopping contractions, which are expected to be geometric counterparts of universal MMAs, have been constructed and studied in \cite{The-length-classification-of-threefold-flops-via-noncommutative-algebras,Centres-and-semisimple-deformations-of-contraction-algebras} and related works. Thus it seems possible to construct universal MMAs using universal flopping contractions. However, it would be interesting to give a purely algebraic construction without relying on geometry. At present, this remains a conjecture.
\appendix
\section{Ascent for some ring properties}
\begin{setup}
    In this appendix, we let $(R,\mm)$ be a Cohen--Macaulay local ring and let $0\ne t\in \mm$ be a regular element of $R$.
\end{setup}
\begin{lem}\label{appendixlemma}
   For any $\m\ne P\in\Spec R$, there exists $Q\in\Spec R$ such that $\height Q=\height P+1$ and $(t)+P\subseteq Q$.
\end{lem}
\begin{proof}
If $t \in P$, then the claim is immediate from CMness. Now assume that $t\notin P$. Then $\overline{t}\in R/P$ is non-zero. Let $\overline{Q}\in\Spec R/P$ be a minimal prime of $\overline{t}$ in $R/P$. By Krull's principal ideal theorem, $\height\overline{Q}=1$. Regard $Q$ as a prime ideal of $R$. Then $\height Q=\height P+1$ since $R$ is Cohen--Macaulay. Hence $Q$ is the desired prime ideal.
\end{proof}
\begin{prop}
    We fix $n\in\ZZ_{\ge 0}$. If $R/tR$ satisfies $(R_n)$, then $R$ also satisfies $(R_n)$.
\end{prop}
\begin{proof}
The case of $n\ge \dim R$ is well-known. Hence we may assume that $n<\dim R$. Let $P\in\Spec R$ with $\height P\le n$. Since $n<\dim R$, we have $P\ne \mm$. Therefore, by Lemma \ref{appendixlemma}, there exists $Q\in\Spec R$ such that $\height Q=\height P+1$ and $t\in Q$. By the catenary property, $\height Q/tR=\height Q -1\le n$. Thus, by assumption, $R_{Q}/tR_{Q}=(R/tR)_{Q/tR}$ is regular. Hence $R_{Q}$ is regular, and therefore the localization $R_{P}$ is regular.
\end{proof}
\begin{lem}
    If $R/tR$ is a domain, then $R$ is also a domain.  
\end{lem}
\begin{proof}
Suppose that $ab=0$ for some $a,b\in R$. We show that either $a=0$ or $b=0$. Assume that $a,b\neq 0$. By Krull intersection theorem, there exist integers $r,s\geq 0$ such that
\[a\in t^rR\setminus t^{r+1}R
\quad
b\in t^sR\setminus t^{s+1}R.\]
Thus, we can write
\[a=t^ra' \quad b=t^sb'\]
for some $a',b'\notin tR$. Then $0=ab=t^{r+s}a'b'$.
Since $t\in\reg R$, $a'b'=0$. Hence $\overline{a'}\,\overline{b'}=0$ in $R/tR$. On the other hand, since $a',b'\notin tR$, both $\overline{a'}$ and $\overline{b'}$ are nonzero. This is a contradiction.
\end{proof}
\begin{cor}\label{normal +1}
    If $R/tR$ is a normal domain, then $R$ is also a normal domain.
\end{cor}
\bibliographystyle{alpha}
\bibliography{ref}
\nocite{*}
\end{document}